\documentclass[11pt,a4paper]{article}

\usepackage{ifthen} 
\usepackage{calc} 

\usepackage[german,english]{babel} 

\usepackage[all]{xy}
\usepackage{amsmath}
\usepackage{amssymb}
\usepackage{latexsym}
\usepackage{enumerate}
\usepackage{theorem}
\usepackage{mathrsfs}
\usepackage{helvet}

\usepackage{lastpage} 
\usepackage{prettyref} 
\usepackage[pagebackref]{hyperref} 

\input xy 

\newboolean{finalversion} 
\setboolean{finalversion}{true}

\theoremheaderfont{\scshape}

\newrefformat{chap}{\hyperref[{#1}]{Chapter~\ref*{#1}}}
\newrefformat{sec}{\hyperref[{#1}]{Section~\ref*{#1}}}

\newtheorem{theorem}{Theorem}[section]
\newrefformat{thm}{\hyperref[{#1}]{Theorem~\ref*{#1}}}
\newtheorem{definition}[theorem]{Definition}
\newrefformat{def}{\hyperref[{#1}]{Definition~\ref*{#1}}}
\newtheorem{lemma}[theorem]{Lemma}
\newrefformat{lem}{\hyperref[{#1}]{Lemma~\ref*{#1}}}
\newtheorem{proposition}[theorem]{Proposition}
\newrefformat{prop}{\hyperref[{#1}]{Proposition~\ref*{#1}}}
\newtheorem{corollary}[theorem]{Corollary}
\newrefformat{cor}{\hyperref[{#1}]{Corollary~\ref*{#1}}}
{\theorembodyfont{\upshape}
\newtheorem{examplecore}[theorem]{Example}}
\newrefformat{ex}{\hyperref[{#1}]{Example~\ref*{#1}}}

\newtheorem{remark}[theorem]{Remark}
\newrefformat{rem}{\hyperref[{#1}]{Remark~\ref*{#1}}}

\newcommand{\Spec}{\ensuremath{\operatorname{Spec}}}

\newenvironment{proof}{\noindent\textsc{Proof:}}{\hspace*{\fill}
$\blacksquare$\par\vspace{.1cm}} 

\newenvironment{example}{\begin{examplecore}}{\hspace*{\fill}
$\square$\par\vspace{.1cm}\end{examplecore}}

\newcommand{\mylabel}[1]{\label{#1}\ifthenelse{\boolean{finalversion}}{
  }{\marginpar{\tiny #1}}}

\title{Valuation Theory, Riemann Varieties and the Structure of integral Preschemes}

\author{Stefan G\"unther}

\date{March,31, 2012}

\begin{document}

\maketitle
                    
\section{Abstract}
\small{In this work we show that the classical subject of general valuation theory and Zariski-Riemann varieties has a much wider scope of applications than commutative algebra and desingularization theory.\\
We construct and investigate birational projective limit objects appropriate for the study of countably many birational models at one time.\\
 We use nonseparated Riemann varieties to investigate the birational structure of integral preschemes satisfying the existence condition of the valuative criterion of properness.}\\

\tableofcontents

\section{Introduction}
Basically, the study of varieties up to birational equivalence is the study of algebraic function fields. Valuation theory is one means to do this. The study of  discrete algebraic rank one valuations can be considered as the study of Cartier divisors on some birational model $X$ of $K=K(X)$\,. The study of linear series on varieties is actually equivalent to looking for rational functions having a value bounded below at some finite set of discrete algebraic rank one valuations.\\
We will show in this work that the classical subject of general valuation theory and Zariski-Riemann varieties has  a much wider scope of applications than commutative algebra and desingularization theory.\\
Even if one is interested in particular birational models, for instance the so called minimal models, the process of finding these involves apriori an infinite set of birational models and the definitions of the standard classes of singularities considered in the Log Minimal Model program involve considering divisors on all sufficiently high birational models of a given variety $X$. So it is desirable to have some geometric object that "lies above", i.e., dominates all birational models under consideration. There is a well defined object, depending only on the function field, called the Zariski-Riemann variety $R(K/k)$ of $K$, which is a locally ringed space, all of whose local rings are the valuation rings of the function field. It dominates each model of the function field. So whenever arguments involve considering infinitely many birational models of a variety $X$ at one time such as an apriori infinite log flip sequence, one could try to implement  these arguments directly on the Riemann variety or a countable birational limit object in the category of locally ringed spaces and then use the geometric properties of $R(X/k)$\, such as ,e.g, quasicompactness  to carry on. For instance, it is easy to show using quasicompactness that a $b$-Cartier-divisor on $K(X)$\, is simply a Cartier divisor on the locally ringed space $R(K(X))$.\\
Preschemes are considered by the broad majority of algebraic geometers as "pathological " or ill behaved phenomena of Grothendieck algebraic geometry. Most of them know just the example of the projective line with one point doubled and bear in mind that one cannot say anything useful about them in full generality.\\
In this work, we will show that this is not the case.\\
First, nonseparated schemes and nonseparated algebraic spaces occur naturally as moduli spaces of nonseparated moduli functors, such as the moduli functor of all projective nonsingular varieties of a fixed dimension.\\
Secondly, in dimension larger than one, by patching together arbitrary complete models of a given function field along the open loci where they are isomorphic gives many interesting examples of integral preschemes that are not of the standard kind obtained by doubling subvarieties of a given complete variety.\\
 After fixing the necessary notations we start by giving a short review of classical valuation theory such as can be found in \cite{ZaSam}[Volume 2, Chapter VI].\\ 
In section four, we show that projective limits of preschemes exist in the category of locally ringed spaces. We will study geometric properties of Riemann varieties and more generally of projective limit objects  of preschemes in the category of locally ringed spaces such as the structure of closed subsets.\\
In \ref{thm:Chevalley} we prove that the projective limit $\mbox{projlim}_{X_j\longrightarrow X_i}X_j$\, over a system of quasicompact preschemes is quasicompact with respect to either the projective limit of the Zariski or constructible topologies on the preschemes $X_i$\,. This is a generalization of the classical theorem of Chevalley saying that the Riemann variety of a function field is quasicompact and of a result of F.V. Kuhlmann saying that the Riemann variety is still quasicompact with respect to the constructible topology (see \cite{Kuhl}[Appendix, Theorem 36,p.21]). On the contrary, we show that the space of discrete algebraic rank one valuations of a function field of transcendence degree larger than one is not quasicompact.\\
In \ref{thm:A310} we prove that the projection maps $p_i: \mbox{projlim}_{X_j\longrightarrow X_i}X_j\longrightarrow X_i$\, are closed if one assumes that the transition morphisms $p_{ij}: X_j\longrightarrow X_i$\, are all proper. This may be rephrased by saying that "the projective limit of proper morphisms is again proper". \\
In section five we will use nonseparated Riemann varieties to study phenomena of integral preschemes satisfying the existence condition of the valuative criterion of properness, which carry a rich birational geometric structure.\\
We hope that our investigations will be helpful to answer the following question:\\
Given an integral prescheme $X$\, of finite type /$k$ that satisfies the existence condition of the valuative criterion of properness, under which circumstances can there always be found an open complete subscheme $X'\subset X$?\\
A positive answer to this question is undoubtedly helpful in order to compactify families of algebraic objects if these form a coarsely representable moduli functor that satisfies the existence condition of the valuative criterion of properness.\\
 We give a criterion for an integral prescheme to satisfy the existence condition of the valuative criterion of properness and show that for an integral prescheme $X$ of finite type to be complete it suffices that each discrete algebraic rank one valuation has exactly one center on $X$\,.\\
We study integral preschemes $X$\, by means of the basic diagram in the category of locally ringed spaces $$X\stackrel{p_X}\longleftarrow R(X)\stackrel{q_X}\longrightarrow R(K(X)),$$ where $R(X)$ is the Riemann variety of the integral prescheme $X$, $p_X$\, is the canonical projection, $R(K(X))$\, is the Riemann variety of the function field and $q_X$\, is the canonical local isomorphism. We introduce and study  the subsets $S_n(X)\subset R(K(X))$\, of all valuations having under the map $q_X$\, at least $n$ preimages in $R(X)$\,. Our main result is \ref{thm:T2068} (with most of the work being done in \ref{thm:A311}), asserting that for an integral prescheme of finite type satisfying the existence condition of the valuative criterion of properness, the subsets $S_n(X)$\, are Zariski closed in $R(K(X))$\, for all $n\in \mathbb N$\,. By closedness of the map $p_X$\,, the sets $T_n(X):=p_X(q_X^{-1}(S_n(X)))$\, are Zariski closed subsets of $X$. These sets consist of all $x\in X$\, for which there exists a valuation $\nu$\, having beside $x$ at least  $n-1$ other centers on $X$. \\
We introduce on each integral prescheme $X$ the canonical stratification $(X_n)_{n\in \mathbb N}$\, defined by $X_n:=T_n(X)\backslash T_{n+1}(X)$\,. Each $X_n$\, is by definition a constructible set. We show that for $X$ of finite type this is a finite stratification.\\
We hope that our main result \ref{thm:T2068} has applications to the study of nonseparated moduli spaces and also makes them worth to be considered. We have the following in mind. The moduli functor $\mathcal M$\, of all projective  reduced but not nessessariliy irreducible varieties  of general type with semi-log canonical singularities that are minimal models, i.e., with $K_X$\, big and nef, should satisfy the existence condition of the valuative criterion of properness. Thus its coarse moduli space  $M$\, should be a (nonseparated) prescheme or nonseparated algebraic space of finite type satisfying the existence condition of the valuative criterion of properness. Its canonical stratification $(M_n)_{n\in \mathbb N}$\, should have the following meaning. If $m\in M_n$\,, there exist points $m_1=m,m_2,...,m_n\in M_n$\, all being centers of the same valuation. These points represent complete varieties $X_1=X,X_2,...,X_n$\,. These should be all birationally equivalent to $X$\, and should consist of all minimal models of $X$. Thus $n$\, could be interpreted as the number of minimal models of a given complete variety \,$X$.
\section{Notations and Conventions} 
Throughout this work $k$ denotes a base field of characteristic zero.\\
$K$ denotes a function field over $k$.\\
If $\nu$\, is a Krull valuation of $K$ we denote by \\
$A_{\nu}$\, the corresponding valuation ring,\\
 $\mathfrak{m}_{\nu}$\, its maximal ideal and \\
 $k_{\nu}$\ or $\kappa(\nu)$\, the residue field.\\
Krull valuations of a field we will always denote by greek letters $\nu,\mu,\omega,...$\, whereas for scheme points we will use latin as well as greek letters $x,y,x,\eta,\xi,...$\,.\\
By a (birational) model for $K$ we understand a complete integral normal variety $X$ with $K(X)=K$, $K(X)$ denoting the field of rational functions on $X$.\\
$R(K/k)$\, or simply $R(K)$\, denotes the Zariski-Riemann variety of the field extension $K/k$\,.\\
For a finitely generated $k$-subalgebra $A$ of $K$\,,  $R^K\Spec A$\,  denotes the subset of all Krull valuations $\nu\in R(K/k)$\, with $A\subset A_{\nu}$\,. If the function field in question is fixed and clear from the context or is the quotient field of $A$, we will use the simpler notation $R\Spec A$\,.\\
For $X/k$ any integral prescheme $R(X/k)$\, or simply $R(X)$\, denotes throughout the Riemann variety  $X/k$\,.\\
For a field extension $K/k$\,, $R_{Ab}(K/k)$\, denotes the subset of all Abhyankar places and $R^{m,n}(K/k)$\, denotes the subset of all places of rational rank $m$ and rank $n$. The set of discrete rank one algebraic places we simply denote by $R^1(K/k)$\,.\\   
If $X$ is an integral scheme and $\nu\in R(K(X))$\, we denote by $c_X(\nu)$\, the unique center of $\nu$\, on $X$, if it exists, and that is then the unique scheme point $\eta$\, of $X$ such that the valuation ring $A_{\nu}$\, dominates $\mathcal O_{X,\eta}$\,. We put $C_X(\nu)=\overline{{c_X(\nu)}}$\,. By saying that $\nu$\, has center on $\Spec A$\, we mean that $A\subset A_{\nu}$\,. If this is the case, then the scheme point corresponding to the prime ideal $\mathfrak{m}_{\nu}\cap A$\, corresponds to the point that is the center of $\nu$\, on $\Spec A$\,.\\
If $A$\, is an integral $k$-algebra inside some field $K$, $A^n$\, denotes the normalization of $A$ inside its quotient field and $A^{n,K}$\, the normalization of $A$ inside $K$.\\
If $(X,\mathcal O_X)$\, is a locally ringed space and $K$ is an abelian group, we denote by $\underline{K}_X$\, or simply by $\underline{K}$\, the constant sheaf associated to the presheaf that takes as value the abelian group  $K$ on every open subset of $X$.\\
If $\omega$\, is a rational top differential form of the function field $K(X)$\,, and $X$\, an integral normal prescheme, then $K_X^{\omega}$\, denotes the divisor of zeroes and poles of the rational section corresponding to $\omega$\,of the canonical reflexive sheaf of $X$.\\
By the existence condition of the valuative criterion of properness we mean the following property of an integral prescheme $X/k$:
Given an arbitrary field $L$ and a valuation $\nu\in R(L/k)$\,, for each morphism $$f:\Spec L\longrightarrow X/\Spec k$$ there exists at least one morphism $$g: \Spec A_{\nu}\longrightarrow X/\Spec k$$ such that $f=g\circ j$\, where $j: \Spec L\longrightarrow \Spec A_{\nu}$\, is the canonical morphism corresponding to the inclusion $j^{\sharp}: A_{\nu}\longrightarrow L$\,.
 \section{Review of valuation theory}
  We will use in this paper only well known facts about valuation theory of function fields as can be found in \cite{ZaSam}[Volume II, chapter VI], which we will review for the convenience of the reader. For a good readable account of the basics of valuation theory, see \cite{Val}. We will only be concerned with Krull valuations of a field $K$, or, more generally of the relative situation of a field extension $K/k$\,, $k$ being a base field of characteristic zero.\\
   We recall that a Krull valuation is a homomorphism of abelian groups 
   $$\nu:K^*=K\backslash \{0\}\longrightarrow \Gamma,$$ $K^*$\, with the multiplication and $\Gamma$\, being totally ordered such that for all $f,g\in K$\, one has $$\nu(f+g)\geq \rm{min}(\nu(f),\nu(g)),$$ where $\nu(0)$\, is considered to be $\infty$\,.\\
    The set $A_{\nu}\subset K$\, consisting of all $f\in K$\, such that $\nu(f)\geq 0$\, (including $f=0$\,) is easily seen to be a ring, not necessarily noetherian, called the valuation ring of $\nu$\, and the subset $\mathfrak{m}_{\nu}\subset A_{\nu}$\, consisting of elements $f$\, with $\nu(f)>0$\, is the maximal ideal.
     The field $k_{\nu}:=A_{\nu}/\mathfrak{m}_{\nu}$\, is the usual residue field and the image of the homomorphism $\nu$\, in $\Gamma$\, is denoted by $\Gamma_{\nu}$\, and is called the value group of $\nu$\,.\\
      Two valuations are said to be equivalent, if the associated valuation rings are equal.\\
    Valuation rings inside a field can be characterized as being maximal subrings of $K$ with respect to the partial order of domination between local rings, i.e., the relation $(A,m)<(B,n)$\, iff there is a local homomorphism $(A,m)\longrightarrow (B,n)$\, being the identity on generic points (see \cite{Val}[chapter 1, p.1]) .\\
     From this maximality property a valuation ring, $A_{\nu}$\, has the property that for each $f\in K^*$\, at least one of the two elements $f$ and $f^{-1}$\, is in $A_{\nu}$\,. It follows, that the abelian group $K^*/A_{\nu}^*$\, is totally ordered with respect to the relation $$[f]<[g]\,\, \mbox{if}\,\,[gf^{-1}]\in A_{\nu}.$$
     There is always  a canonical isomorphism of totally ordered abelian groups $$(*):K^*/A_{\nu}^*\cong \Gamma_{\nu}\,\, \mbox{sending}\,\, [f]\,\, \mbox{to}\,\, \nu(f).$$
      If $L\subset K$\, is any subfield, we may restrict $\nu$\, to $L$ and obtain a valuation of $L$, denoted by $\nu\mid_L$\,.\\
      It is known (see \cite{ZaSam}[Vol.I, chapter 7, Theorem 11] that each valuation $\nu$ of a field $K$ can be extended to a valuation $\mu$\, of any field extension $K\subset L$\,, algebraic or not, such that $\mu\mid_K=\nu$\,.\\
A valuation $\nu$\, of $K$\, will be called  over a subfield $k\subset K$\,, if $\nu\mid_k$\, is the trivial valuation, i.e., $\nu(c)=0$\, for all $c\in k^*$\,. $k$ will usually be the fixed base field.  If we assume the restriction to be trivial, we will simply write $\nu\in R(K/k)$\,.\\ 
We denote the dimension of the $\mathbb Q$-vector space $\Gamma_{\nu}\otimes_{\mathbb Z}\mathbb Q$\, by $rr(\nu)$\, and call it the rational rank of the valuation $\nu$. The Krull dimension of the valuation ring $A_{\nu}$\, will be called the rank of $\nu$\,, denoted by $r(\nu)$\,. We always have an inequality $r(\nu)\leq rr(\nu)$\,. The transcendence degree $\mbox{trdeg}(k_{\nu}/k)$\, will be called the dimension of $\nu$\, denoted by $\dim(\nu)$\,.\\
A segment $\Delta$\, of a totally ordered abelian group $\Gamma$\, is a symmetric subset such that 
$$ \alpha,\beta\in \Delta, \alpha <\gamma< \beta \vee \beta<\gamma<\alpha \Rightarrow \gamma\in \Delta.$$
By \cite{Val}[chapter 1.2, Proposition 1.6, p.4] there is a one-to-one inclusion reversion correspondence between segments of $\Gamma_{\nu}$\, and ideals of $A_{\nu}$\,. If $\mathfrak{a}\subset A_{\nu}$\, is an ideal, then the associated segment $\Delta(\mathfrak{a})$\, is defined by $$\Delta(\mathfrak{a}):=\{\gamma\in \Gamma\mid \gamma<\nu(a)\forall a\in \mathfrak{a}\vee \gamma >-\nu(a)\forall a\in \mathfrak{a}\}.$$
If $\Delta\subset \Gamma$\, is a segment, then the associated ideal $\mathfrak{a}(\Delta)$\, is defined by $$\mathfrak{a}(\Delta):=\{a\in A_{\nu}\mid \nu(a)> \delta\forall \delta\in \Delta\}.$$
We recall that the order rank of a totally ordered abelian group $\Gamma$\, is the length of the maximal chain of convex subgroups
$$0\subsetneq \Gamma_1\subsetneq \Gamma_2\subsetneq ...\subsetneq \Gamma_n=\Gamma.$$
By \cite{Val}[chapter 1.2,Corollary Theorem 1.7,p.5] the Krull dimension of $A_{\nu}$\, equals the order rank of $\Gamma_{\nu}$\,. 
If $K$ happens to be of finite transcendence degree over $k$, all the above defined numbers are finite and there is always the fundamental Abhyankar inequality $$rr(\nu)+\dim(\nu)\leq \rm{trdeg}(K/k).$$ A place $\nu$\, for which equality holds is called an Abhyankar place. If $K$ happens to be finitely generated over $k$, then for an Abhyankar place $\nu$\,, the value group $\Gamma_{\nu}$\, is known to be a finite (free) $\mathbb Z$-module and the residue field $k_{\nu}:=A_{\nu}/\mathfrak{m}_{\nu}$\, is known to be finitely generated over $k$. Both properties may fail in general for arbitrary places of a function field $K/k$\, and such places always exist if $\rm{trdeg}(K/k)\geq 2$\,.\\
It is known that every totally ordered abelian group $\Gamma$\, with finite rational rank can be order imbedded into some $\mathbb R^n$\,with the lexicographic product order of the usual order on the factors $\mathbb R$\,. We always have $n\geq r(\Gamma)$\, and it is known that the smallest $n\in \mathbb N$\, for which there exist such an embedding is equal to the  rank, i.e., if $\Gamma=\Gamma_{\nu}$, then $n$ equals the Krull dimension of $A_{\nu}$\,, hence the rank of $\nu$\,. \\
We will consider in this paper exclusively fields $K$ of finite transcendence degree over a base field $k$\,. We denote by $R^{n,l}(K/k)$\, the subset of all valuations $\nu$\, of $R(K/k)$\, with $rr(\nu)=n$\, and $r(\nu)=l$\,. We denote by $R_{Ab}^{n,l}(K/k)\subset R^{n,l}(K/k)$\, the subset of all such Abhyankar places. \\
The space of all algebraic discrete rank one valuations ($R^{1,1}_{Ab}(K/k)$) we will simply denote by $R^1(K/k)$\,.
\\
If $A_{\nu}\subset K$\, is a valuation ring and $A_{\nu_1}\subset k_{\nu}$\, is another valuation ring, we can form the composed ring $$A_{\nu\circ \nu_1}:=p^{-1}(A_{\nu_1})$$
 where $p:A_{\nu}\longrightarrow k_{\nu}$\, is the canonical residue map. As the notation suggests, $A_{\nu\circ \nu_1}$\, is again a valuation ring of $K$ and the corresponding valuation is called the composition of $\nu$\, with $\nu_1$\,, denoted by $\nu\circ \nu_1$\,. Because $$A_{\nu\circ \nu_1}\backslash\mathfrak{m}_{\nu\circ n_1}\subset A_{\nu}\backslash\mathfrak{m}_{\nu},$$
  there is a canonical surjection $$K^*/A_{\nu\circ \nu_1}^*\longrightarrow K^*/A_{\nu}^*,$$
   whose kernel is canonically isomorphic to $p^{-1}(k_{\nu_1}^*)/p^{-1}(A_{\nu_1}^*)$\, which is again canonically isomorphic to $k^*_{\nu_1}/A_{\nu_1}^*$\,. \\
   Using the canonical isomorphisms $(*)$\,, we get a canonical exact sequence of ordered abelian groups $$(**)0\longrightarrow \Gamma_{\nu_1}\longrightarrow \Gamma_{\nu\circ \nu_1}\longrightarrow \Gamma_{\nu}\longrightarrow 0,$$ such that $\Gamma_{\nu_1}\longrightarrow \Gamma_{\nu\circ \nu_1}$\, is an order inclusion and $\Gamma_{\nu_1}$\, is a  convex subgroup of $\Gamma_{\nu\circ \nu_1}$\, in the sense that if 
   $$\alpha,\beta\in \Gamma_{\nu_1}, \gamma\in \Gamma_{\nu\circ \nu_1}\,\, \mbox{and}\,\, \alpha<\gamma<\beta,\,\, \mbox{then}\,\, \mbox{also}\,\, \gamma\in \Gamma_{\nu_1}.$$ The order on $\Gamma_{\nu}$\, is isomorphic under these isomorphisms to the induced quotient order on $\Gamma_{\nu\circ \nu_1}/\Gamma_{\nu_1}$\,.\\ Every totally ordered abelian group must necessarily be a torsion free $\mathbb Z$-module, hence every such sequence possesses a (noncanonical) order preserving splitting $$\Gamma_{\nu}\longrightarrow \Gamma_{\nu\circ \nu_1}.$$
    If $\Gamma_{\nu}$\, happens to possess a finite $\mathbb Z$-basis $\gamma_i, i=1,...,n$\,, then any choice of elements $t_1,...,t_n$\, with $\nu(t_i)=\gamma_i$\,, which we sloppily call a $\mathbb Z$-basis of $\nu$\,, induces a splitting of $(**)$\,. Namely, writing $\gamma\in \Gamma_{\nu}$\, as $[f]\in K^*/A_{\nu}^*$\, under the above isomorphism $(*)$\,, writing $$\nu(f)=\sum_in_i\nu(t_i), n_i\in \mathbb Z,$$
     we first send $[f]$\, to $$(\nu_1(p(f\cdot\prod t_i^{-n_i}))\in \Gamma_{\nu_1}$$ (observe $f\cdot \prod_it_i^{-n_i}\in A_{\nu}^*$\,). This is indeed a group homomorphism. We define the isomorphism $$\phi(\underline{t}):\Gamma_{\nu\circ \nu_1}\cong \Gamma_{\nu}\oplus \Gamma_{\nu_1}$$ by sending $[f]$\, to $(\nu(f), \nu_1(p(f\cdot \prod_it_i^{-n_i})))$\,.\\
 \paragraph{Examples of zero dimensional valuations}
 \begin{example}\mylabel{ex:E60}(Monomial valuations)\\
Suppose $(A,\mathfrak{m},k)$\, is a regular local ring essentially of finite type and the residue field $k$ is algebraically closed ( and of characteristic zero). Then $(A,m,k)$\, is a smooth local $k$-algebra and there are elements  $r_1,...,r_n\in A$\, such that the induced morphism $\phi:k[x_1,...,x_n]_{(x_1,...,x_n)}\longrightarrow A, x_i\mapsto r_i, i=1,...,n$\, is etale (see \cite{SGA}[SGA I, Expose II, Definition 1.1, p.29]). $r_1,r_2,...,r_n$\, is then necessarily  a regular sequence generating the maximal ideal $\mathfrak{m}$\,. Let $\widehat{A}$\, be the completion of $A$ along the maximal ideal. As $\phi$\, is etale, we have 
$$\widehat{A}\cong \widehat{k[x_1,x_2,...x_n]_{(x_1,x_2,...,x_n)}}\cong k[[x_1,x_2,...,x_n]].$$
We have thus an inclusion $A\hookrightarrow k[[x_1,...,x_n]]$\, where $r_i$\, maps to $x_i$\, for $i=1,...,n$\,.  Now choose positive real numbers $c_i>0$\, and if 
$$p(r_1,...,r_n)=\sum_I\underline{r}^I$$ is a homogenous polynomial define $$\nu(p)=\mbox{min}_I(\prod_{j=1}^ni_j\cdot c_j)>0.$$
 Now if $p$\, is a power series in $\widehat{A}$\,, define $$\nu(p)=\mbox{inf}_n\nu(p_n),$$
 where $p_n$\, is the homogeneous part of $p$ of degree $n$. If $c=\mbox{min}_{i=1}^n c_i$\,, then \\
$\nu(p_n)\geq nc$\,, which tends to infinity. Thus the above infemum is for every power series actually a minimum. Thus we have $\nu(p)\in \mathbb Z(c_1,...,c_n)\subset \mathbb R$\,. It is clear that $\nu(p+q)\geq \mbox{min}(\nu(p),\nu(q))$\, as power series are added term by term. Similarly, there is a monomial $m_p$\, occurring in the power series of $p$ and a a monomial $m_q$\, in the power series of $q$ such that $\nu(p)=\nu(m_p)$\, and $\nu(q)=\nu(m_q)$\,. Then the monomial $m_pm_q$\, occurs in the power series of $pq$\, and is there the monomial with minimal value. Thus $$\nu(pq)=\nu(m_pm_q)=\nu(m_p)+\nu(m_q)=\nu(p)+\nu(q).$$
 We thus can extend $\nu$\, to $k((r_1,...,r_n))$\, by $\nu(\frac{p}{q})=\nu(p)-\nu(q)$\, to get a valuation of $k((r_1,...,r_n))$\,. Restricting to $K(A)\subset k((r_1,...r_n))$\, we get a valuation of $K(A)$\, with center above $A$\,. It has finitely generated value group. If $c_1,...,c_n$\, are linearly independent over $\mathbb Q$\, then we get a zero dimensional Abhyankar place of rank one and rational rank $n$.
\end{example}
 \begin{example}\mylabel{ex:E61} (Flag valuations)
These are the nonnoetherian valuations that are most commonly known to nonspecialists.\\
Let $K/k$\, be a function field and $X$ be a model of $K$. Consider a flag of subvarieties
$$X=D_0\supsetneq D_1\supsetneq ...\supsetneq D_n\supsetneq 0$$
such that
\begin{enumerate}[i]
\item each $D_i, i=0,...,n$\, is integral.
\item The local ring $\mathcal O_{D_{i},\eta_{D_{i+1}}}$\, is one dimensional and regular.
\end{enumerate}
It follows that $\mathcal O_{D_i,\eta_{D_{i+1}}}$\, is a discrete algebraic rank one valuation  ring with corresponding valuation $\nu_i$\, in the function field $k(D_i)$ for all $i=1,....,n$\,.\,
We can form the successive composition $\nu:=\nu_0\circ \nu_1\circ ...\circ \nu_{n-1}$\, and get a valuation of the function field $K(X)$\,. Since under composition the rational rank and the rank are additive  (see \cite{Val}[chapter 1.2, Corollary to Proposition 1.11, p.9]), we have $rr(\nu)=n$\, and $r(\nu)=n$\, hence we have got an Abhyankar place in $R^{n,n}_{Ab}(K(X))$\,. The value group $\Gamma_{\nu}$\, is isomorphic to $\mathbb Z^n$\, but in a noncanonical way! Such an isomorphism depends on the choices of local parameters for the discrete algebraic valuation rings in the function fields $k(D_i)$\,.
\end{example}
 \begin{example}\mylabel{ex:E62}(Mixed valuations) In example (2) we might as well take an incomplete flag 
$$D_0=X\supsetneq D_1\supsetneq ...\supsetneq D_k$$
 in order to construct an Abhyankar place $\nu\in R^{k,k}_{Ab}(K(X))$\,. It has dimension $\dim(\nu)=\mbox{trdeg}(\kappa(\nu))=n-k$\,. We can take an arbitrary monomial valuation $\nu_0$\, of the function field $\kappa(\nu)$\, (example (1)). and form the composed valuation $\mu:=\nu\circ \nu_0$\,.\\ By the additivity of rank and rational rank (see \cite{Val}[chapter 1.2, Corollary to Proposition 1.11, p.9]), we have got an Abhyankar place in $R^{n,k+1}_{Ab}(K(X))$\,.
\end{example}
\section{Riemann varieties and projective limit objects in the category of preschemes}
\paragraph{Projective limit objects in the category of preschemes}
It is known that projective limits in the category of preschemes do not exist in general. They exist in the case of affine transition morphisms. In this section, we will consider projective limit objects of preschemes with proper transition morphisms. We show that they exist in the broader category of locally ringed spaces. Our main motivation to consider them is the description of Zariski-Riemann varieties as birational projective limit objects. Our main results are that projective limits of quasicompact preschemes are again quasicompact and that projective limits of proper morphisms are universally closed.\\
For each directed set $I$ and each collection 
$$X_i,i\in I, f_{ij}: X_j\longrightarrow X_i\,\, \mbox{for}\,\, i<j$$
 with $X_i$\, a prescheme and $f_{ij}$\, a proper morphism, we want to construct  the projective limit object $$\widehat{X}:=\mbox{projlim}_{i\in I}X_i$$ in the category of locally ringed spaces.\\
 As a point set, we define $$\mid \widehat{X}\mid:=\mbox{projlim}_{i\in I}\mid X_i\mid$$ and  we put on $\mid \widehat{X}\mid$\, the coarsest topology such that all the natural projections 
 $$p_{X_i}=p_i: \mid\widehat{X}\mid\longrightarrow \mid X_i\mid$$ become continuous. I.e., one declares the set 
$$\left\{\bigcap_{i\in \mathcal F}p_{X_{i}}^{-1}(U_{i})\mid \mathcal F\subset I\,\,\mbox{finite}\,\,, U_{i}\subset X_{i}\,\, \mbox{Zariski}\,\, \mbox{open}\right\} $$ to be  a basis of open subsets. \\
 This is called the Zariski topology on $\widehat{X}$\,.\\
  On $\widehat{X}$\, one has the directed system of sheaves of  rings $p_{X_i}^{-1}\mathcal O_{X_i}$,  and homomorphisms of sheaves of rings 
  $$\phi_{ij}:=p_{X_j}^{-1}(f_{ij}^{\sharp}): p_{X_i}^{-1}(\mathcal O_{X_i})\longrightarrow p_{X_j}^{-1}(\mathcal O_{X_j})$$ for each proper birational morphism $$f_{ij}: X_j\longrightarrow X_i.$$ 
   We define the structure sheaf $\mathcal O_{\widehat{X}}$\, to be the direct limit $\lim_{i\in I}p_{X_i}^{-1}\mathcal O_{X_i}$\, of the inductive system $((p_{X_i}^{-1}(\mathcal O_{X_i}))_{i\in I}, \phi_{ij})$\,. \\
   For each point $\widehat{x}\in \widehat{X}$\, that is given by an indexed set of points $x_i\in X_i$\, $f_{ij}(x_j)=x_i$\, the stalk $\mathcal O_{\widehat{X},\widehat{x}}$\, is the inductive limit of the local rings $\mathcal O_{X_i,x_i}$\, under the local homomorphisms $$p_j^{-1}(f_{ij}^{\sharp}): \mathcal O_{X_i,x_i}\longrightarrow \mathcal O_{X_j,x_j}$$  and so itself local. Therefore, the space $(\widehat{X},\mathcal O_{\widehat{X}})$ is a locally ringed space.\\
   \begin{lemma}\mylabel{lem:L730} Let $((X_i)_{i\in I},f_{ij})$\, be as above an indexed system of preschemes $X_i$\, and proper morphisms $f_{ij}$\,. Then the above defined locally ringed space $(\widehat{X},\mathcal O_{\widehat{X}})$\, is the projective limit of this system in the category of locally ringed spaces.
   \end{lemma}
   \begin{proof}
   The proof is straight forward. Let $(Y,\mathcal O_Y)$\, be any locally ringed space and $g_i:(Y,\mathcal O_Y)\longrightarrow (X_i,\mathcal O_{X_i})$\, morphisms of locally ringed spaces commuting with the $f_{ij}$:$$g_i=f_{ij}\circ g_j\,\,\forall i,j\in I, i<j.$$
 As point sets, we get a unique map $g:\mid Y\mid \longrightarrow \mid\widehat{X}\mid$\,such that $g_i=p_i\circ g$\, since $\mid\widehat{X}\mid$\, is the projective limit in the category of sets. Let for some finite set $\mathcal F\subset I$\, $$\bigcap_{i\in \mathcal F}{p_i^{-1}(U_i)}, U_i\subset X_i\,\,\mbox{Zariski}\,\,\mbox{open}$$ be a basic open subset in $\widehat{X}$\,. Since $g_i$\, is continuous, $g_i^{-1}(U_i)=g^{-1}p_i^{-1}(U_i)$\, is open in $Y$. Hence $$ \bigcap_{i\in \mathcal F} g^{-1}(p_i^{-1}(U_i))=g^{-1}(\bigcap_{i\in \mathcal F}p_i^{-1}(U_i))$$ is open in $Y$. Thus $g$ is continuous.\\
   We have homomorphisms of sheaves $g_i^{\sharp}:g_i^{-1}\mathcal O_{X_i}\longrightarrow \mathcal O_Y$\,, commuting with the homomorphisms  $$g_j^{-1}(f_{ij}^{\sharp}): g_i^{-1}\mathcal O_{X_i}\longrightarrow g_{j}^{-1}\mathcal O_{X_j}.$$ Hence we get a unique homomorphism of sheaves
   $$ g^{\sharp}:\lim_{i\in I}g_i^{-1}\mathcal O_{X_i}=g^{-1}\lim_{i\in I}p_i^{-1}\mathcal O_{X_i}=g^{-1}\mathcal O_{\widehat{X}}\longrightarrow \mathcal O_Y$$ such that $g_i^{\sharp}=p_i^{\sharp}\circ g^{\sharp}$\,.\\
Thus $(g,g^{\sharp})$\, is a morphism of locally ringed spaces.   
   \end{proof}
   Of course, the main example we have in mind is where $I$ is the set of all proper  birational models of a function field $K/k$\, and $i<j$\, iff there is a proper birational morphism $f_{ij}: X_j\longrightarrow X_i$\,. In this case $\widehat{X}$\, is called the Zariski-Riemann variety of the function field $K/k$\,, to be introduced in the next paragraph.\\
 \begin{lemma}\mylabel{lem:L620} Let $$((X_i)_{i\in I}, f_{ij})\,\, \mbox{and}\,\, ((Y_i)_{i\in I}, g_{ij})$$ be two projective systems of preschemes indexed by the same directed set $I$. Let $(h_i: X_i\longrightarrow Y_i)_{i\in I}$\, be a collection of morphisms of preschemes such that for $i<j: h_i\circ f_{ij}=g_{ij}\circ h_j$\,.\\
 Then the induced map of the projective limits 
 $$h:\mbox{projlim}_{i\in I}X_i\longrightarrow \mbox{projlim}_{i\in I}Y_i$$
  is a morphism of locally ringed spaces.
  \end{lemma}
  \begin{proof} The projections $p_l:\mbox{projlim}_{i\in I}X_i\longrightarrow X_l$\, are morphisms of locally ringed spaces for all $l\in I$\, as $p_l$ is by definition of the projective limit topology continuous and there is a canonical map $$p_l^{\sharp}: p_l^{-1}\mathcal O_{X_l}\longrightarrow \lim_{i\in I}p_i^{-1}\mathcal O_{X_i}=\mathcal O_{\mbox{projlim}_{i\in I}X_i}.$$
  Then, $ h_l\circ p_l: \mbox{projlim}_{i\in I}X_i\longrightarrow Y_l$\, are also morphisms of locally ringed spaces for all $l\in I$\,. By \ref{lem:L730}, $\mbox{projlim}_{i\in I}Y_i$\, is a projective limit object in the category of locally ringed spaces, hence there is a unique morphism of locally ringed spaces $$(h,h^{\sharp}): \mbox{projlim}_{i\in I}X_i\longrightarrow \mbox{projlim}_{i\in I}Y_i$$ satisfying $q_l\circ h= h_l\circ p_l, \forall l\in I$\, where $q_l:\mbox{projlim}_{i\in I}Y_i\longrightarrow Y_l$\, is the projection.
  \end{proof} 
 
 We prove the following generalization of a theorem of Chevalley(\cite{Val}[Theorem 2.7]), saying that the Riemann variety of a function field is quasicompact. We will give a new proof of this classical theorem using the projective limit representation of the Riemann variety being introduced in the next section and the associated constructible topology. \\
Recall that the constructible topology on a scheme $X$ is given by the basis of constructible-open sets 
\begin{align*}\{A\cap U\mid A\subset X\,\,\mbox{Zariski}\,\,\mbox{closed}\,\,\mbox{and}\,\,
U\subset X\,\,\mbox{Zariski}\,\,\mbox{open}\}
\end{align*}
We need the following 
\begin{lemma}\mylabel{lem:L23} For each quasicompact scheme $X$, $\mid X\mid$\, with the constructible topology is a Hausdorff quasicompact space.
\end{lemma}
\begin{proof} Let $\eta,\xi\in \mid X\mid$\,.If $\eta\in \overline{\{\xi\}}$\,, then $\overline{\{\xi\}}\backslash\overline{\{\eta\}} $ is Zariski open in $\overline{\{\xi\}}$\, and there is a Zariski open subset $U\subset X$\, such that $$(\overline{\{\eta\}}) \cap(\overline{\{\xi\}}\cap U)=\emptyset.$$
$\overline{\{\eta\}}$\, is a constructible neighbourhood of $\eta$\, and $\overline{\{\xi\}}\cap U$\, is a constructible neighbourhood of $\xi$\,. If $\eta\notin \overline{\{\xi\}}$\,, then  $$(X\backslash\overline{\{\xi\}})\cap \overline{\{\xi\}}=\emptyset $$ and $(X\backslash\overline{\{\xi\}})$\, is a constructible neighbourhood of $\eta$\, and $\overline{\{\xi\}}$\, is a constructible neighbourhood of $\xi$\,. Hence $\mid X\mid$\, is Hausdorff. \\
That $X$ is quasicompact we prove by induction on the dimension.\\
 For $\dim(X)=0$\,, $\mid X\mid $ with the constructible topology is a finite discrete set of points and this is a quasicompact space.\\
 Let $\dim(X)=n+1$\, and let $X=\bigcup_{i\in I}U_i$\, be a constructible open covering. Let $\eta_i,i=1,...,d$\, be the generic points of the irreducible components of $X$\,. There are indices $i_1,...,i_d$\, such that the generic points $\eta_j\in U_{i_j}, j=1,...,d$\,.\\ $U_{i_j}$\, must contain a Zariski open subset $V_j, j=1,...,d$. $X\backslash\bigcup_{j=1}^dV_j$\, is again a quasicompact scheme of dimension at most $n$, hence by the induction hypothesis, we have that for the constructible covering
 $$X\backslash(\bigcup_{j=1}^dV_j) = \bigcup_{i\in I}(U_i\cap (X\backslash\bigcup_{j=1}^dV_j))$$
 
  there is a finite subcovering $$X\backslash(\bigcup_{j=1}^dV_j)=\bigcup_{j=d+1}^n(U_{i_j}\cap (X\backslash\bigcup_{j=1}^dV_j)).$$
   Hence $$X\subset\bigcup_{j=d+1}^nU_{i_j}\cup \bigcup_{j=1}^dV_j\subset \bigcup_{j=1}^nU_{i_j}.$$
   \end{proof}
   Now, if $f:Y\longrightarrow Z$\, is a morphism of preschemes, then $f$ is also continuous for the constructible topologies on $Y$ and $Z$\,. Namely, if $A\subset Z$\, is closed and $U\subset Z$\, is open, then $f^{-1}(A\cap U)=f^{-1}(A)\cap f^{-1}(U)$\, and $f^{-1}(A)$\, is closed and $f^{-1}(U)$\, is open in $Y$, since $f$ is continuous for the Zariski topologies on $Y$ and $Z$.\\
   \begin{definition} Let $\widehat{X}:=\mbox{projlim}_{i\in I}X_i$\, be a projective limit of preschemes in the category of of locally ringed spaces.
    We define the constructible topology on $\widehat{X}$\, to be the projective limit topology over the constructible topologies of the preschemes $X_i,i\in I$\,. \\
    I.e., for  $\mathcal F\subset I$\, a finite subset, we declare the set
    $$\left\{ \bigcap_{i\in \mathcal F}p_{i}^{-1}(A_{ij}\cap U_{ij})\mid A_{ij}\subset X_{i}\,\,\mbox{Zariski closed and}\,\, U_{ij}\subset X_{i}\,\,\mbox{Zariski open}\right\}$$
    to be a basis for the constructible topology on $\widehat{X}$. 
    \end{definition}
 \begin{theorem}\mylabel{thm:Chevalley}(Theorem of Chevalley)  Let $\widehat{X}:=\mbox{projlim}_{X_j\longrightarrow X_i}X_j$ be a projective limit of an indexed system $((X_i)_{i\in I}, f_{ij})$\,  with $X_i$\,  quasicompact preschemes  over the base field $k$ with proper and quasicompact transition morphisms $f_{ij}$\, in the category of locally ringed spaces . Then $\mid \widehat{X}\mid$\, with the Zariski topology as well as the constructible topology defined as above is a quasicompact topological space.\\
\end{theorem}
\begin{proof} We first treat the case where each $X_i, i\in I$\, is a quasicompact scheme. On each  $X_i,i\in I$\, we put the constructible topology where a basis of open subsets is given by intersections $A\cap U$\, with $A$ Zariski-closed and $U$ Zariski-open in $X$.
 Since each space $X_i,i\in I$ is by \ref{lem:L23} a Hausdorff quasicompact space with the constructible topology, the projective limit sits inside the topological product $\prod_{i\in I}X_i$\, as a closed subspace.\\
   Namely, we have 
   \begin{align}\mbox{projlim}_{X_j\longrightarrow X_i}X_j=\bigcap_{i<j}p_{ij}^{-1}\Gamma_{f_{ij}}\subset \prod_{i\in I}X_i,
   \end{align}
   where $$p_{ij}: \prod_{i\in I}X_i\longrightarrow X_i\times X_j$$
   is the canonical projection and $$\Gamma_{f_{ij}}\subset X_i\times X_j$$ is the graph of the transition  morphism $f_{ij}$\,.  It is a well known fact from basic topology that the graph $\Gamma$\, of a continuous map of Hausdorff spaces $g:Y\longrightarrow Z$\, with $Y$ quasicompact is a closed subset of $Y\times Z$\, with the product topology. For the convenience of the reader, we include a proof.\\ Let $(y_0,z_0)\in Y\times Z\backslash\Gamma_f$\,. Because the topological product $Y\times Z$\, is again Hausdorff, for each $y\in Y$\, there are  open neighbourhoods  $U_y$\, of $(y,f(y))$\, and $V_y$\, of $(y_0,z_0)$\, such that $U_y\cap V_y=\emptyset$\,. By quasicompactness of $\Gamma_f\cong Y$\, there exists a finite open subcovering $\Gamma_f\subset \bigcup_{i=1}^n U_{y_i}$\, of the covering $\Gamma_f\subset \bigcup_{y\in Y}U_y$\,. We have 
   $$\Gamma_f\cap \bigcap_{i=1}^nV_{y_i}\subseteq \bigcup_{i=1}^nU_{y_i}\cap \bigcap_{i=1}^nV_{y_i}\subseteq \bigcup_{i=1}^n(U_{y_i}\cap V_{y_i})=\emptyset.$$
   Hence for each point $(y_0,z_0)\in Y\times Z\backslash\Gamma_f$\, there exists an open neighbourhood $U=\bigcap_{i=1}^nV_{y_i}$\, with $U\cap \Gamma_f=\emptyset$\, and hence $\Gamma_f$\, is closed in $Y\times Z$\,.  \\
   By Tychonov's theorem, the  space $\prod_{i\in I}X_i$\, is  quasicompact when equipped with the product of the quasicompact constructible topologies of the $X_i, i\in I$\,. Hence $\widehat{X}$\, sits as a closed subset inside the quasicompact space $\prod_{i\in I}X_i$\,  and is thus again quasicompact with respect to the constructible topology. As the Zariski topology is coarser than the constructible topology, $\widehat{X}$\, is also quasicompact with respect to the Zariski topology.\\
   If each $X_i,i\in I$ is in general only a quasicompact prescheme, fix $i_0\in I$\, and  let $X_{i_0}=\bigcup_{i=1}^n \Spec A_i$\,be a finite affine open covering. Then $$\widehat{X}=\bigcup_{i=1}^n p_{i_0}^{-1}(\Spec A_i).$$
   As each $f_{ij}$\, is assumed to be proper and quasicompact, each $f_{i_0j}^{-1}(\Spec A_i)$\, is a quasicompact scheme, so $p_{i_0}^{-1}(\Spec A_i)$\, is a projective limit of quasicompact schemes  for $i=1,...,n$\,. If $$\widehat{X}=\bigcup_{j\in J}U_j$$ is any open covering, either for the constructible or the Zariski topology, then we consider the induced coverings $$p_{i_0}^{-1}(\Spec A_i)\subset \bigcup_{j\in J}U_j, i=1,...,n$$ which have by the previously treated case finite open subcoverings $$p_{i_0}^{-1}(\Spec A_i)\subset \bigcup_{j=1}^{n_i}U_{i_j}.$$
   Putting together these finitely many finite open coverings of the $p_{i_0}^{-1}(\Spec A_i)$\, we get a finite open subcovering $$\widehat{X}\subseteq \bigcup_{i=1}^n\bigcup_{j=1}^{n_i}U_{i_j}.$$
   \end{proof}
   We prove the following application of the proceeding theorem.
   \begin{proposition}
   \mylabel{prop:P555}
   Let $((X_i)_{i\in I},f_{ij})$\, be a projective system of preschemes of finite type with proper  surjective transition morphisms of finite type.\\
 Let $\widehat{X}:=\mbox{projlim}_{i\in I}X_i$\, be its projective limit in the category of locally ringed spaces. Then for all $i\in I$\,, the projection morphisms $p_i:\widehat{X}\longrightarrow X_i$\, are also surjective.
   \end{proposition}
   \begin{proof} Suppose to the contrary that there exists $i_0\in I$\, such that $$p_{i_0}:\widehat{X}\longrightarrow X_{i_0}$$ is not surjective. Let $x\in X_{i_0}\backslash p_{i_0}(\widehat{X})$\, be any point. Let for all $j>i_0$\, $Z_j:=f_{i_0j}^{-1}(\{x\})$\,. Then as each $f_{i_0j}$\, is surjective, the $Z_j, j\geq i_0$\, are all nonempty. Obviously, for $l>j\geq i_0$\, we have $$f_{jl}(Z_l)=f_{jl}(f_l^{-1}(\{x\})\subset f_j^{-1}(\{x\})=Z_j.$$
   As each $f_{jl}$\, is surjective, we even have $f_{jl}(Z_l)=Z_j$\,. So $((Z_j)_{j\in I, j>i_0}, f_{jl}\mid_{Z_l})$\, forms again a projective system of preschemes of finite type with proper surjective transition morphisms of finite type. Each $Z_j$\, is even a scheme of finite type over $\kappa(x)$ because $f_{i_0j}$\, was assumed to be proper and of finite type. Each $$\widehat{z}\in \mbox{projlim}_{j\in I, j\geq i_0}Z_j\subset \mbox{projlim}_{j\in I}X_j$$ satisfies $p_{i_0}(\widehat{z})=x$\,. But by assumption, there is no such $\widehat{z}$\,. So it follows $$\mbox{projlim}_{j\in I,j\geq i_0}Z_j=\emptyset.$$
   We apply the formula (1) with notation as in the proof of \ref{thm:Chevalley} to this projective system:
   $$\mbox{projlim}_{j\in J, j\geq i_0}Z_j=\bigcap_{j<l}p_{jl}^{-1}(\Gamma_{f_{jl}\mid_{Z_l}})\subset \prod_{j\in I, j\geq i_0}Z_j.$$ 
   We equip each $Z_j, j\geq i_0$\, with the constructible topology and the above product with the product topology. Then by \ref{lem:L23} each $Z_j$\, is a quasicompact Hausdorff space and the product $\prod_{j\in I, j\geq i_0}Z_j$\, is again Hausdorff and quasicompact by the theorem of Tychonoff. We have deduced that $$\bigcap_{j<l}p_{jl}^{-1}(\Gamma_{f_{jl}\mid_{Z_l}})=\emptyset.$$
   By quasicompactness, there are finitely many indices $j_i,l_i, i=1,...,n$\, such that
   $$\bigcap_{i=1}^np_{j_il_i}^{-1}(\Gamma_{f_{j_il_i}\mid_{Z_{l_i}}})=\emptyset.$$
   We now consider the finite projective system $\mathcal P$\, consisting of all $Z_{j_i},Z_{l_i}, i=1,...,n$\, and all morphisms $f_{ij}$\, of the projective system $(Z_i, f_{ij})$\, between any two of these objects.  The relations for a point $z\in \prod_{i=1}^nZ_{j_i}\times Z_{l_i}$\, to lie in $$\bigcap_{i=1}^np_{j_il_i}^{-1}(\Gamma_{f_{j_il_i}\mid_{Z_{l_i}}})$$ are part of the relations to lie in the projective limit of the system $\mathcal P$\,. So we have 
   $$\mbox{projlim}_{Z_i\in \mathcal P}Z_i\subset\bigcap_{i=1}^np_{j_il_i}^{-1}(\Gamma_{f_{j_il_i}\mid_{Z_{l_i}}})\subset\prod_{i=1}^nZ_{j_i}\times Z_{l_i}.$$
   So it follows $\mbox{projlim}_{Z_i\in \mathcal P}Z_i=\emptyset$\,. Let $d\in I$\, such that $d>j_i, d>l_i\forall 1\leq i\leq n$\,. Then the morphisms $f_{j_id}\mid_{Z_d}$\, and $f_{l_id}\mid_{Z_d}, i=1,...,n$\, commute with all transition morphisms of the projective system $\mathcal P$\, because $(Z_i, f_{ij})$ is a projective system and for $i<j<d$ $f_{id}=f_{ij}\circ f_{jd}$\, holds. So we get a well defined morphism $$Z_d\longrightarrow \mbox{projlim}_{Z_i\in \mathcal P}Z_i.$$ As $Z_d$ is nonempty  we have reached a contradiction. Hence our assumption was false and $p_i:\widehat{X}\longrightarrow X_i$\, is surjective for all $i\in I$\,.
   \end{proof}
 \begin{definition}
\mylabel{def:D103} Let $\widehat{X}:=\mbox{projlim}_{i\in I}X_i$\, be as above a projective limit of preschemes in the category of locally ringed spaces with proper transition morphisms of finite type. Let  $\mathcal Z\subset \widehat{X}$\, be any subset of $\widehat{X}$\,. For each $i\in I$  we call $\mathcal Z_{X_i}:=p_i(\mathcal Z)$\, the trace of $\mathcal Z$\,on $X_i$.
\end{definition}
\begin{theorem}\mylabel{thm:A310} With the above notation, for all $i\in I$\, the morphism $p_i:\widehat{X}\longrightarrow X_i$\, is a closed map for the Zariski topologies on $\widehat{X}$\, and $X_i,i\in I$\,.\\
For any closed subset $\mathcal Z\subset \widehat{X}$\,
 one has $\mathcal Z=\cap_{i\in I}p_i^{-1}\mathcal Z_{X_i}$\, where the intersection can be replaced by the intersection over any final subset $I_0\subset I$\, .\\
Conversely, if one is given for each $i\in I$ a Zariski closed subset $ Z_i\subset X_i$\, such that for each $i<j$ $f_{ij}(Z_{j}) = Z_i$\,, then $\mathcal Z:=\cap_{i\in I}p_i^{-1}( Z_i)$\, is Zariski closed in $\widehat{X}$\, with trace $(\mathcal Z)_{X_i}=Z_i$\, for each $i\in I$.
\end{theorem}
\begin{proof} Let $\mathcal Z\subset \widehat{X}$ be a closed subset. By definition of the Zariski topology,  there is $I_0\subset I$ and for each $i\in I_0$\, a closed subset $Z_i\subset X_i$\, such that $\mathcal Z=\bigcap_{i\in I_0}p_i^{-1}(Z_i)$\,.\\
Let $I_1\subset I$\, be the final subset of all $i\in I$\, for which there is an $j\in I_0$\, with $j\leq i$\,. We define  
$$Z_i':=\bigcap_{j\leq i, j\in I_0}f_{ji}^{-1}(Z_j)\,\,\mbox{for}\,\, i\in I_1.$$ Then 
\begin{align*}\bigcap_{i\in I_1}p_i^{-1}(Z_i')=\bigcap_{i\in I_1}p_i^{-1}(\bigcap_{j\leq i,j\in I_0}f_{ji}^{-1}(Z_j))\\=\bigcap_{i\in I_1}\bigcap_{j\leq i, j\in I_0}p_i^{-1}(f_{ji}^{-1}(Z_j))= \bigcap_{i\in I_1}\bigcap_{j\leq i, j\in I_0}p_j^{-1}(Z_j)\\
=\bigcap_{j\in I_0}p_j^{-1}(Z_j)=\mathcal Z.
\end{align*}
So we can assume that $I_0=I_1$\, and replace each $Z_i$ with $Z_i'$\,.\\ 
For $i\in I$\, arbitrary, we define
$$Z_i':=\bigcap_{j\geq i, j\in I_1}f_{ij}(Z_j).$$
As each $f_{ij}$\, is by assumption a proper morphism and all $Z_i$\, are Zariski closed, so are the $Z_i', i\in I$\,.\\
We have
\begin{align*}
\bigcap_{i\in I}p_i^{-1}(Z_i')=\bigcap_{i\in I}p_i^{-1}(\bigcap_{j\geq i, j\in I_1}f_{ij}(Z_j))\\
= \bigcap_{i\in I}\bigcap_{j\geq i, j\in I_1}p_i^{-1}(f_{ij}(Z_j))\supseteq \bigcap_{i\in I}\bigcap_{j\geq i, j\in I_1}p_j^{-1}(Z_j)\\
=\bigcap_{j\in I_1}p_j^{-1}(Z_j)=\mathcal Z.
\end{align*}
On the other hand, 
\begin{align*}
\bigcap_{i\in I}p_i^{-1}(Z_i')=\bigcap_{i\in I}p_i^{-1}(\bigcap_{j\geq i, j\in I_1}f_{ij}(Z_j))\\
\subseteq \bigcap_{i\in I_1}p_i^{-1}(\bigcap_{j\geq i, j\in I_1}f_{ij}(Z_j))\subseteq \bigcap_{i\in I_1}p_i^{-1}(Z_i)\\
=\mathcal Z.
\end{align*}
Thus $$\bigcap_{i\in I}p_i^{-1}(Z_i')=\bigcap_{i\in I_1}p_i^{-1}(Z_i)$$ and we may replace the $Z_i$\, by the $Z_i'$\, and are now in the following situation:\\
For each $i\in I$\, we are given a Zariski closed subset $Z_i\subset X_i$\, such that 
$$\mathcal Z=\bigcap_{i\in I}p_i^{-1}(Z_i).$$
I claim that $,\forall i\in I :Z_i=\mathcal Z_{X_i}$\,.\\
We have 
\begin{align*}\mathcal Z_{X_i}=(\bigcap_{j\in I}p_j^{-1}(Z_j))_{X_i}\subseteq (p_i^{-1}(Z_i))_{X_i}\\
=p_i(p_i^{-1}(Z_i))=Z_i,
\end{align*}
thus $\mathcal Z_{X_i}\subseteq Z_i$\,.
Obviously one has $\mathcal Z\subset \bigcap_{i\in I}p_i^{-1}(\mathcal Z_{X_i})$\, as   $$\mathcal Z\subset p_i^{-1}(p_i(\mathcal Z))\,\,\forall i\in I.$$
So we have the following chain of relations
$$\mathcal Z\subseteq \bigcap_{i\in I}p_i^{-1}(\mathcal Z_{X_i})\subseteq \bigcap_{i\in I}p_i^{-1}(Z_i)=\mathcal Z.$$
We conclude that equality holds at each place, i.e., $\mathcal Z=\bigcap_{i\in I}p_i^{-1}(\mathcal Z_{X_i})$\, and $$\bigcap_{i\in I}p_i^{-1}(\mathcal Z_{X_i})=\bigcap_{i\in I}p_i^{-1}(Z_i).$$
 These intersections are the same as the intersections restricted to any cofinal directed subsystem of the directed system .\\ For, for both systems $$(\mathcal Z_{X_i},i\in I, f_{ij}\mid_{\mathcal Z_{X_j}})\,\,\mbox{and}\,\, (Z_i,i\in I, f_{ij}\mid_{Z_j})$$ we have $f_{ij}(\mathcal Z_{X_j})\subset \mathcal Z_{X_i}$\, and $f_{ij}(Z_j)\subset Z_i$\,. Let $I_0\subset I$\, be a cofinal subsystem. For each $j\in I$\, there is $i_j\in I_0$\, with $j<i_j$\,.$$\mbox{ As}\,\, f_{ji_j}(Z_{i_j})\subset Z_j\,\, \mbox{we have}\,\, p_{i_j}^{-1}(Z_{i_j})\subseteq p_j^{-1}(Z_j).$$ Thus
$$\bigcap_{i_j\in I_0}p_{i_j}^{-1}(Z_{i_j})\supseteq \bigcap_{j\in I}p_j^{-1}(Z_j)\supseteq \bigcap_{i_j\in I_0}p_{i_j}^{-1}(Z_{i_j})$$ and it follows $\bigcap_{j\in I}p_j^{-1}(Z_j)=\bigcap_{i\in I_0}p_i^{-1}(Z_i)$\,. The case with the $\mathcal Z_{X_i}$\, is the same.\\
From this and $$\mathcal Z_{X_i}\subseteq \overline{\mathcal Z_{X_i}}\subseteq Z_i$$
we get $$\bigcap_{i\in I}p_i^{-1}(\mathcal Z_{X_i})=\bigcap_{i\in I}p_i^{-1}(\overline{\mathcal Z_{X_i}}).$$
Suppose there is $i\in I$\, such that $\overline{\mathcal Z_{X_i}}\backslash \mathcal Z_{X_i}\neq \emptyset$\,.\\
Let $z_i\in \overline{\mathcal Z_{X_i}}\backslash \mathcal Z_{X_i}$\, be any point. Then as $p_i(\mathcal Z)=\mathcal Z_{X_i}$\, we have $p_i^{-1}(\{z_i\})\cap \mathcal Z=\emptyset$\,. Since each $f_{ij}$\, is assumed to be proper and of finite type, the fibre product 
$$\mathcal Z\times_{X_i}\{z_i\}=\mbox{projlim}_{j\geq i}(X_j\times_{X_i}\{z_i\})$$
 is a projective limit of schemes of finite type over the base field $\kappa(z_i)$\, and by \ref{thm:Chevalley} quasicompact. \\
From $$\emptyset =p_i^{-1}(\{z_i\})\cap \mathcal Z=p_i^{-1}(\{z_i\})\cap \bigcap_{j\in I}p_j^{-1}(\overline{\mathcal Z_{X_j}})$$
it follows that there are indices $j_1,...,j_n$\, such that already $$p_i^{-1}(\{z_i\})\cap \bigcap_{k=1}^np_{j_k}^{-1}(\overline{\mathcal Z_{X_{j_k}}})=\emptyset.$$
Let $j_0\in I$\, such that $j_0\geq j_k, k=1,...,n\,\,\mbox{and}\,\, j_0\geq i$\,. Then $$p_{j_0}^{-1}(\overline{\mathcal Z_{j_0}})\subseteq \bigcap_{k=1}^np_{j_k}^{-1}(\overline{\mathcal Z_{X_{j_k}}})$$ and hence $$\emptyset =p_i^{-1}(\{z_i\})\cap p_{j_0}^{-1}(\overline{\mathcal Z_{X_{j_0}}}))=p_{j_0}^{-1}(f_{ij_0}^{-1}(\{z_i\})).$$
We have $$f_{ij_0}(\mathcal Z_{X_{j_0}})= f_{ij_0}(p_{j_0}(\mathcal Z))=p_i(\mathcal Z)=\mathcal Z_{X_i}$$ and hence    $$f_{ij_0}\mid_{\overline{\mathcal Z_{X_{j_0}}}}:\overline{\mathcal Z_{X_{j_0}}}\longrightarrow \overline{\mathcal Z_{X_i}}$$ is surjective because $f_{ij_0}$\, is proper, $f_{ij_0}(\overline{\mathcal Z_{X_{j_0}}})$\, is closed and contains $\mathcal Z_{X_i}$.\\
The set $f_{ij_0}^{-1}(\{z_0\})$\, is therefore nonempty. By \ref{prop:P555}, the projection morphism $$p_{j_0}: \mbox{projlim}_{i\in I}\overline{\mathcal Z_{X_i}}\longrightarrow \overline{\mathcal Z_{X_{j_0}}}$$ of the projective system $$((\overline{\mathcal Z_{X_i}})_{i\in I},f_{ij}\mid_{\overline{\mathcal Z_{X_j}}})$$
with surjective transition morphisms is itself surjective.\\
 Therefore, the set $p_{j_0}^{-1}(f_{ij_0}^{-1}(\{z_i\}))$\, has to be also nonempty, a contradiction.\\
 Thus for each $i\in I$\, $\mathcal Z_{X_i}=\overline{\mathcal Z_{X_i}}$\,.\\
 Thus the traces $\mathcal Z_{X_i}$\, are Zariski closed for each Zariski closed subset $\mathcal Z\subset \widehat{X}$\, and each $i\in I$\,. Fixing $i$ and letting $\mathcal Z$\, vary it follows that the projections $p_i:\widehat{X}\longrightarrow X_i$\, are closed maps.\\ 

To prove the last statement, as arbitrary intersections of closed subsets are again closed, it is clear that $\mathcal Z$\, is Zariski closed in $R(K(X))$\,. Obviously for each complete $X_i$\, we have $$p_i(\mathcal Z_{X_i})=p_i(\bigcap_{j\in I}p_j^{-1}(Z_j))\subset p_i(p_i^{-1}(Z_i))=  Z_i.$$
 Conversely, we consider $((Z_i)_{i\in I}, g_{ij}:=f_{ij}\mid_{Z_j})$\, as a projective system of preschemes locally of finite type with proper surjective transition morphisms of finite type. By \ref{prop:P555} it follows that for all $i\in I$\, the projection $g_i: \widehat{Z}=\lim_{j\in I}Z_j\longrightarrow Z_i$\, , which is the restriction of the projection $p_i: \widehat{X}\longrightarrow X_i$\,, is surjective.\\
 We have $\widehat{Z}\subset \widehat{X}$\, and even because of $g_i(\widehat{Z})=p_i(\widehat{Z})\subset Z_i$\, that $$\widehat{Z}\subset \bigcap_{j\in I}p_j^{-1}(Z_j).$$ 
 It follows that $p_i\mid_{\mathcal Z}: \bigcap_{j\in I}p_j^{-1}(Z_j)\longrightarrow Z_i$\, is surjective for all $i\in I$\, and $\mathcal Z_{X_i}=Z_i$\,.  
\end{proof}
\begin{corollary}\mylabel{cor:C10} Let $\mathcal Z\subset \widehat{X}$\, be a Zariski closed subset. Let $Z_i:=(\mathcal Z)_{X_i}$\,. Then $\mathcal Z\cong \mbox{projlim}_{i\in I}Z_i$\, as topological spaces.
\end{corollary}
\begin{proof} First, by the proceeding theorem, each $Z_i\subset X_i$\, is a Zariski closed subset of $X_i$.\\
Secondly, also by the proceeding theorem, we have $\mathcal Z=\bigcap_{i\in I}p_i^{-1}(Z_i)$\,. In the course of the above proof we had already shown that there is an inclusion $$\mbox{projlim}_{i\in I}Z_i\hookrightarrow \bigcap_{i\in I}p_i^{-1}(Z_i).$$
 To show the other direction, let $(z_i)_{i\in I}\in \bigcap_{i\in I}p_i^{-1}(Z_i)$\, be given. Because $$(z_i)_{i\in I}\in \mbox{projlim}_{i\in I}X_i\,\, \mbox{we have}\,\, f_{ij}(z_j)=z_i.$$ Furthermore, $$\forall l\in I, p_l((z_i)_{i\in I})\in (\mathcal Z)_{X_l}=Z_l$$
again by the proceeding theorem. Hence $(z_i)_{i\in I}$\, defines an element in $\mbox{projlim}_{i\in I}Z_i$\,.
\end{proof}
\paragraph{Zariski-Riemann varieties of integral preschemes  of finite type}

Preschemes are the most general objects of Grothendieck algebraic geometry. Among them, the schemes are the separated preschemes. Separatedness and properness are most conveniently checked by applying the valuative criterion of separatedness and properness (see \cite{Ha}[part II, chapter 4, Theorem 4.3, p.97, Theorem 4.7, p.101]. As Riemann varieties, defined in the sequel, are the global objects of all valuations of a function field, it is reasonable  to use them for the study of   phenomena in the theory of integral preschemes.\\ 
The notion of the Zariski-Riemann variety of a function field is classical and goes back to O. Zariski. We extend the definition to Riemann varieties of integral preschemes. An integral prescheme can look much different from the basic example of the affine line with the origin doubled. In higher dimensions, we may glue different birational models along the locus where they are isomorphic. Going up to the level of Riemann varieties, this phenomenon vanishes. "Nonseparated Riemann varieties", i.e., Riemann varieties of  arbitrary integral preschemes look all like the "separated" Riemann variety of the function field with closed subsets being doubled or in general $n$-tuplified. This opens up the possibility to study the structure of integral preschemes via the simpler object of their associated Riemann varieties.\\  
Riemann varieties  are locally ringed spaces all of its local rings are valuation rings of a given function field. They were originally invented  in the context of desingularization theory in order to pass from local uniformizations to global resolutions. Although they are not schemes, in some sense they are higher dimensional analogues of smooth curves (see \cite{Val}[XXX2.2]).\\
 One the one hand, one may associate to each integral prescheme of finite type a Riemann variety, on the other hand, one can associate to an arbitrary field extension $K/k$\, such an object. The two cases overlap if $K/k$ is a function field since as we show later the Riemann variety $R(K/k)$\, is canonically isomorphic to the Riemann variety of any complete model of $K/k$\,. We will need the more general case where $K/k$\, is not necessarily finitely generated but only of finite transcendence degree in our applications, mainly in the situation where $K$ is the residue field $\kappa(\nu)$\, of a valuation ring $A_{\nu}$ of some function field $L/k$\,. Before we  recall the definitions we give the proof of a fundamental lemma that we will need throughout this work.
\begin{lemma}
\mylabel{lem:L104} Let $X_i,i=1,...,n$\, be a collection of integral schemes of finite type with quotient field $K=K(X_i)$\,. Then there is a normal projective model $Y$\,of $K/k$ and open subsets $U_i\subset Y, i=1,...,n $\, and proper surjective morphisms\\
 $p_i:U_i\longrightarrow X_i$\,. 
\end{lemma}
\begin{proof} First, we compactify each $X_i$ to a complete integral scheme of finite type $Y_i$ which is always possible by the compactification theorem of Nagata (see \cite{Nagata}).\\
 Put $Z_1:=Y_1$\, and, inductively, choose $Z_{k+1}$ to be a complete model of $K/k$\, dominating both $Z_k$\, and $X_{k+1}$\,, such as does the closure of $$U\stackrel{\Delta}\hookrightarrow U\times U\hookrightarrow Z_k\times X_{k+1}$$ where $U\subset X_{k+1},Z_k$\, is a common open subset of the birationally equivalent varieties $X_{k+1}$\, and $Z_k$\, and $\Delta$ is the diagonal. Put $Y':=Z_n$\, and take as $Y$ any projective model lying over $Y'$.\\
By construction, there are proper birational morphisms $q_i:Y\longrightarrow Y_i$\,. Put $U_i:=q_i^{-1}(X_i)$\, and $p_i: U_i\longrightarrow X_i$\,, the restriction of $q_i$\, to $U_i$\,.
\end{proof}
\begin{definition}(1)
\mylabel{def:D99}\,\, Let $X$ be an integral prescheme locally of finite type. Let $I$ be the set of all integral normal preschemes $X'$\, together with a proper birational morphism $f_i: X'\longrightarrow X$\,. We define $i<j$\, if there exists a proper birational morphism $f_{ij}: X_j\longrightarrow X_i$\, such that $f_i\circ f_{ij}=f_j$\,.\\
We define the Zariski-Riemann variety $R(X)$\, of the prescheme $X$ to be the projective limit $\mbox{projlim}_{i\in I}X_i$\, over the directed system $((X_i)_{i\in I},f_{ij})$ in the category of locally ringed spaces.
\end{definition}
We give now the more classical definition.
\begin{definition}(2)\,\,
\mylabel{def:D100} Let $\Spec A$ be an integral affine variety. We construct a locally ringed space $R(\Spec A)$\, whose points are the  valuations $\nu$\, of $K(A)$\,with $A\subset A_{\nu}$\,.\\
 A basis of open subsets will be the sets $R\Spec A'$\, of  valuations $\nu$\,of $K(A)$\, with $A'\subset A_{\nu}$\, where $A'$\, is a finitely generated $k$-algebra  with the same quotient field as $A$ and there is an inclusion of finitely generated $k$-algebras  $A\longrightarrow A'$\,.\\
  The sheaf of rings will be the one induced by the presheaf inside the constant sheaf $\underline{K(X)}$, that takes the value $(A')^n$\, on the open set $R(\Spec A')$\, for $A'\subset K(X)$\, finitely generated over $k$ and $A\subset A'$\,.  \\
If $X$ is an integral  prescheme locally of finite type over $k$ with affine open cover $X=\bigcup_{i\in I} \Spec A_i$\,, the gluing datum $U_{ij}: =\Spec A_i\cap_X\Spec A_j$\, determines a gluing datum 
 $$\mathcal U_{ij}:=R(U_{ij}):=\left\{\nu\in R\Spec A_i\mid \nu\,\mbox{has}\,\,\mbox{center}\,\,\mbox{in}\,\,\Spec A_i\cap \Spec A_j\right\}$$
 along which the Riemann varieties $(R(\Spec A_i))_{i\in I}$\, glue in a natural manner to a Riemann variety $R(X)$\,.
\end{definition}
\begin{definition}(3) Let $K/k$\, be an arbitrary field extension. We define the Riemann variety $R(K/k)$\, as a point set to be the set of all Krull valuations of $K$ being trivial when restricted to $k$\,. As basic open sets we define the subsets $$R^K\Spec A:=\{\nu\in R(K/k)\mid A\subset A_{\nu}\}$$ for each finitely generated (!) $k$-subalgebra $A$ of $K$.\\
We define a locally ringed space structure on $R(K/k)$\, by defining the structure sheaf $\mathcal O_{R(K/k)}$\, to be the sheaf associated to the presheaf $\mathcal O_{R(K/k)}'$\, generated inside the constant sheaf $\underline{K}$\, by the requirement that it takes the value $A^{n,K}$ on the open subset $R^K\Spec A$\,.
\end{definition} 

\begin{remark} By saying that the presheaf $\mathcal O_{R(K/k)}'$\, is generated by the requirement 
$$\mathcal O_{R(K/k)}'(R^K\Spec A)=A^{n,K}$$
 we mean that it is the intersection of all presheaves $\mathcal F\subset \underline{K(X)}$\, inside the presheaf $\underline{K(X)}$\, with $\mathcal F(R^K\Spec A)\supseteq A^{n,K}$\,.
\end{remark}
From the first definition and \ref{thm:Chevalley}, we have the following immediate 
\begin{corollary}\mylabel{cor:C11} Let $X$ be an integral prescheme of finite type. Then the Riemann variety 
$R(X)$\,either with the Zariski topology or the constructible topology is a quasicompact topological space.\\
\end{corollary}
\begin{proof}
\end{proof}
Although arbitrary valuation rings are more difficult to handle in that they are not noetherian, this fact makes them worth to be considered because as we shall see later, restricting to discrete algebraic rank one valuations, the space $R^1(K(X))$\, is definitely not quasicompact. \\
\begin{remark}
\mylabel{rem:R100}  F.V. Kuhlmann proved in \cite{Kuhl} the result, saying that the Riemann variety is still quasicompact with respect to the constructible topology (see\cite{Kuhl}[Appendix, Theorem 36, p.21] by using model theoretic methods. 
\end{remark}
\begin{proposition}\mylabel{prop:P400}
 \begin{enumerate}[i]
\item The above notion of the Riemann variety $R(X)$\, of an integral prescheme of finite type as a locally ringed space in (2)  and of $R(K/k)$\, in (3) is well defined.
\item If $f:X_1\longrightarrow X_2$\, is a proper birational morphism of integral preschemes, inducing the identity on the fields of rational functions then $f$\, induces a morphism $R(f): R(X_1)\longrightarrow R(X_2)$\, that maps a valuation $\nu$\, with center $x_1\in X_1$\, to itself as an element in $R(X_2)$\, with center $f(x_1)\in X_2$\,. $R(f)$\, is a homeomorphism of Zariski topological spaces.
 \item The assignment $R^K\Spec A\mapsto A^{n,K},$\, $A\subset K(X)$\, a finitely generated $k$ algebra, defines a presheaf $\mathcal O_R'$\, of $k$-algebras on $R^K\Spec A$\, whose associated sheaf $\mathcal O_R$\, satisfies $\mathcal O_R(R^K\Spec A)= A^{n,K}$\, for each  finitely generated $k$-subalgebra $A\subset K(X)$\,.
\item The stalk at a valuation $\nu$ is  the valuation ring $A_{\nu}$\, associated to $\nu$\,.\\
\item The definitions (1) and (2) are equivalent.
\end{enumerate}
\end{proposition}
\begin{proof} For a proof in the particular case of $X$ being a complete integral scheme of finite type, see also  \cite{Val}[Theorem 2.11].\\
\begin{enumerate}[i]
\item First, we have to show that if $A_1,A_2\subset K(X)=K$\, are finitely generated  and $R\Spec A_1=R\Spec A_2$\,, then $A_1^{n,K}=A_2^{n,K}$\,  such that the above assignment is well defined. By \cite{Matsumura}[chapter 4, paragraph 12, Theorem 12.4.(i), p.88], we have $$A_1^{n,K}=\bigcap_{\nu\in R\Spec A_1}A_{\nu}=\bigcap_{\nu\in R\Spec A_2}A_{\nu}=A_2^{n,K}$$
the intersections taken inside $K(X)$\,. So the above assignment in (2) and (3) is well defined.\\
\item Let $R(f): R(X_1)\longrightarrow R(X_2)$\, be defined as above. If $\Spec A\subset X_2$\, is an affine open subset, then $R\Spec A$\, is an open subset of $R(X_2)$\,and $R(f)^{-1}(R\Spec A)$\, consists of all $\nu\in R(X_1)$\, such that $f(c_{X_1}(\nu))\in \Spec A$\,,i.e., 
$$c_{X_1}(\nu)\in f^{-1}(\Spec A)\,\, \mbox{or}\,\, \nu\in R(f^{-1}(\Spec A)).$$
 By the continuity of $f, f^{-1}(\Spec A)$\, is open in $X_1$\, and by the definition of the Zariski topology on $R(X_1)$\,, $$R(f)^{-1}(R\Spec A)=R(f^{-1}(\Spec A))$$ is open in $R(X_1)$\,. Thus $R(f)$\, is continuous.  By the properness of $f$, each valuation $\nu\in R(X_2)$\, has a unique lift to a valuation $\nu'\in R(X_1))$\,. It follows that $f$\, is bijective. For each open affine subset $\Spec A\subset X_2$\,, we can consider $R\Spec A$\, and $Rf^{-1}(\Spec A)$\, canonically as subsets of $R(K(X_1)/k)$\,. By the above arguement, these two sets are identified as the same subset of $R(K(X_1)/k)$\,. They both have the same topology, for if $$f^{-1}(\Spec A)=\bigcup_{j=1}^n\Spec B_j,$$ then by definitions (2) and (3), we have $$R\Spec A=\bigcup_{j=1}^nR\Spec B_j$$ as subsets of $R(K(X_1)/k)$\,. It follows that $R(f)$\, is a bijective local homeomorphism, hence a homeomorphism.\\
\item If $R\Spec A\subset R\Spec B, A,B\subset K(X)$\, finitely generated over $k$, then again by \cite{Matsumura}[chapter 4, paragraph 12, Theorem 12.4 (i), p.88] $$B^{n,K}=\bigcap_{\nu\in R\Spec B}A_{\nu}\subset \bigcap_{\nu\in R\Spec A}A_{\nu}=A^{n,K},$$ so the above assignment defines a presheaf $\mathcal O_R'$. \\
We show next, that restricted to any $R\Spec B$\,, $\mathcal O_R'$\,  is already a sheaf.\\
 Let $$R\Spec B=\bigcup_{i\in I}R\Spec A_i$$ be an open affine covering. By quasicompactness of $R\Spec B$\, (\ref{thm:Chevalley}) there is a finite subcovering $$R\Spec B=\bigcup_{i=1}^nR\Spec A_i.$$
  By \ref{lem:L104} we find a normal complete model $Y$ plus open subsets\\
 $U_i,i=1,...,n$\, and $U$ plus proper birational morphisms 
 $$p_i:U_i\longrightarrow \Spec A_i, i=1,...,n\,\, \mbox{and}\,\, q:U\longrightarrow \Spec B.$$ We must have $U=\bigcup_{i=1}^n U_i$\, since by the second part of this proposition, we have $R\Spec B=R(U)$\, and $R\Spec A_i=R(U_i), i=1,...,n$\,.\\
 From the properness of the morphisms $p_i,i=1,...,n$\, and $q$ and the normality of $Y$\, we further conclude 
\begin{align*}\Gamma(U_i,\mathcal O_Y)=(A_i)^{n,K}=\mathcal O_R'(R\Spec A_i),i=1,...,n\\ \mbox{and}\,\,\Gamma(U,\mathcal O_Y)=B^{n,K}=\mathcal O_R'(R\Spec B).
\end{align*}
  The sheaf condition for $\mathcal O_R'$\, restricted to $R\Spec B$\, reduces to the sheaf condition for $\mathcal O_Y$\, restricted to the open subset $U$ with respect to the covering $U=\bigcup_{i=1}^nU_i$\,. As it is satisfied for the structure sheaf $\mathcal O_Y$\, it is also satisfied for $\mathcal O_R'\mid_{R\Spec B}$\,. Thus $\mathcal O_R'$\, restricted to $R\Spec B$\, is already a sheaf.\\ 
\item Let $\nu\in R(K(X))$\, be a point. By definition, the stalk $\mathcal O_{R,\nu}$\, is equal to the inductive limit $\lim_{\nu\in R\Spec B}B^{n,K}$\,. The transition maps are injective since everything sits inside $K(X)$\,. As by \cite{Matsumura}[chapter 4, paragraph 12, Theorem 12.4 (i), p.88]:
$$B^{n,K}=\bigcap_{\mu\in R\Spec B}A_{\mu}$$ and since the intersection of all $R\Spec B$\, with $\nu\in R\Spec B$\,  is precisely  $\{\nu\}$\,, we have
$$\mathcal O_{R,\nu}= \lim_{\nu\in R\Spec B}B^{n,K}=\lim_{\nu\in R\Spec B}(\bigcap_{\mu\in R\Spec B}A_{\mu})\subseteq A_{\nu}.$$ 
Conversely, if $f\in A_{\nu}$\, and $B\subset A_{\nu}$\,, then also $B[f]\subset A_{\nu}$\, and so $\nu\in R^K\Spec B[f]$\,. We have 
$$\mathcal O_{R,\nu}=\lim_{\nu\in R^K\Spec B}B^{n,K}=\lim_{\nu\in R^K\Spec B[f]}(B[f])^{n,K}.$$ So $f\in \mathcal O_{R,\nu}$\, and $A_{\nu}\subset \mathcal O_{R,\nu}$\,. Hence the stalk $\mathcal O_{R,\nu}$\,  is equal to $A_{\nu}$\,.\\
\item We have to show, that 
\begin{enumerate}[a]
\item The projective limit in (1) is as a topological point space the same as $R(X)$\, defined in (2).
\item The structure sheaf of the projective limit evaluated at $R^{K(X)}\Spec A\subset R(X)$\, , $A\subset K(X)$\, a finitely generated $k$-subalgebra of $K(X)$\, is equal to $A^{n,K(X)}$.
\end{enumerate}
\begin{enumerate}[a]
 \item Let $X=\bigcup_{i\in I}\Spec A_i$\, be a Zariski open covering. By \cite{Val}[chapter 1.1, Definition p.1,Theorem 1.2, p.2], for each $$\nu\in R^{K(X)}\Spec A_i^n\subset R(K(X))$$ ($A_i^n$ denotes the normalization of $A_i$), the valuation ring $A_{\nu}$\, is the inductive limit of all local rings $B_j\subset A_{\nu}$\, essentially of finite type.  Let $\nu\in R\Spec A_i$\,.
 Now if $B_j$\, is any such local ring, let $C_j\subset K(X)
$\, be a finitely generated $k$-algebra such that $B_j$\, is some localization of $C_j$\,.\\
 By \ref{lem:L104}, there is a complete (projective)  normal model $Y$ with quotient field $K(X)$\, and open subsets $U,V\subset Y$ plus proper birational morphisms $$p: U\longrightarrow \Spec A_i$$ and $$q: V\longrightarrow \Spec C_j.$$
 The unique lifts of the centers of $\nu$\, in $\Spec A_i$\, and $\Spec C_j$\, via the proper birational morphisms $p,q$ have to coincide since $Y$ was assumed to be complete.\\
  Since $U$ and $\Spec A_i$\, are quasiprojective by \cite{Ha}[chapter II, 7, Theorem. 7.17, p.166] there is an ideal $\mathfrak{a}_i\subset A_i$\, such that $$p: U= \mbox{Bl}_{\mathfrak{a}_i}\Spec A_i\longrightarrow \Spec A_i$$ is a blowing up morphism.
   Let $\mathcal I_i\subset \mathcal O_X$\, be an ideal sheaf such that $\mathcal I(\Spec A_i)=\mathfrak{a}_i$\,.  $$(\mbox{e.g.,}\,\mathcal I_i=j_{i,*}\mathfrak{a}_i\cap \mathcal O_X \,\,\mbox{where}\,\, j_i: \Spec A_i\longrightarrow X$$ is the open immersion.)\\
Let $Z:=\mbox{Bl}_{\mathcal I_i}X$\, with canonical projection $\pi: Z\longrightarrow X$\,. We then have $U\subset Z$\, and $\pi\mid_U=p$\,. Now the local ring $\mathcal O_{Z,u}$\, of the center of $\nu$\, on $U\subset Z$ dominates the local ring $B_j$\,.\\
We have shown that in the projective system $((X_i)_{i\in I},f_{ij})$\, of integral preschemes $X_i,$\, proper and birational over $X$, the system of the local rings $\mathcal O_{X_i,x_i}$\, of the centers of a valuation $\nu\in R\Spec A_i$\, is final in the inductive system of all local rings essentially of finite type under $A_{\nu}$\,. Thus, the stalk of the structure sheaf of the projective limit taken at the point $\widehat{x}$ corresponding to the projective system of points $(x_i)$, being the inductive limit of the local rings $\mathcal O_{X,x}$\, is equal to $A_{\nu}$\, and $\widehat{x}=(x_i)$\, may be identified with $\nu$\,.\\
Conversely, each point $\widehat{x}$\, of the projective limit $\widehat{X}$\, is represented by a projective system of points $(x_i)_{i\in I}, x_i\in X_i$\, such that the local rings $\mathcal O_{X_i,x_i}$\, form an inductive system ordered by the relation of domination. By \cite{Val}[chapter 1.1, Definition p.1, Theorem 1.2 p.2] there is a valuation $$\nu\in\Spec A_i\subset  R(K(X))$$ whose valuation ring $A_{\nu}$\, dominates each local ring of this inductive system. By the above, \,$\lim_{i\in I}\mathcal O_{X_i,x_i}=A_{\nu}$\, and the point $\widehat{x}$\, can be identified with $\nu$\,.\\
So as a point space, $\mid\widehat{X}\mid$\, coincides with $R(X)$\, from definition (1).\\
We now look at the topologies.\\
A basic open subset of the projective limit $\widehat{X}=\mbox{projlim}_{X_j\longrightarrow X}X_j$\, is given by subsets of the form 
$$ \bigcap_{i=1}^np_i^{-1}(\Spec B_i), \,\, \Spec B_i\subset X_i\,\,\mbox{Zariski open in}\,\, X_i$$
and $ X_i\longrightarrow X$\, a proper and birational morphism. \\
  Let  $\Spec C$\, be an integral affine scheme of finite type with $A_i\subset C$\,.\\
Let $\nu\in R\Spec C$\,. By what we have just proven, we may identify $\nu$\, as a point 
$$\widehat{x}=(x_i)_{i\in I}\in \mbox{projlim}_{i\in I}X_i.$$
 It follows that there is $j\in I$\, such that $\mathcal O_{X_j,x_j}$\, dominates the $k$-algebra $C$\,, because $A_{\nu}$\, dominates $C$ and $A_{\nu}=\lim_{i\in I}\mathcal O_{X_i,x_i}$\,. As $\mathcal O_{X_j,x_j}$\, is essentially of finite type, there is an affine open $\Spec B_j\subset X_j$\, with $x_j\in \Spec B_j$\, and $C\subset B_j$\,. It follows $\nu\in R\Spec B_j$\,. I.e., for all $\nu\in \Spec C$\, there is $j\in I$\, and an open $\Spec B_j\subset X_j$\, such that $\nu\in \Spec B_j$\,.\\
Therefore, the topology in definition (2), whose basis is given by all  $$R\Spec C\subset R\Spec A_i$$ is generated  by the basis of open neighbourhoods $$p_i^{-1}(\Spec B_i)\subset R(X)$$ which is a basis for the topology in definition (1). As the topology in (1) is apriori coarser than the topology in (2), the topologies from (1) and (2) coincide. 
\item Let $X_i$\, be a normal integral prescheme proper and birational over $X$ and $\Spec A_i\subset X_i$\, be an open subset. For each integral normal prescheme $X_j$\,, proper and birational over $X_i$\, we have
\begin{align*}p_j^{-1}\mathcal O_{X_j}( p_i^{-1}(\Spec A_i))=p_j^{-1}\mathcal O_{X_j}(p_j^{-1}\circ p_{ij}^{-1}(\Spec A_i))\\
=\mathcal O_{X_j}(p_{ij}^{-1}(\Spec A_i))=(A_i)^{n,K(X)},
\end{align*}
 since $p_{ij}$\, is proper and $X_j$\, is assumed to be normal. Thus the inductive limit $$\lim_{j\in J}p_j^{-1}\mathcal O_{X_j}(p_i^{-1}(\Spec A_i))$$ becomes eventually constant and is equal to $(A_i)^{n,K(X)}$\,.\\
  This holds for each $\Spec A_j\subset X_j\longrightarrow X_i,\, A_i\subset A_j$\,. Thus the presheaves defined in (1) and (2) on $R\Spec A_i$\, coincide and the presheaf in (2) is already shown to be a sheaf.  
\end{enumerate}
\end{enumerate} 
\end{proof}
\begin{corollary}
\mylabel{cor:C333} Let $K/k$\, be a function field. Then the Riemann variety $R(K/k)$\, is isomorphic to the projective limit $\mbox{projlim}_{K(X)=K}X$\, over all complete normal birational models of the function field $K/k$\,.\\
In particular, $R(K/k)$\, is isomorphic to the Riemann variety $R(X)$ for any complete model of $K/k$\, and as such quasicompact.
\end{corollary}
\begin{proof} Let $X$ be any complete model of the function field $K/k$\,. Comparing the definitions (2) and (3), we see that $R(K/k)$\, and $R(X)$\, have the same local structure, they are locally isomorphic to locally ringed spaces of the form $(R^K\Spec A, \mathcal O_{R^K\Spec A})$\, for $A$ a finitely generated $k$-algebra. As $X$ is assumed to be separated, it follows that there is an injective morphism of locally ringed spaces 
$$\phi:(R(X),\mathcal O_{R(X)})\hookrightarrow (R(K/k),\mathcal O_{R(K/k)})$$ which is a local isomorphism and sends $R^K\Spec A\subset R(X)$\, to $R^K\Spec A\subset R(K/k)$\,, where $A$ is a finitely generated $k$-algebra lying over $X$\, i.e., there is an open affine subset $\Spec B\subset X$\, such that $B\subset A$\,.\\ Now,  if $\nu\in R(K/k)$\,is any valuation, it follows from the valuative criterion of properness, that $\nu$\, has exactly one center on $X$, say $x\in X$\,. Then $\mathcal O_{X,x}\subset A_{\nu}$\, and there is an open affine subset $\Spec B\subset X$\, such that $$B\subset \mathcal O_{X,x}\subset A_{\nu}$$
($\mathcal O_{X,x}$\, is the localization of the finitely generated $k$-algebra at some prime ideal.) Thus $\phi$\, is also surjective. It follows that $\phi$\, is a bijective local isomorphism and hence an isomorphism.\\
By \ref{prop:P400}, the  definitions (1) and (2) for $R(X)$\, are equivalent and it follows 
$$R(K/k)\cong R(X)\cong \mbox{projlim}_{X'\longrightarrow X}X'$$
the projective limit taken over all $X'$ lying proper and birational over $X$.\\
As the directed system of all $X'$ lying proper and birational over $X$ is cofinal in the directed system of all models of $K/k$\,, it follows $$R(K/k)\cong \mbox{projlim}_{K(X)\cong K}X,$$ the projective limit taken over all models of $K/k$\,. By \ref{thm:Chevalley}, the last quantity is quasicompact. 
\end{proof}
 \begin{lemma}
 \mylabel{lem:L100} Let $X_1$ and $X_2$ be integral preschemes  and let $f: X_1\longrightarrow X_2$\, be a proper surjective morphism. Then the induced map  $R(f): R(X_1)\longrightarrow R(X_2)$\, is a morphism in the category of locally ringed spaces. If $f$ is proper and birational, then $R(f)$\, is the identity morphism.\\
 If $k\subset K\subset L$\, is a tower of field extensions, then the map $$\phi_{LK}: R(L/k)\longrightarrow R(K/k)$$ that sends each valuation $\nu\in R(L/k)$\, to its restriction $\nu\mid_K$\, to the subfield $K$\,, is a surjective morphism of locally ringed spaces.
 \end{lemma}
 \begin{proof} 
 We prove first the second assertion. The surjectivity of the map $\phi_{LK}$\, follows from the extension theorem for valuations with respect to arbitrary field extensions (see \cite{ZaSam}[Vol II, chapter VI, Theorem 11, p.26]). \\
We show that the map $\phi_{LK}$\, is continuous. Let $A\subset K$\, be a finitely generated $k$-algebra such that $R^K\Spec A\subset R(K/k)$\, is a basic open subset.\\
 Then $\phi_{LK}^{-1}(R^K\Spec A)$\, consists of all extensions  of valuations $\nu$\, of $K$ with $A\subset A_{\nu}$\, to valuations $\mu$\, of $L$ with $\mu\mid_K=\nu$\,. Since $A_{\nu}\subset A_{\mu}$\, we have $A\subset A_{\mu}$\,. Conversely, if $A\subset A_{\mu}$\, then, as $A\subset K$\, we also have $A\subset A_{\mu}\cap K=A_{\nu}$\,. Hence $\phi_{LK}^{-1}(R^K\Spec A)=R^L\Spec A$\,and $\phi_{LK}$\, is continuous.\\
 By \ref{prop:P400} we have $$\mathcal O_{R(K/k)}(R\Spec A)=A^{n,K}\,\, \mbox{and}\,\, \mathcal O_{R(L/k)}(R\Spec A)=A^{n,L}.$$
 Let $B\subset L$\, be any finitely generated $k$ algebra. Then $$\phi_{LK}^{-1}\mathcal O_{R(K/k)}(R\Spec B)=\lim_{A\subset B\cap K}A^{n,K}= (B\cap K)^{n,K}$$ and 
\begin{align*}\phi_{LK}^{\sharp}(R\Spec B): \phi_{LK}^{-1}\mathcal O_{R(K/k)}(R\Spec B)=(B\cap K)^{n,K}\\
\longrightarrow \mathcal O_{R(L/k)}(R\Spec B)=B^{n,L}
\end{align*} is the natural inclusion
 $$ (B\cap K)^{n,K}\hookrightarrow (B\cap K)^{n,L}\hookrightarrow B^{n,L}.$$
 Thus $$(\phi_{LK},\phi_{LK}^{\sharp}):(R(L/k),\mathcal O_{R(L/k)})\longrightarrow (R(K/k),\mathcal O_{R(K/k)})$$ is a morphism of locally ringed spaces. \\
We first prove the first assertion in the case where $X_1$\, and $X_2$\, are integral schemes of finite type. The morphism $f$\, induces an inclusion of fields of rational functions $K(X_2)\subset K(X_1)$\,. By the first part of this lemma, there is an induced morphism of Riemann varieties $R(K(X_1)/k)\longrightarrow R(K(X_2)/k)$\,. $R(X_1)\subset R(K(X_1)/k)$\, is the subset of all valuations $\nu\in R(K(X_1)/k)$\, that have center on $X_1$\,, i.e., there exists a scheme point $x_1\in X_1$\, such that $\mathcal O_{X_1,x_1}\subset A_{\nu}$\,. Then $$\mathcal O_{X_2,f(x_1)}\subset \mathcal O_{X_1,x_1}\subset A_{\nu},$$ and $$\mathcal O_{X_2,f(x_1)}\subset A_{\nu}\cap K(X_2)=A_{\nu\mid_{K(X_2)}}.$$
 So  $R(\phi_{K(X_1),K(X_2)})(\nu)=\nu\mid_{K(X_2)}$\, has center on $X_2$\,. Thus the map $R(f)$\, is nothing but the restriction of the morphism $R(\phi_{K(X_1),K(X_2)})$\, to the open subset $R(X_1)\subset R(K(X_1)/k)$\, and so itself a morphism of locally ringed spaces.\\
We turn to the general case. Let $\Spec B\subset X_2$\, be an open affine subset. Since $f$ is proper, $f_B:f^{-1}(\Spec B)\longrightarrow \Spec B$\, is a proper morphism of integral schemes of finite type. As we have just proven, the induced map $R(f_B): R(f^{-1}(\Spec B))\longrightarrow R(\Spec B)$\, is a morphism of locally ringed spaces. Let $$X_2=\bigcup_{i=1}^n\Spec B_i.$$ Then 
$$\forall i=1,...,n\,\,f_{B_i}: f^{-1}(\Spec B_i)\longrightarrow \Spec B_i$$ are morphisms and $X_1=\bigcup_{i=1}^nf^{-1}(\Spec B_i)$\,. Now if $\Spec B\subset \Spec B_i$\, is any open affine subset, then $f_{B}=f_{B_i}\mid_{f^{-1}(\Spec B)}$\, and 
$$R(f_B)=R(f_{B_i})\mid_{R(f^{-1}(\Spec B))}= \phi_{K(X_1),K(X_2)}\mid_{R(f^{-1}(\Spec B))}.$$
We can thus glue the morphisms $R(f_{B_i}), i=1,2,...,n$\, to a morphism $R(f): R(X_1)\longrightarrow R(X_2)$\,.\\
 If $f: X_1\longrightarrow X_2$\, happens to be proper and birational,
 since $$R(X_i)=\mbox{projlim}_{X\longrightarrow X_i}X, i=1,2$$
  and integral preschemes proper and birational over $X_1$\, are final in the directed system of integral preschemes proper and birational over $X_2$\,, it follows that the two projective limits are equal.\\
 So the first part of the lemma is  also proven.
 \end{proof}
 \begin{remark}\mylabel{rem:R500} For each pair of integral complete varieties $X_1,X_2$\, with quotient field $K$, there is a third variety $X$ dominating both $X_1$ and $X_2$\,. By the above lemma, we have $R(X)\cong R(X_1)$\, and $R(X)\cong R(X_2)$\,. It follows $R(X_1)\cong R(X_2)$\,.
 \end{remark}
\begin{proposition}\mylabel{prop:P500} Let $L/k$\, be an arbitrary field extension. Then $$\phi:R(L/k)=\mbox{projlim}_{K\subset L}R(K/k),$$
the projective limit taken over all finitely generated subfields $k\subset K$\, of $L$.
\end{proposition}
\begin{proof} The map $\phi$\, is defined to be the projective limit of the morphisms $$\phi_{LK_i}: R(L/k)\longrightarrow R(K_i/k),\,\, K_i\subset L$$ an arbitrary finitely generated field extension of $k$\,. So $\phi$\, is itself a morphism of locally ringed spaces by \ref{lem:L100} and \ref{lem:L620}. Suppose we are given a system of valuations $\nu_i$\,, one for each finitely generated subfield $K_i\subset L$\, in such a way that $$\mbox{whenever}\,\, K_l\subset K_i,K_j\,\, \mbox{we}\,\, \mbox{have}\,\, \nu_i\mid_{K_l}=\nu_j\mid_{K_l}=\nu_l.$$
Then we may view $(\nu_i)_{i\in I}$\, as an element  
$$(\nu_i)_{i\in I}\in \mbox{projlim}_{K_i\subset L}R(K_i/k)$$
 and each of its elements can be described in this way. The valuations $\nu_i$\, glue to a valuation $\nu$\, of $L/k$\, since $L$ is the union of its finitely generated subfields. We have $\phi(\nu)=(\nu_i)_{i\in I}$\,. Thus $\phi$\, is surjective. \\
 The morphism $\phi$\, is also injective, because if $\phi(\nu)=\phi(\mu)$\, for $\nu,\mu\in R(L/k)$\, then $\nu$\, and $\mu$\, coincide on each finitely generated subfield $K_i$\, of $L$ and therefore also on $L$, since $L=\bigcup_{i\in I}K_i$\,. Thus $\phi$\, is bijective. \\
A basic open subset of the projective limit is given by $p_{i_0}^{-1}(R^{K_{i_0}}\Spec A)$\,, $A$ a finitely generated $k$-subalgebra of some $K_{i_0}$\,. As shown earlier in the proof of \ref{lem:L100}, $$\phi_{LK_{i_0}}^{-1}(R^{K_{i_0}}\Spec A)=R^L\Spec A.$$
 We have by definition
\begin{align*}\mathcal O_{\mbox{projlim}_{K_i\subset L}R(K_i/k)}(R^{K_{i_0}}\Spec A)=\lim_{i\geq i_0}p_i^{-1}\mathcal O_{R(K_i/k)}(R^{K_{i_0}}\Spec A)\\
=\lim_{i\geq i_0}\mathcal O_R(R^{K_i}\Spec A)=\lim_{i\geq i_0}A^{n,K_i}=A^{n,L}.
\end{align*}
Also, $\mathcal O_{R(L/k)}(R\Spec A)=A^{n,L}$\,. The homomorphism $$\phi^{\sharp}: \phi^{-1}\mathcal O_{\mbox{projlim}_{K_i\subset L}R(K_i/k)}(R^{K_{i_0}}\Spec A)\longrightarrow \mathcal O_{R(L/k)}(R^{K_{i_0}}\Spec A)$$ is induced by the inclusions $$A^{n,K_i}=\mathcal O_{R(K_i/k)}(R^{K_i}\Spec A)\longrightarrow \mathcal O_{R(L/k)}(R^L\Spec A)=A^{n,L}.$$ Hence $\phi^{\sharp}$\, is the identity homomorphism and $\phi$\, is an isomorphism.
\end{proof}
\begin{proposition}\mylabel{prop:P401} Let $K/k$\, be an arbitrary field extension. Then the Riemann variety $R(K/k)$\, is equal to the projective limit over all complete integral schemes of finite type $X/k$ with field of rational functions $K(X)\subseteq K$\, and proper  surjective transition morphisms $f: Y\longrightarrow X$\,of finite type, inducing the canonical inclusion $K(X)\subset K(Y)\subset K$\, on the level of rational functions.
\end{proposition}
\begin{proof} If in the projective system of all complete integral schemes of finite type $X/k$\, with quotient field $K(X)\subset K$\, and proper surjective transition morphisms we consider for each subfield $K_i\subset K$\,, finitely generated over $k$, the projective subsystem of all $X$ with quotient field $K_i$\,, we get for each $i_0\in I$\, a map
$$\phi_{i_0}:\mbox{projlim}_{K_i\subset K, K(X_{ij})=K_i}X_{ij}\longrightarrow \mbox{projlim}_{K(X_{i_0j})=K_{i_0}}X_{i_0j}=R(K_{i_0}/k),$$
which is the canonical projection onto that subsystem and by \ref{lem:L620} a morphism of locally ringed spaces. The last equality is the statement of \ref{cor:C333}. \\
If $k\subset K_{i_0}\subset K_{i_1}$\, is a tower of field extensions with $K_{i_0}$ and $K_{i_1}$\, finitely generated over $k$, by \ref{lem:L100}, there is a well defined morphism of locally ringed spaces 
\begin{align*}\phi_{i_1i_0}:=\phi_{K_{i_1}K_{i_0}}: R(K_{i_1}/k)=\mbox{projlim}_{K(X_{i_1j})=K_{i_1}}X_{i_1j}\\
 \longrightarrow \mbox{projlim}_{K(X_{i_0j})=K_{i_0}}X_{i_0j}=R(K_{i_0}/k).
\end{align*}
We have $\phi_{i_0}=\phi_{i_1i_0}\circ \phi_{i_1}$\, because the field extension $K_{i_0}\subset K_{i_1}$\, is finitely generated and each proper morphism $X_{i_1j_1}\longrightarrow X_{i_0j_0}$\, with $K(X_{i_1j_1})=K_{i_1}$\, and $K(X_{i_0,j_0})=K_{i_0}$\, has to be of finite type.
We get 
\begin{align*} \mbox{projlim}_{K_i\subset K, K(X_{ij})=K_i}X_{ij}=\mbox{projlim}_{K_i\subset K}\mbox{projlim}_{K(X_{ij})=K_i}X_{ij}\\
=\mbox{projlim}_{K_i\subset K}R(K_i/k)=R(K/k).
\end{align*}
The last final equality is the contents of \ref{prop:P500}.
\end{proof}  
\begin{proposition}\mylabel{prop:P501} Let $k\subset K\subset L$\, be a tower of field extensions. Then the morphism $\phi_{LK}: R(L/k)\longrightarrow R(K/k)$\, is a closed map, i.e., $\phi_{LK}(\mathcal Z)$\, is closed in $R(K/k)$\, for each closed subset $\mathcal Z\subset R(L/k)$\,.
\end{proposition}
\begin{proof} By \ref{prop:P401} we may write
$$R(L/k)=\mbox{projlim}_{K_i\subset L, K(X_{ij})=K_i}X_{ij},\,\,\mbox{and}\,\, R(K/k)=\mbox{projlim}_{K_i\subset K, K(X_{ij})=K_i}X_{ij}$$ the projective limit taken over all finitely generated subfields $K_i$\, of $L$ and $K$, respectively. The morphism $\phi_{LK}$\, is the projection onto the projective subsystem defining $R(K/k)$\,.\\
Let the closed subset $\mathcal Z\subset R(L/k)$\, be given. Let $Z_{ij}:=(\mathcal Z)_{X_{ij}}$\, be the trace of $\mathcal Z$\, on $X_{ij}$\,, a complete model of $K_i=K(X_{ij})\subset L$\, finitely generated. By \ref{thm:A310}, for each $X_{ij},$\, $Z_{ij}\subset X_{ij}$\, is a Zariski closed subset.\\
 By \ref{cor:C10} we have $$\mathcal Z\cong \mbox{projlim}_{K_i\subset L, K(X_{ij})=K_i}Z_{ij}.$$
  Under this identification, $$\phi_{LK}(\mathcal Z)\cong \mbox{projlim}_{K_i\subset K, K(X_{ij})=K_i}Z_{ij}.$$
 By \ref{cor:C10} $$\mbox{projlim}_{K_i\subset K, K(X_{ij})=K_i}Z_{ij}\cong\bigcap_{K_i\subset K, K(X_{ij})=K_i}\pi_i^{-1}(Z_{ij})$$
 and this is obviously  a Zariski closed subset of $R(K/k)$\,.
 \end{proof}
  
  For any integral prescheme $X$ locally of finite type over $k$ there is a natural map of sets $p_X:R(X)\longrightarrow X$\,. Each $\nu\in R(X)$\, is mapped to $c_X(\nu)\in X$\,.\\  
   Furthermore, there is a natural map of sets $q_X: R(X)\longrightarrow R(K(X))$\, sending each valuation of $R(X)$\, to itself as a point of $R(K(X))$\,.\\
   We have the following 
   \begin{lemma} 
   \mylabel{lem:L101}For each integral prescheme $X$, the maps $q_X: R(X)\longrightarrow R(K(X))$\,  and $p_X: R(X)\longrightarrow X$\,are morphisms of locally ringed spaces. The morphism $q_X$ is a local isomorphism  and the morphism $p_X$\, is universally closed.
   \end{lemma}
   \begin{proof} First, we show the statement for $p_X$\,.\\
    By definition, $R(X)$\, is the projective limit object over all birational models, proper over $X$ in the category of locally ringed spaces. The map $p_X$\, appears to be the projection map onto the factor $X$, which is by definition of the projective limit topology continuous. As $\mathcal O_{R(X)}$\, is defined to be the inductive limit over all $p_Y^{-1}(\mathcal O_Y)$\, the homomorphism $p_X^{\sharp}: p_X^{-1}\mathcal O_X\longrightarrow \mathcal O_{R(X)}$\, is simply the injection into the inductive limit.\\
    If $Y\longrightarrow X$\, is any morphism of preschemes, then $$Y\times_XR(X)=\lim_{K(X_i)=K(X)}Y\times_XX_i$$ is again a projective limit of preschemes with proper transition morphisms and by \ref{thm:A310} the base changed map $p_X\times_XY$\, is closed. Thus $p_X$\, is universally closed.\\
    Now we show the statement for the map $q_X$\,.\\
    Let $R\Spec B\subset R(K(X))$\, be an open subset with $B\subset K(X)$\, normal.  $$\mbox{If}\,\,X=\bigcup_{i\in I} \Spec A_i\,\,\mbox{then}\,\,R(X)=\bigcup_{i\in I}R\Spec A_i.$$ We show that 
$$q_X^{-1}(R\Spec B)\cap R\Spec A_i\,\,\mbox{is open in}\,\, R(X)$$ where from continuity follows. \\
    By \ref{lem:L104} for each $i\in I$\, there exists a complete integral normal model $Y_i$ of $K(X)$\, , open subsets $V_i,U_i\subset Y_i$\, and proper morphisms $p_i:V_i\longrightarrow \Spec B$\, and $q_i:U_i\longrightarrow \Spec A_i$\,. We have $$R\Spec B=R(V_i)\subset R(K(X))\,\,\mbox{and}\,\, R\Spec A_i=R(U_i)\subset R(X).$$
      I claim that the  restriction of $q_X$\, to $R(\Spec A_i)$\,is an isomorphism.\\
 Indeed, $\Spec A_i$\, is a separated integral scheme and a valuation $$\nu\in R(\Spec A_i)\subset R(X)$$ is sent to the valuation $$\nu\in R(\Spec A_i)\subset R(K(X)).$$
 So $q_X$\, restricted to $R(\Spec A_i)\subset R(X)$\, is the identity morphism.\\
      It follows, that $\nu\in R\Spec A_i\subset R(X)$\, is mapped into $R\Spec B\subset R(K(X))$\, precisely when $\nu\in R(U_i)\cap R(V_i)=R(U_i\cap V_i)$\,. Thus $$q_X^{-1}(R\Spec B)\cap R\Spec A_i=R(U_i\cap V_i)$$ which is open in $R(X)$\, by definition of the projective limit topology. Thus $q_X$\,is continuous.\\
Let $R\Spec C\subset R(X)$\, be an open subset, $C\subset K(X)$\, normal. By definition, $$q_X^{-1}\mathcal O_{R(K(X))}(R\Spec C)=\lim_{q_X^{-1}R\Spec B_i\supseteq R\Spec C}\mathcal O_R(K(X))(R\Spec B_i)=\lim_{i\in I}B_i,$$ $B_i\subset K(X)$\, normal. We have to construct a homomorphism of $k$-algebras $$q_X^{\sharp}(R\Spec C):\lim_{i\in I}B_i\longrightarrow C.$$
By \ref{lem:L104}, there are complete normal models $Y_i$\, of $K(X)$\,, open subsets $V_i,U_i\subset Y_i$\, and proper birational morphisms $V_i\longrightarrow \Spec C$\, and $U_i\longrightarrow \Spec B_i$\,. Since $q_X$\, restricted to $R\Spec C$\, is the identity map, $$R\Spec C\subset R\Spec B_i\subset R(K(X)),$$
 hence $V_i\subset U_i$\,. The homomorphism $B_i\longrightarrow C$\, is then simply the sheaf restriction homomorphism
$$B_i=\Gamma(U_i,\mathcal O_{Y_i})\longrightarrow \Gamma(V_i,\mathcal O_{Y_i})=C.$$
 The induced map $\lim_{i\in I}B_i\longrightarrow C$\, is the required map. As $R\Spec C$\, is the final element in the inductive limit, this is the identity homomorphism.    
   \end{proof}
\begin{corollary}\mylabel{cor:C88} Let $K/k$\, be a function field, $\nu\in R(K/k)$\, an arbitrary valuation and $X$ any model of $K/k$\,. Let $p_X(\nu)=\eta$\,. Then $p_X(\overline{\{\nu\}})=\overline{\{\eta\}}$\,.
\end{corollary}
\begin{proof} As $\eta\in p_X(\overline{\{\nu\}})$\, and by \ref{lem:L101} the morphism $p_X$\, is closed, we have $\overline{\{\eta\}}\subseteq p_X(\overline{\{\nu\}})$\,. To show the converse, we have to show, that for each specialization $x$ of $\eta$\,, there is a specialization $\nu\circ \mu, \mu\in R(\kappa(\nu)/k)$\, of $\nu$\, such that $p_X(\nu\circ \mu)=x$\,. Now, there is an inclusion of fields over $k$, $\kappa(\eta)\subset \kappa(\nu)$\, and $\overline{\eta}$\, is a complete model of $\kappa(\eta)$\,. There are canonical continous surjective maps 
$$p_{\overline{\eta}}:R(\kappa(\eta)/k)\twoheadrightarrow  \overline{\eta}\,\,\mbox{and}\,\, \phi_{\kappa(\nu)\kappa(\eta)}: R(\kappa(\nu)/k)\twoheadrightarrow R(\kappa(\eta)/k),$$
where the first map is surjective because over each scheme point $y$ of $X$ there is a valuation with center $y$ on $X$ and the second is surjective by general valuation theory (each valuation $\mu'$\, of $\kappa(\eta)/k$\, can be extended to a valuation $\mu$\, of $\kappa(\nu)/k$). Now choose $\mu'$\, such that $c_{\overline{\eta}}(\mu')=x$\,. Then,
$c_X(\nu\circ \mu)= c_{\overline{\eta}}(\mu')=x$\,.
\end{proof}
   We give the following application of the concept of the Riemann variety.\\
  If $f_{ij}:X_j\longrightarrow X_i$\, is a proper birational morphism of models of a function field $K/k$\,, there is an induced homomorphism of groups of Cartier divisors $$f_{ij}^*:\mbox{CDiv}(X_i)\longrightarrow \mbox{CDiv}(X_j).$$ Thus the system 
\begin{align*}X_i\mapsto\mbox{CDiv}(X_i), X_i\,\,\mbox{a model of}\,\, K(X)\\
              f_{ij}:(X_j\longrightarrow X_i)\mapsto f_{ij}^*: (\mbox{CDiv}(X_i)\longrightarrow \mbox{CDiv}(X_j))
\end{align*}
is a contravariant functor from the category of all models of $K/k$ plus proper birational morphisms to the category of abelian groups and gives an inductive system of abelian groups.\\
  Recall that the abelian group of $b$-Cartier divisors of a function field $K/k$\,, denoted by $\mbox{b-CDiv}(K(X)/k)$\,  is defined to be the inductive limit over that inductive system:
$$\mbox{b-CDiv}(K(X)/k):=\lim_{K(X_i)=K(X)}\mbox{CDiv}(X_i).$$
    Let $$\mbox{CDiv}((R(K/k),\mathcal O_{R(K/k)}))=\mbox{CDiv}(R(K/k)):=\Gamma(R(K/k), \underline{K^*}/\mathcal O_{R(K/k)}^*)$$ denote the group of Cartier divisors on the locally ringed space $(R(K/k),\mathcal O_{R(K/k)})$\,. 
   \begin{proposition}
   \mylabel{prop:P222} Let $K/k$\, be a function field. Then there is an isomorphism of abelian groups $$\alpha:\mbox{b-CDiv}(K/k)\longrightarrow \mbox{CDiv}((R(K/k),\mathcal O_{R(K/k)}))$$ defined by $$\alpha: [(X,D)]\mapsto p_X^*(D),$$
where $D$ is a Cartier divisor on the model $X$, $[(X,D)]$\, indicates the equivalence class of $D$ in the inductive limit and $p_X: R(K(X))\longrightarrow X$\, is as usual the projection morphism.
   \end{proposition} 
\begin{proof} Let $X$ be any model. The map 
\begin{align*}p_X^*: \mbox{CDiv}(X)=\Gamma(X,\underline{ K(X)^*}/\mathcal O_X^*)\\
 \longrightarrow \mbox{CDiv}(R(K(X)/k))=\Gamma(R(K(X)/k), \underline{K(X)^*}/\mathcal O_{R(K(X)/k)}^*) 
\end{align*}
is a homomorphism of abelian groups that comes from the sheaf homomorphism $$p_X^{\sharp}: p_X^{-1}\mathcal O_X^*\longrightarrow \mathcal O_{R(K(X)/k)}^*$$
in the following way: we have a homomorphism of sheaves $$p_X^{\sharp}:p_X^{-1}(\mathcal O_X)\longrightarrow \mathcal O_{R(K(X))}.$$
For $f_{ij}: X_j\longrightarrow X_i$\, a proper birational morphism of models, we have $$f_{ij}^{\sharp}(\mathcal O_{X_i}^*)\subset \mathcal O_{X_j}^*,$$
since units are mapped to units and we get 
$$\mathcal O_{R(K(X))}^*=\lim_{i\in I} p_{X_i}^{-1}O_{X_i}^*,$$ the inductive limit taken over the system of maps $$p_i^{-1}(f_{ij}^{\sharp})\mid_{p_i^{-1}\mathcal O_{X_i}^*}: p_i^{-1}\mathcal O_{X_i}^*\longrightarrow p_j^{-1}\mathcal O_{X_j}^*,$$ since each unit in $\lim_{i\in I}p_{X_i}^{-1}\mathcal O_{X_i}$\, has to lie in some  $(p_{X_i}^{-1}\mathcal O_{X_i})^*=p_{X_i}^{-1}(\mathcal O_{X_i}^*)$\,.\\
From the universal property of the inductive limit, we get a homomorphism of sheaves $p_X^{-1}(\mathcal O_X^*)\longrightarrow \mathcal O_{R(K(X))}^*$\, which induces on global sections 
$$\Gamma (R(K(X)),\underline{K(X)^*}/p_X^{-1}\mathcal O_X^*)\longrightarrow \Gamma(R(K(X)),\underline{K(X)^*}/\mathcal O_{R(K(X))}^*).$$ Finally, we have a homomorphism of abelian groups 
$$\Gamma(X,\underline{K(X)^*}/\mathcal O_X^*)\longrightarrow \Gamma(R(K(X)), \underline{K(X)^*}/p_X^{-1}\mathcal O_X^*)$$ which sends the cycle $(f_i)_{i\in I}$\, subordinate to the open covering $X=\bigcup_{i\in I}\Spec A_i$\, to the cycle $(f_i)_{i\in I}$\, subordinate to the open covering $R(K(X))=\bigcup_{i\in I}R\Spec A_i$\,. The composition of both morphisms is equal to $p_X^*$\,. Thus each $p_X^*$\, is a homomorphism of abelian groups. For each proper birational morphism $f_{ij}:X_j\longrightarrow X_i$\, we have $f_{ij}\circ p_{X_j}=p_{X_i}$\,. We get a homomorphism of abelian groups 
$$ \mbox{b-CDiv}(K(X)/k)=\lim_{i\in I} \mbox{CDiv}(X_i)\stackrel{\lim p_{X_i}^*}\longrightarrow \mbox{CDiv}(R(K(X)/k)).$$
This is nothing else than the homomorphism $\alpha$, for, if $D\in \mbox{CDiv}(X)$\, corresponds to the cycle $(f_i)_{i\in I}$\, subordinate to the open covering $X=\bigcup_{i\in I}\Spec A_i$\,, then $p_X^*(D)$\, corresponds to the cycle $(f_i)_{i\in I}$\, subordinate to the open covering $R(K(X)/k)=\bigcup_{i\in I}R\Spec A_i$\, which is the above constructed homomorphism. It follows that $\alpha$\, is a homomorphism of abelian groups.\\
We show that $\alpha$\, is injective. Let $[(X,D)]\in \mbox{b-CDiv}(K/k)$\, be a $b$-Cartier divisor with $\alpha([(X,D)])=0$\,. Let $D$ correspond to the cycle $(f_j),j=1,...,n$\, subordinate to the open covering $X=\bigcup_{j=1}^n\Spec B_j$\,. Then there is an open covering $$R(K/k)=\bigcup_{i\in I}R\Spec A_i,$$ finer than the covering $R(K(X)/k)=\bigcup_{j=1}^nR\Spec B_j$\,with, say  $B_j\subset A_{i_j}, i_j\in I$\, such that the cycle $((f_j,R\Spec A_{i_j}))_{j\in I}$\, is trivial.\\ Because of the theorem of Chevalley we may assume that the index set $I$ is finite, say $I=\{1,...,m\}$\,. By \ref{lem:L104} there is a complete normal model $Y$ of $K/k$\, , open subsets $U_i\subset Y, i=1,...,m$\, plus proper birational morphisms $$p_i:U_i\longrightarrow \Spec A_i, i=1,...,m$$ plus a proper morphism $q: Y\longrightarrow X$\,. Since the $R\Spec A_i$\, cover $R(K/k)$\, and $p_i$\, is proper, the $U_i,i=1,...,m$\, cover $Y$. The covering $Y=\bigcup_{i=1}^mU_i$\, is finer than the covering $$Y=\bigcup_{j=1}^nq^{-1}(\Spec B_i).$$ By assumption, the cycle $((f_{i_j},U_j))_{j=1}^m$\, is trivial and represents $q^*(D)$\,. But the homomorphism $q^*: \mbox{CDiv}(X)\longrightarrow \mbox{CDiv}(Y)$\, is injective as $q$ is proper and birational (observe that we are at the cycle level and not on the level of rational equivalence classes). Thus $D=0$\, and $\alpha$\, is injective.\\
Now we show surjectivity of the homomorphism $\alpha$\,. A Cartier divisor $\mathcal D$\, on the locally ringed space $(R(K(X)/k),\mathcal O_{R(K(X)/k)})$\, is given by a Zariski open covering $$R(K/k)=\bigcup_{i\in I}\mathcal U_i$$
 and a collection of rational functions $$f_i\in \Gamma(\mathcal U_i,\underline{K(X)}), i\in I$$ such that $$\frac{f_i}{f_j}\in \Gamma(\mathcal U_i\cap \mathcal U_j, \mathcal O_{R(K/k)}^*),$$ by definition of the quotient sheaf $\underline{K(X)}^*/\mathcal O_R^*$\,. Since open subsets of the form $R^K\Spec B, B\subset K$\, a finitely generated $k$-algebra, form a basis for the topology for $R(K/k)$\,, we may assume that there are finitely generated normal $k$-algebras $A_i\subset K$\, such that $\mathcal U_i=R\Spec A_i$\,. By the theorem of Chevalley, there are finitely many indices $i_1,...,i_n$\, such that $$R(K/k)=\bigcup_{j=1}^nR\Spec A_{i_j}.$$
By \ref{lem:L104}, there is a complete normal model $Y$ of $K/k$, open subsets\\ $U_j\subset Y,j=1,...,n$\, plus proper birational morphisms $$p_j:U_j\longrightarrow \Spec A_{i_j},j=1,...,n.$$ By our well known arguement, we have $R\Spec A_{i_j}=R(U_j), j=1,...,n$\, and $$R\Spec A_{i_j}\cap R\Spec A_{i_l}=R(U_j\cap U_l), 1\leq j<l\leq n.$$
We have $\mathcal O_{R(K/k)}(R(U_j\cap U_l))=\mathcal O_Y(U_j\cap U_l)$\, by \ref{prop:P400} and thus 
$$\mathcal O_{R(K/k)}(R(U_j\cap U_l))^*=\mathcal O_Y(U_j\cap U_l)^*.$$
So for $1\leq j<l\leq n$\, we have $\frac{f_{i_j}}{f_{i_l}}\in \mathcal O_Y(U_j\cap U_l)^*.$\\
Therefore, the collection $(f_{i_j}), j=1,...,n$\, can be considered as a section  $$D\in\Gamma(Y, \underline{K(Y)^*}/\mathcal O_Y^*)=\mbox{CDiv}(Y).$$ Thus we have $\alpha([(Y,D)])=\mathcal D$\,. Since $\mathcal D$\, was arbitrary, $\alpha$\, is surjective.  
\end{proof}

\begin{definition}
\mylabel{def:D102} An abstract Riemann variety is a locally ringed space $(\mathcal Z,\mathcal O_{\mathcal Z})$ which has an open cover by locally ringed spaces of the form $(R(\Spec A_i),\mathcal O_{R(\Spec A_i)})$\,, where each $A_i$\, is an integral domain of finite type with the same quotient field $K$.
\end{definition}
One might wonder whether for each abstract Riemann variety $(\mathcal Z,\mathcal O_{\mathcal Z})$\, there is an integral prescheme $X$ locally of finite type such that $(\mathcal Z,\mathcal O_{\mathcal Z})\cong R(X)$\,.\\
That this is not the case is shown by the following
\begin{example}
\mylabel{ex:E200} Let $K/k$\, be a function field of transcendence degree larger than one and $\nu\in R(K/k)$\, be a closed point. Then define $(\mathcal Z,\mathcal O_{\mathcal Z})$\, by glueing two copies of $R(K/k)$\, along $R(K/k)\backslash\{\nu\}$\,. \\
Suppose there was an integral prescheme $X$ locally of finite type with $R(X)=\mathcal Z$\,.\\ Then there are open affine subsets $\Spec A_1,\Spec A_2\subset X$\, such that one copy $\nu_1$ of $\nu$\, is in $R\Spec A_1$\, and the other copy $\nu_2$ of $\nu$\, is in $R\Spec A_2$\,.\\
 By \ref{lem:L104} there is a complete integral scheme $Y$ plus open sets $U,V\subset Y$\, plus proper morphism $p_1:U\longrightarrow \Spec A_1$\, and $p_2:V\longrightarrow \Spec A_2$\,.\\
  The unique lifts of the centers $x_i$ of $\nu_i$\, on $\Spec A_i$\, via $p_i, i=1,2$\, to $Y$ must coincide and be equal to the center $y$\, of $\nu$\, on $Y$\, since $Y$ was assumed to be proper. As $\nu$\, is a closed point in $R(K/k)$\,, $y$\, is a closed point of $Y$\,. Hence $\mathcal O_{Y,y}$\, dominates the local rings $\mathcal O_{\Spec A_i,x_i},i=1,2$\,.\\
  Now let $\mu\neq \nu$\, be any valuation with center $y$ on $Y$. $\mu\neq \nu$\, always exists because $\mathcal O_{Y,y}$\, is essentially of finite type and not a valuation ring. Then it has centers $x_i$\, on $\Spec A_i, i=1,2$\,. Since $x_1\neq x_2$\,, there is one copy $\mu_1$\, of $\mu$\, in $R\Spec A_1\subset R(X)$\, and a second copy $\mu_2$\, of $\mu$\, contained in $R\Spec A_2\subset R(X)$\,. But obviously, $\mathcal Z\backslash\{\nu_1,\nu_2\}$\, is separated, so it cannot contain two copies of the same valuation, a contradiction.
  \end{example} 
  If we refer to the space of discrete algebraic rank one valuations (that is, valuations that correspond to prime divisors on some blow-up), we will use the notation $R^1(K(X))$\, instead of the correct notion $R^{1,1}_{Ab}(K(X))$\,. \\
  We prove the following 
 
\begin{lemma}\mylabel{lem:L1000} Let $K/k$\, be an algebraic function field of characteristic zero with $\mbox{trdeg}(K/k)\geq 2$\,. Then $R^1(K/k)$\, is not quasicompact.
\end{lemma}
\begin{proof} Start with a smooth complete model $X$ of $K/k$, remove one point $P$, let $U_1:=X\backslash\{P\}$\,, and inductively, if one has defined $X_n$, $P_n$\, and  $U_n:=X_n\backslash\{P_n\}$, then one defines $X_{n+1}:=\mbox{Bl}_{P_n}X_n$\,, chooses a point $P_{n+1}$\, on the exceptional divisor of $X_{n+1}$\, and defines $U_{n+1}:=X_{n+1}\backslash\{P_{n+1}\}$\,. Then we claim that $\bigcup_nR^1(U_n)=R^1(K(X))$\, but this open cover has no finite subcover since $U_n\subsetneq U_{n+1}$\,. Note that $\bigcup_nR(U_n)\neq R(K(X))$\,, because exactly the valuations are missing whose valuation rings dominate  the direct limit of the local rings of the points $P_n$\, on $X_n$\,. These always exist because there is a maximal local ring dominating the local ring $$\lim_{n\in \mathbb N}\mathcal O_{X_n,P_n}=\bigcup_{n\in \mathbb N}\mathcal O_{X_n,P_n}\subset K(X).$$
We have to show that there is no discrete algebraic rank one valuation having center $P_n$\, on each $X_n$.\\
Recall that the discrepancy $a(X,\nu_F)$\, of a variety $X$ at a discrete algebraic rank one valuation $\nu_F$\, is defined by choosing any rational top differential form $\omega\in \Lambda^n(\Omega(K(X)/k))$\, and putting 
$$a(X,\nu_F):=\nu_F(\omega)-\nu_F(K_X^{\omega}),$$ where $K_X^{\omega}$\, is the rational section of $\Lambda^n(\Omega(\mathcal O_X/k))$\, corresponding to $\omega$\, and $\nu_F(\omega)$\, is calculated by considering $\omega$\, as a rational section of $\Lambda^n(\Omega(A_{\nu_F}/k))$\, and taking the usual valuation of a Cartier divisor. For closer information about discrepancies see \cite{Matsuki}[chapter 4, p. 163-191].\\ Suppose there is such a valuation $\nu_E$\,. As every $X_n$\, is smooth, it is known that $X_n$ has so called  terminal singularities, which by definition means that for all  $\nu\in R^1(X)$\, which are exceptional over $X$ one has $a(X_n,\nu)>0$\, (see \cite{Matsuki}[chapter 4, section 4.1,Theorem 4-1-3, p.166]). If one considers the blow-up $p_n:X_{n+1}\longrightarrow X_n$, then $$K^{\omega}_{X_{n+1}}-p_n^*K^{\omega}_{X_n}= E_n$$ (the canonical divisors $K_{X_{n+1}}$ and $K_{X_n}$\,are taken with respect to the same rational differential $n$-form $\omega$). As $\nu_E$\, has center on $E_n$\, (namely the point $P_n$), we have $\nu_E(E_n)>0$\, because the local equation for $E_n$ $(t_n=0)$\, lies in $\mathfrak{m}_{X_n,P_n}\subset\mathcal O_{X_n,P_n}$\, which is contained in the maximal ideal $\mathfrak{m}_{\nu_E}\subset A_{\nu_{E}}$\,. This implies that the  discrepancy $$a_n:=a(X_n,\nu_E)=\nu_E(\omega)-\nu_E(K^{\omega}_{X_n}) >0 $$   strictly diminishes with growing $n$. So the $(a_n)_{n\in \mathbb N}$ form an infinite  strictly monotonously decreasing sequence of positive integers and this is impossible.\\
So we conclude $R^1(K(X))\subset \bigcup_{n\in \mathbb N}R^1(U_n)$\, is an infinite open covering which has no finite subcovering because of $R^1(U_n)\subsetneq R^1(U_{n+1})$\,.
\end{proof}
\begin{remark}\mylabel{rem:R21}
This example also shows that, unlike in the case of finite type schemes, not every open subset of $R(X)$\, is quasicompact. One has to take $$\mathcal U:=\bigcup_{n}R(U_n).$$ The open covering of $\mathcal U$\, by the $R(U_n)$\, has obviously no finite subcovering.
\end{remark}
A scheme of finite type over a finite or countable base field has only countably many subschemes which may be seen as follows. The affine coordinate rings are countable and there are only countably many ideals in a countable noetherian ring as every ideal is generated by finitely many elements and there are only countably many possibilities to choose finitely many elements out of a countable set. As $X$ has a finite cover by affine open subsets $\Spec A_i, i=1,...,n$\,, the collection of all subschemes of $X$ which corresponds to the collection of ideal sheaves $(\mathcal I_j\subset \mathcal O_X)_{j\in I}$\, can be injectively mapped to 
$$\prod_{i=1}^n(\mbox{ideals of}\Spec A_i)$$ and this is a finite product of countable sets and thus is at most countable.\\
On the contrary we prove the following
\begin{lemma}\mylabel{lem:L22} If one is given an integral complete scheme $X$ over a finite field $k$ with $\dim (X)>1$, then the set of points of $R(X/k)$\, is uncountable. 
\end{lemma}
\begin{proof} Choose any closed point $P$ on $X$ . As in the example in the course of the proof of \ref{lem:L1000}, the set of sequences of points $(P_n)$\,, starting with $P=P_0$\, is uncountable, because every exceptional divisor is of dimension at least one  and so has  at least $2$ closed points.  which we can take as $P_n$\,. Thus the set of valuations over $P_0$\, contains a subset that is mapped surjectively to the set of  arbitrary sequences of zeroes and ones, i.e., arbitrary binary numbers and this set is uncountable.
\end{proof}
The Riemann variety $R(K(X))$\, is a Zariski topological space, i.e., every irreducible closed subset has a generic point. The closure of $\{\nu\}$\,, where $\nu$\, is a valuation, is the set of all valuations composed with $\nu$\, (see \cite{Val}[Theorem 2.6]). The closure  $\overline{\{\nu\}}$\, is isomorphic to the Riemann variety $R(\kappa(\nu))$\,, where $\kappa(\nu)$\, denotes the residue field of $\nu$\,. The generic point of the whole Riemann variety is the point corresponding to the trivial valuation of $K(X)$\,.\\
\begin{proposition}\mylabel{prop:P888} Let $k\subset K\subset L$\, be a tower of finitely generated field extensions and $\nu\in R(L/k)$\,. The Zariski topological dimension of the irreducible closed subset $\overline{\{\nu\}}$\, equals the dimension of $\nu$\,, i.e., the transcendence degree $\mbox{trdeg}(\kappa(\nu)/k)$\,.\\ We have $\dim(\phi_{LK}(\overline{\{\nu\}}))\leq\dim(\overline{\{\nu\}})$\,.\\
If $Y$ is any model of $L$  and $\mathcal Z$\, a proper Zariski closed subset of $R(L/k)$\,, then $p_Y(\mathcal Z)$\, is a proper Zariski closed subset of $Y$.
\end{proposition}
\begin{proof} Let  $\dim(\overline{\{\nu\}})=d$\,. Then there is a chain $$\nu=\nu_1,\nu_2,...,\nu_{d+1}$$ of valuations of $R(L/k)$\, with $\nu_{i+1}\in \overline{\{\nu_i\}},i=1,...,d$\,. This means precisely that there are valuations $\mu_i\in \kappa(\nu_i)$\, such that $\nu_{i+1}=\nu_i\circ \mu_i, i=1,...,d$\,. We have $\kappa(\mu_i)=\kappa(\nu_{i+1})$\,, which is the residue field of the valuation $\nu_{i+1}$\,. It follows $$\mbox{trdeg}(\kappa(\nu_{i}))\geq \mbox{trdeg}(\kappa(\nu_{i+1}))+1\,\, \mbox{and} \,\,\mbox{trdeg}(\kappa(\nu_1))=\mbox{trdeg}(\kappa(\nu))\geq d.$$
 Suppose the transcendence degree of $\kappa(\nu)$\, is srictly larger than $d$, say $d+l,l>0$\,. Let $x_1,...,x_{d+l}$\, form a transcendence basis of $\kappa(\nu)$\,. We consider affine\\
 $(d+l)$-space $\mathbb A^{d+l}_k$\, with field of rational functions $k(x_1,...,x_{d+l})$\,. Let 
$$H_1:(x_1=0),H_2:(x_2=0),...,H_{d+l}:(x_{d+l}=0)$$
 be the coordinate hyperplanes. Put $$D_k:=\cap_{i=1}^kH_i, k=1,...,d+l.$$ By \ref{ex:E61} and \ref{ex:E62}, the incomplete flag $$D_1\supseteq D_2... \supseteq D_k,k\leq d+l$$
 gives rise to a valuation $\mu_k$\, of $k(x_1,...,x_{d+l}), k=1,...,d+l$\, and we have \\ $\mu_{k+1}\in \overline{\{\mu_k\}}$\,.\\
 We construct inductively a chain of valuations $\nu_1,...,\nu_{d+l}$\, of $\kappa(\nu)$\, with $$\nu_{k+1}\in \overline{\{\nu_k\}}\,\,\mbox{and} \,\,\phi_{\kappa(\nu)k(x_1,...,x_{d+l})}(\nu_k)=\mu_k, k=1,...,d+l.$$ Since by the extension theorem of valuations the morphism $\phi_{\kappa(\nu)k(x_1,...,x_{d+l})}$\, is surjective, we find $\nu_1\in R(\kappa(\nu))$\, with image $\mu_1$\, in $R(k(x_1,...,x_{d+l}))$\,. Inductively, suppose we have constructed a chain $\nu_1,...,\nu_k, k<d+l$\, with 
$$\phi_{\kappa(\nu)k(x_1,...,x_{d+l})}(\nu_i)=\mu_i, i=1,...,k\,\, \mbox{and} \,\,\nu_{i+1}\in \overline{\{\nu_i\}}, i=1,...,k-1.$$ Let $\mu_{k+1}=\mu_k\circ \omega_k, \omega_k\in R(\kappa(\mu_k))$\,. Again by the extension theorem of valuations, we can find $\xi_k\in R(\kappa(\nu_k))$\, with $\phi_{\kappa(\nu_k)\kappa(\mu_k)}(\xi_k)=\omega_k$\,. Put $\nu_{k+1}:=\nu_k\circ\xi_k$\,. Then $$\nu_{k+1}\in \overline{\{\nu_k\}}\,\, \mbox{and} \,\,\phi_{\kappa(\nu)k(x_1,...,x_{d+l})}(\nu_{k+1})=\mu_{k+1}.$$
 Considering the chain $\nu, \nu_1,...,\nu_{d+l}$\,, it follows $\dim(\overline{\{\nu\}})\geq d+l$\,, a contradiction. This proves the first assertion.\\
The second assertion follows easily from the first: a valuation $\nu$\, of $R(L/k)$\, is mapped to the restriction $\nu\mid_K$\, and we have an inclusion of residue fields $$\kappa(\nu\mid_K)\subseteq \kappa(\nu).$$ So it follows
$$ \dim(\phi_{LK}(\overline{\{\nu\}}))=\dim(\overline{\{\nu\mid_K\}})=\mbox{trdeg}(\kappa(\nu\mid_K))\leq \mbox{trdeg}(\kappa(\nu))=\dim(\overline{\{\nu\}}).$$
To prove the last statement, by \ref{lem:L101}, $p_Y(\mathcal Z)$\, is closed. The Zariski topological dimension of $\mathcal Z$\, is defined to be $$\dim(\mathcal Z)=\sup_{\nu\in \mathcal Z}\dim(\overline{\{\nu\}}).$$
If $\eta\in p_Y(\mathcal Z)$\, there is $\nu\in \mathcal Z$\, with $p_Y(\nu)=\eta$\,. Then $\kappa(\eta)\subseteq \kappa(\nu)$\, and by the first part 
\begin{align*}\dim(Y)=\dim(R(L/k))>\dim(\mathcal Z)\geq \dim(\overline{\{\nu\}})\\
=\mbox{trdeg}(\kappa(\nu)/k)\geq \mbox{trdeg}(\kappa(\eta)/k)=\dim(\overline{\{\eta\}}),
\end{align*}
so the generic point of $Y$ cannot be in $p_Y(\mathcal Z)$\, and $p_Y(\mathcal Z)$\, has to be a proper Zariski closed subset of $Y$.
\end{proof}  
Observe that unlike in the case of noetherian schemes not every closed subset is a finite union of irreducible ones (see \ref{lem:L41}).\\
But we have the following lemma characterizing Zariski closed subsets of $R(K(X))$\, being a finite union of irreducible ones.
\begin{lemma}\mylabel{lem:L41} A Zariski closed subset $\mathcal Z\subset R(K(X))$\, is a finite union of irreducible closed subsets if and only if there is a complete model $Y$ such that for all proper and birational morphisms $Y'\longrightarrow Y$\, the induced morphism $\mathcal Z_{Y'}\longrightarrow \mathcal Z_{Y}$\, is generically finite and the number of irreducible components of $\mathcal Z_{Y'}$\, is bounded above by a natural number $n$.
\end{lemma}
\begin{proof} As by \cite{Val}[Theorem 2.6], each irreducible closed subset has a generic point, if $\mathcal Z$\, is a finite union of irreducible closed subsets $\mathcal Z_i$\, and $\eta_i$\, is the generic point of $\mathcal Z_i$\,, according to \cite{Val}[Proposition 2.3], for each $i$ there is a complete model $Y_i$\, such that if $y_i$\, is the center of $\eta_i$\, on $Y_i$\,, the residue field extension $\kappa(y_i)\subset \kappa(\eta_i)$\, is algebraic.\\ Choosing $Y$ dominating all $Y_i$\,, if $Y'\longrightarrow Y$\, is a proper birational morphism and $y'_i$\, is the center of $\eta_i$\, on $Y'$, one has an inclusion of residue fields $$\kappa(y_i)\subset \kappa(y'_i)\subset \kappa(\eta_i),$$ hence $\kappa(y_i)\subset \kappa(y_i')$\, is algebraic (and finite).\\
 It follows that $y_i$ and $y_i'$\, are the generic points of $\mathcal Z_{i,Y}$\, and $\mathcal Z_{i,Y'}$\,, respectively. The morphism $$\mathcal Z_{i,Y'}\longrightarrow \mathcal Z_{i,Y}$$ corresponds to a finite extension of quotient fields, hence it is generically finite.\\ The number of irreducible components of $\mathcal Z_{Y'}$\, is bounded by the number of irreducible components of $\mathcal Z$\,.\\
To prove the converse, assume that the complete model $Y$\, satisfies the assumptions of the lemma. Choose $Y$ furthermore such that the number of irreducible components of $\mathcal Z_Y$\, takes on its maximal value. Let $y_i, i=1,...,n$\, be the generic points of the (finitely many) irreducible components of $\mathcal Z_Y$\,.  For each proper birational morphism $Y'\longrightarrow Y$\,, the fibres over the generic points $y_i$\, of $\mathcal Z_Y$\,
$$\mathcal Z_{Y'}\times_Y\Spec \kappa(y_i)$$
 are finite and nonempty, because $\mathcal Z_{Y'}\longrightarrow \mathcal Z_{Y}$\, was assumed to be generically finite and is by definition surjective. Over each $y_i$\, must lie a point $y_i'$\,. This therefore must be the generic point of some irreducible component of $\mathcal Z_{Y'}$. But the number of  the irreducible components of $\mathcal Z_Y$\, was assumed to be maximal, and therefore over each point $y_i$\, there is exactly one point $y_i'$\, in $$\mathcal Z_{Y'}\times_Y\Spec \kappa(y_i).$$
  Hence in the limit, there is exactly one point $\eta_i\in \mathcal Z$\, lying over $y_i$\,.\\ Consider the Zariski closed subset $\mathcal Z':=\cup_i\overline{\{\eta_i\}}$\,. For all proper birational morphisms\\
 $Y'\longrightarrow Y$\, we have by \ref{cor:C88} by what was said above \\
$$\mathcal Z'_{Y'}=\mathcal Z_{Y'} = \bigcup_i\overline{\{y_i'\}}.$$
 By \ref{thm:A310},we have
 $$\mathcal Z=\bigcap_{Y'\longrightarrow Y}p_{Y'}^{-1}(\mathcal Z_{Y'})=\bigcap_{Y'\longrightarrow Y}p_{Y'}^{-1}(\mathcal Z_{Y'}')=\mathcal Z'.$$ 
Thus $\mathcal Z$\, is a finite union of irreducible closed subsets.
\end{proof}  
\begin{remark}
\mylabel{rem:R101} The most simple examples for which the assumptions of the above lemma are not satisfied are closed subsets of the form $\mathcal Z=p_Y^{-1}(A)$\,, $A\subset Y$\, a closed subset . For each blow-up sequence $Y''\stackrel{q}\longrightarrow Y'\stackrel{p}\longrightarrow Y$\, $q$\, with center in $p^{-1}(A)$\, of codimension greater than one, the exceptional divisor of $q$ belongs to $\mathcal Z_{Y''}$\, but does not lie generically finite over $\mathcal Z_{Y'}$\,.\\
So closed subsets of the form 
$$\mathcal Z=p_Y^{-1}(A), A\subsetneq Y\,\,\mbox{Zariski}\,\,\mbox{closed}$$
are not finite unions of irreducible closed subsets.
\end{remark}
\section{The structure of integral preschemes}
All we know about preschemes which are not schemes is that they are not separated. We will investigate in this section the closer structure of nonseparatedness, or "the degree of nonseparatedness" at particular points of an integral prescheme. If $x\in X$\, and we can find a valuation $\nu\in R(K(X)/k)$\, with $n$ distinct centers $x_1=x, x_2,...,x_n$\, we could say that the degree of nonseparatedness of the prescheme $X$ at the point $x$ is at least $n$. We will make sense of this heuristic idea and introduce a stratification $(X_n)_{n\in \mathbb N}$\, on each integral prescheme $X$ which we show to consist of constructible sets if $X$ is of finite type and satisfies the existence condition of the valuative criterion of properness. At each $x\in X_n$\, the "degree of nonseparatedness" will then be exactly $n$\,.\\
We start by giving a structure theorem for integral preschemes locally of finite type with quotient field $K$.\\
Let $A_i,A_j\subset K$\, be two finitely generated $k$-algebras with quotient field $K$. We have the associated Riemann varieties $R\Spec A_i,R\Spec A_j\subset R(K)$\,. We define the Zariski open subset  $$\mathcal U_{ij}:=R\Spec A_i\cap R\Spec A_j$$ of $R(K(X))$\,.\\
 We prove the following
\begin{lemma}
\mylabel{lem:L112} With notation as above, if there is a birational equivalence 
$$\phi:\Spec A_i\supseteq U_i\cong U_j\subseteq \Spec A_j,$$
 then $U_i=U_j=:U$\, and $p^{-1}(U)\subseteq \mathcal U_{ij}$\,. There exists a largest such open subset which we call $U_{ij}$\, and this is the largest open subset contained in $p(\mathcal U_{ij})$\,. 
\end{lemma}
\begin{proof} Since $U_i$ and $U_j$\, have identical quotient field $K,$\, the structure sheaves $\mathcal O_{U_i}$\, and $\mathcal O_{U_j}$\, are subsheaves of the constant sheaf $\underline{K}$\,. Namely, if $\Spec B_j\subset U_j$\, is an affine open subset, then $\phi^{-1}(\Spec B_j)=:\Spec B_i$\, is again an affine open subset of $U_i$\,. The isomorphism $\phi^{\sharp}: B_j\longrightarrow B_i$\, induces the identity on quotient fields and since $B_i,B_j$\, are integral, $\phi^{\sharp}$\, must also be the identity. Hence from $R(U_i)\subset R\Spec A_i$ and $R(U_j)\subset R\Spec A_j$\, it follows that 
$$R(U)=R(U_i)=R(U_j)\subseteq \mathcal U_{ij}.$$
 As every birational equivalence of integral schemes with the same quotient field is the identity on some open subset, each such collection of birational equivalences glues to a birational equivalence on the union. By Zorn's Lemma (or simply by quasicompactness) there exists a largest such subset $U_{ij}$ with $p^{-1}(U_{ij})\subset \mathcal U_{ij}$\, or $U_{ij}\subseteq p(\mathcal U_{ij})$\,.
\end{proof}
\begin{remark}
\mylabel{rem:R112} The subset $p(\mathcal U_{ij})\subset \Spec A_i,\Spec A_j$\, is not itself necessarily open.
 By \ref{lem:L104} there is a complete normal model $Y$ of $K/k$\, and open subsets $V_i,V_j\subset Y$\, and proper birational morphisms $$q_i:V_i\longrightarrow \Spec A_i\,\, \mbox{and}\,\, q_j: V_j\longrightarrow \Spec A_j.$$
 By the properness of $q_i,q_j$\, and \ref{lem:L100} we have $$R\Spec A_i=R(V_i)\,\, \mbox{and}\,\, R\Spec A_j=R(V_j).$$
 It follows $\mathcal U_{ij}= R(V_i\cap V_j)$\,.
  The maps $$p_i: \mathcal U_{ij}\longrightarrow \Spec A_i, p_j:\mathcal U_{ij}\longrightarrow \Spec A_j$$ factor over $Y$, so $$p_i(\mathcal U_{ij})=q_i(V_i\cap V_j)\,\, \mbox{and}\,\, p_j(\mathcal U_{ji})= q_j(V_i\cap V_j).$$ As the proper morphisms $q_i$ and $q_j$ are closed but not necessarily open, $p_i(\mathcal U_{ij}),p_j(\mathcal U_{ij})$\, are in general only  constructible subsets of $\Spec A_i,\Spec A_j$\,,respectively, by \cite{Ha}[chapter II, Exercise 3.19, p.94] .
\end{remark}
\begin{proposition}
\mylabel{prop:P100} Let $K/k$\, be a function field and let $A_i\subset K, i\in I$\, be a collection of finitely generated $k$-algebras with quotient field $K$\,. With notation as above, let furthermore for each pair $(i,j)\in I\times I$\,  a Zariski open subset $U_{ij}'\subset U_{ij}$\, be given.\\
Assume that for each pair $(i,j)\in I\times I$\, there is a $k\in I$\, such that $U_{ik}'\cap U_{jk}'\neq \emptyset$\,(*). Then there is an integral prescheme $X$ locally of finite type with $K(X)=K$\, such that
\begin{enumerate}[1]
\item $\Spec A_i  \subset X$\, is a Zariski open subset for all $i\in I$\,,
\item $\Spec A_i\cap \Spec A_j=U_{ij}'$\,, the intersection taken inside $X$ and
\item $X=\bigcup_{i\in I}\Spec A_i$\,.
\end{enumerate}
Conversely, if $X$ is an integral prescheme locally of finite type with quotient field $K$ and $X=\bigcup_{i\in I}\Spec A_i$\, is a Zariski open covering, if we define $$U_{ij}'=\Spec A_i\cap \Spec A_j,$$ then the data $$A_i\subset K, i\in I\,\, \mbox{and}\,\, U_{ij}', (i,j)\in I\times I$$ define by the first part an integral prescheme $X'$\, with $X'\cong X$.\\ 
\end{proposition}
\begin{proof} The data $$\Spec A_i, i\in I, U_{ij}', (i,j)\in I\times I\,\,\mbox{and}\,\, \phi_{ij}:U_{ij}\cong U_{ji}$$ define gluing data that satisfy trivially the cocycle condition because $\phi_{ij}$\, are the identity morphisms. This is because all sheaves are subsheaves of the fixed function field $K$. The condition (*) ensures that the prescheme $X$ in connected. For, if $X_i,X_j$\, are two irreducible components of $X$ and $\Spec A_i\subset X_i$\, and $\Spec A_j\subset X_j$\, are two Zariski open subsets, by (*) there is an index $k$\, such that
 $$\emptyset \neq U_{ik}'\cap U_{jk}'\subset \Spec A_i\cap \Spec A_j\subset X_i\cap X_j.$$
 Because each $A_i$ is an integral domain, $X$ is integral.\\
Conversely, the gluing datum extracted out of the the open covering of $X$ satisfies $U_{ij}'\subset U_{ij}$, since $$R(U_{ij}')\subset R\Spec A_i\cap R\Spec A_j\,\, \mbox{and}$$
$$\forall i,j,k: \emptyset \neq \Spec A_i\cap_X\Spec A_j\cap_X\Spec A_k=U_{ik}'\cap_XU_{jk}'.$$
 Thus $X$ and the prescheme $X'$\, from the first part are constructed out of the same gluing datum and consequently $X'\cong X$\,.
\end{proof}
We carry on by  giving proofs of some versions of the valuative criterion of properness.
\begin{lemma}\mylabel{lem:L2062} Let $X$ be an integral prescheme such that  the natural morphism $q_X: R(X)\longrightarrow R(K(X))$\, is an isomorphism. Then $X$ is a complete scheme which is a variety if $X$ is assumed to be locally of finite type.\\
 If $q_X$ is assumed only to be injective, then $X$ is a separated scheme. 
\end{lemma}
\begin{remark}
\mylabel{rem:R102} This lemma can be seen as the statement of the fact, that in order to verify the valuative criterion of properness and separatedness for an integral prescheme $X$, it suffices to consider valuations with quotient field $K(X)$\,. This is actually well known. See \cite{EGA}[EGA I, chapter 5, Proposition (5.5.4); EGA II, chapter 7, $\rm{Th\acute{ e}or\grave{e}me}$ 7.8.3].
\end{remark}
\begin{proof} We want to verify the valuative criterion of properness. Say, we are given a morphism $x: \Spec L\longrightarrow X$\, and $\nu$ a valuation of $L/k$\,. If $L=K(X)$, in this case the valuative criteria of properness and separatedness follow immediately from the bijectivity or injectivity of the map $q_X$, respectively.\\
 In general, $x$ factors as $$\Spec L\longrightarrow \Spec \kappa(\eta)\hookrightarrow X$$
  for some scheme point $\eta$\,. Let $\nu_{\eta}$\, be the restriction of $\nu$\, to $\kappa(\eta)$\,. There is in any case  a valuation $\omega$\, of $K(X)$ dominating $\mathcal O_{X,\eta}$\, and having residue field $\kappa(\omega)$\, with $\kappa(\eta)\subset \kappa(\omega)$\,.\\
 By the extension theorem for valuations (see \cite{ZaSam}[Volume II, chapter VI, Theorem 11, p.26]), we can extend the valuation $\nu_{\eta}$\, of $\kappa(\eta)$\, to a valuation $\nu_{\omega}$\, of $\kappa(\omega)$\,.\\
    We form the composite valuation $\omega\circ \nu_{\omega}$\, (see \cite{Val}[XXX1.2, Proposition 1.8, Proposition 1.11, Proposition 1.12]), which is again a valuation of $K(X)$\,.\\
   According to the bijectivity or the injectivity of the map $q$, it has a unique (maximal one) center on $X$, that is, we get a morphism 
 $$\Spec A_{\omega\circ \nu_{\omega} }\longrightarrow\Spec \mathcal O_{X,\eta}\hookrightarrow X$$ corresponding to the point $$\omega\circ \nu_{\omega}\in R(K(X)).$$
  According to \cite{Val}[Proposition 1.11], $\nu_{\omega}$ corresponds to a prime ideal $\mathfrak{p}$\, in $A_{\omega\circ \nu_{\omega}}$\, and the factor ring $A_{\omega\circ \nu_{\omega}}/\mathfrak{p}$\, is isomorphic to the valuation ring $A_{\nu_{\omega}}$\, of $\kappa(\omega)$\,. Restricting to residues we get an extension 
  $$\Spec A_{\nu_{\omega}}\longrightarrow X$$
  of the morphism
  $$ \Spec \kappa(\omega)\longrightarrow \Spec \kappa(\eta)\longrightarrow X.$$
  It factors over a morphism $$\Spec A_{\nu_{\eta}}\longrightarrow X.$$ 
  Composing with $\Spec A_{\nu}\longrightarrow \Spec A_{\nu_{\eta}}$\, we get the desired extension.\\
In order to show uniqueness of the constructed map $$\Spec A_{\nu}\longrightarrow \overline{\Spec \kappa(\eta)}\hookrightarrow X,$$
 if there where two such maps, then by composing some extension $\nu_{\omega}$\, of $\nu_{\eta}$\, with the same valuation $\omega$, we would get two different maps $$\Spec A_{\omega\circ \nu_{\omega}}\longrightarrow \mathcal \Spec\, \mathcal O_{X,\eta},$$ because the center of a composed valuation $\omega\circ \nu_{\omega}$\, is equal to the center of\linebreak $\nu_{\omega}\in R(\kappa(\omega))$\,. This is mapped to the center of $$\Spec A_{\nu}\longrightarrow \overline{\Spec \kappa(\eta)}\longrightarrow X$$ on the closure of the image of $\Spec \kappa(\nu)$\, in $X$.\\ This  contradicts the uniqueness of valuations with quotient field $K(X)$\,.\\
The assertion that $X$ is complete or separated then follows from the usual valuative criterion  of properness and separatedness, respectively.\\
 That $X$ is of finite type if it is complete and locally of finite type follows from the quasicompactness of the Riemann variety $R(X)\cong R(K(X))$\,: If $X=\bigcup_{i\in I}\Spec A_i$\,, then $$R(X)=\bigcup_{i\in I}R\Spec A_i,$$
  hence by quasicompactness there are finitely many indices $i=1,...,n$\, such that $$R(X)=\bigcup_{i=1}^nR(\Spec A_i).$$
   But this is only possible if $X=\bigcup_{i=1}^n\Spec A_i$\, since otherwise there were a point $$P\in X\backslash (\bigcup_{i=1}^n\Spec A_i)$$
    and a valuation centered over $P$\,(which always exists) and which would not be in $\bigcup_{i=1}^nR(\Spec A_i)$\,.\\  Hence $X$ possesses a finite covering by spectra of finitely generated $k$-algebras and is thus of finite type.
\end{proof}
\begin{remark}
\mylabel{rem:R103} Of course, if one is given a prescheme locally of finite type, it does not need to be of finite type, one may glue infinitely many copies of the same variety along any open subvariety. One might wonder whether a separated scheme locally of finite type has to be of finite type. That this is not the case is illustrated by the following example, see \cite{Hoffmann}.\\
 Let $X=\mathbb P^2$, and consider a line $\mathbb P^1\subset \mathbb P^2$\,. For each closed point $P\in \mathbb P^1,$\, consider the blow-up $X_P:=\rm{Bl}_P\mathbb P^2$\,. Let $L_P$ be the strict transform of $\mathbb P^1$\, and $E_P$\, the exceptional divisor on $X_P$\,. We let $Y_P:=X_P\backslash L_P$\,. The open subsets $U_P:=Y_P\backslash (E_P\cap Y_P)$\, are all isomorphic to $\mathbb P^2\backslash \mathbb P^1$\,, hence we can glue the varieties $Y_P$\, along $U_P$\, to a prescheme $X$ together with a morphism $X\longrightarrow \mathbb P^2$.\\
  If the base field is infinite, this is obviously not of finite type but only locally of finite type because the $Y_p$ are varieties.\\
   But $X$ is separated, for, if there were two morphisms $$p_1,p_2:\Spec A_{\nu}\longrightarrow X,$$
   for some valuation ring $A_{\nu}\subset L$\, with image the same generic point $\Spec L\longrightarrow X$\,, the images of the closed points in $\mathbb P^2$\, must be the same as $\mathbb P^2$\, is separated.\\
     If the image does not lie on $\mathbb P^1$\,, then $p_1=p_2$\,; if the image is a point $P\in \mathbb P^1$\,, then $p_1$ and $p_2$\, must have image in $Y_P$\, and this is separated, hence again $p_1=p_2$\,. That $X$ is separated follows from \ref{lem:L2062}.\\
      The above lemma also shows that $X$ is not compactifyable to a complete scheme $X'$ locally of finite type because $X'$ and hence $X$ would be of finite type.
\end{remark}
 In order to prove that $X$\, satisfies the existence condition of the valuative criterion of properness, it does not suffice to consider valuations with quotient field $K(X)$\,. Consider the following example.
 \begin{example}\mylabel{ex:E1}
  Let $X$\, be the prescheme obtained from $\mathbb P^2$\, by doubling one line $\mathbb P^1$\, and then by removing one closed point $P$\, from one copy $L_1$\, of $\mathbb P^1$\,. Then as $X$ contains $\mathbb P^2$\, as an open subscheme, there is a morphism $\Spec A_{\nu}\longrightarrow X$\, for each valuation $\nu$\, with quotient field $K(X)$\,. But if one considers the generic point $\eta$\, of the copy $L_1$ of $\mathbb P^1$\,, where one has removed the point $P$, if $\mu$\, is the discrete valuation of $\kappa(\eta)$\, corresponding to $P$\,, the canonical morphism $\Spec \kappa(\eta)\longrightarrow X$\, has no extension to a morphism $\Spec A_{\mu}\longrightarrow X$\,.
  \end{example}
Nevertheless we prove the following lemma.
\begin{lemma}\mylabel{lem:L2063} An integral prescheme satisfies the existence condition of the valuative criterion of properness if and only if for each valuation $\omega$\, of $K(X)$\,, each morphism $f:\Spec A_{\omega}\longrightarrow X$\,and each composed valuation $\nu=\omega\circ \mu$\, with $\mu\in R(\kappa(\omega))$\,, there is a morphism $f':\Spec A_{\nu}\longrightarrow X$\, extending $f$ in the sense that $f=f'\circ j$\,, where $j:\Spec A_{\omega}\longrightarrow \Spec A_{\nu}$\, is the natural inclusion.
\end{lemma}
\begin{remark} 
\mylabel{rem:R104}Observe that the assumptions include the case where $\omega$\, is the trivial valuation, where the condition comes down to the fact that every valuation of $R(K(X))$\, has at least one center on $X$\,.
\end{remark}
\begin{proof} First assume that the conditions of the lemma are satisfied. Let \linebreak $\Spec L\longrightarrow X$\, be a morphism factoring over $\Spec \kappa(\eta)$\, for some scheme point $\eta$\, of $X$\, and let $\nu$\, be a valuation of $L$\,. The valuation $\nu$\, induces by restriction a valuation $\nu_{\eta}$\, of $\kappa(\eta)$\,. In any case, there is a valuation $\mu$\, of $K(X)$\, and a morphism $f:\Spec A_{\mu}\longrightarrow X$\, with center $c_X(\mu)=\eta$\, on $X$\,. It follows that there is an inclusion of residue fields $\kappa(\eta)\subset \kappa(\mu)$\,.\\ By \cite{Val}[Proposition 1.14], if $K_1\subset K_2$\, is any field extension, and $\nu$\, is a valuation of $K_1$\,, there is always an extension of $\nu$\, to $K_2$\,. In our case, we can find an extension of the valuation $\nu_{\eta}$\, of $\kappa(\eta)$\, to a valuation $\nu_{\mu}$\, of $\kappa(\mu)$\,. \\
We consider the composed valuation $\omega:=\mu\circ \nu_{\mu}$\,. By assumption, there is a morphism $f':\Spec A_{\omega}\longrightarrow X$\, extending $f$. There is a prime ideal $\mathfrak{p}$\, in $A_{\omega}$\, such that $$A_{\omega}/\mathfrak{p}\cong A_{\nu_{\mu}}$$
 and the composed morphism $$\Spec A_{\nu_{\mu}}\longrightarrow \Spec A_{\omega}\longrightarrow X$$ is an extension of the morphism
  $$\Spec \kappa(\mu)\longrightarrow \Spec \kappa(\eta)\longrightarrow X.$$
 As the valuation $\nu_{\mu}$\, of $\kappa(\mu)$\, induces the valuation $\nu_{\eta}$\, of $\kappa(\eta)$\,, we get an extension of the above morphism 
 $$\Spec A_{\nu_{\mu}}\longrightarrow \Spec A_{\nu_{\eta}}\longrightarrow  X.$$
  As the valuation $\nu$\, of $L$\, induces the valuation $\nu_{\eta}$\, of $\kappa(\eta)$\, we get the desired morphism $$\Spec A_{\nu}\longrightarrow \Spec A_{\nu_{\eta}}\longrightarrow X$$ extending 
  $$\Spec L\longrightarrow \Spec\kappa(\eta)\longrightarrow X.$$\\
Conversely, if $X$ satisfies the existence condition, let $\omega\in R(K(X))$ and\linebreak  $f:\Spec A_{\omega}\longrightarrow X$\, and let $\eta\in X$ be the image of the closed point of $\Spec A_{\omega}$\,. If \linebreak $\mu\in R(\kappa(\omega))$\, is given, by the existence condition there is a morphism $g:\Spec A_{\mu}\longrightarrow  X$\, extending the morphism induced by $f$, $\Spec \kappa(\omega)\longrightarrow X$\,. Now the situation is as follows. Let $x\in X$\, be the image of the closed point of $A_{\mu}$\,. Then there is a prime ideal $\mathfrak{p}$\, in $B:=\mathcal O_{X,x}$\, corresponding to $V:=\overline{\{\eta\}}$\, and we have a local homomorphism $i:B_{\mathfrak{p}}\subset A_{\omega}$\, with induced inclusion of residue fields
 $$r(i):B_{\mathfrak{p}}/\mathfrak{p}B_{\mathfrak{p}}\hookrightarrow \kappa(\omega)$$ as well as an inclusion $j:B/\mathfrak{p}B\subset A_{\mu}$\,. $B/\mathfrak{p}B$\, is a subring of $B_{\mathfrak{p}}/\mathfrak{p}B_{\mathfrak{p}}$\,\,.\\
 The inverse image of $A_{\mu}$\, in $A_{\omega}$\, is by definition  the valuation ring of the composed valuation  $A_{\omega\circ\mu}$\,; the inverse image of $B/\mathfrak{p}B$\, in $B_{\mathfrak{ p}}$\,via the residue homomorphism is by a well known fact of commutative algebra precisely $B$\,. This gives the  inclusion $B\subset A_{\omega\circ\mu}$\, which in turn gives rise to the desired extension $$f':\Spec A_{\omega\circ \mu}\longrightarrow X$$ of $f$.
\end{proof}  
\begin{lemma}\mylabel{lem:L2064} Let $X$ be an integral prescheme  of finite type such that the morphism $q_X: R(X)\longrightarrow R(K(X))$\, induces an isomorphism $R^1(X)\cong R^1(K(X))$\,. Then $q_X$ is an isomorphism and $X$ is complete.
\end{lemma}
\begin{proof} We first show that $q_X$ is injective. Suppose that there are two points $\nu_1,\nu_2\in R(X)$\, of dimension smaller than $\dim X -1$\, having the same image under $q_X$. They cannot have the same image under the morphism $p_X: R(X)\longrightarrow X$\,, since the image would lie in some affine $\Spec A$\, and $q_X$ restricted to $p_X^{-1}(\Spec A)=R(\Spec A)$\, is an isomorphism. Say,
$$p_X(\nu_1)=x_1\in \Spec A_1\,\, \mbox{and}\, \,p_X(\nu_2)=x_2\in \Spec A_2.$$
 We now consider $R(\Spec A_1)$\, and $R(\Spec A_2)$\, naturally as subspaces of $R(K(X))$\, and in it the  subsets  $$V(\{x_1\}):=p_X^{-1}(\{x_1\})\subset R\Spec A_1$$ and $$V(\{x_2\}):=p_X^{-1}(\{x_2\})\subset R\Spec A_2.$$  These are the sets of valuations  lying above  of $x_1$ and $x_2$\,, respectively.\\
 We know that the intersection $V(\{x_1\})\cap V(\{x_2\})$\, inside $R(K(X)/k)$\, is nonempty; it contains the point $q_X(\nu_1)=q_X(\nu_2)$\,. By \ref{lem:L104}, we can choose a complete model $Y$ of $K(X)$\,, open subsets $U_1,U_2\subset Y$\, and proper surjective morphisms   $p_1: U_1\longrightarrow \Spec A_1$\, and $p_2: U_2\longrightarrow \Spec A_2$\,.\\
  As $$V(\{x_1\})\cap V(\{x_2\})\neq \emptyset$$ and 
  $$ (V(\{x_1\})\cap V(\{x_2\})_Y\subseteq V(\{x_1\})_Y\cap V(\{x_2\})_Y$$ we must have 
  $$V(\{x_1\})_Y\cap V(\{x_2\})_Y\neq \emptyset.$$  Let $y\in Y$\, be any scheme point in this intersection such as $c_Y(q_X(\nu_1))$\,.  The discrete algebraic rank one valuation $\omega$\, corresponding to some irreducible component of the preimage of  the exceptional divisor on the blow-up $\rm{Bl}_{\overline{\{y\}}}Y$\, on the normalization $(\mbox{Bl}_{\overline{\{y\}}}Y)^n$\, has center 
  $$y\in V(\{x_1\})_Y\cap V(\{x_2\})_Y.$$ 
   We can also consider $V(\{x_1\})$\, and $V(\{x_2\})$\, as subsets of $R(X)$\, since we have natural inclusions  $$V(\{x_i\})\subset R(\Spec A_i)\subset R(X), i=1,2.$$
    Inside $R(X)$\, their intersection is empty, since the points $x_1$ and $x_2$ are distinct. Over both points lie two different discrete algebraic rank one valuations $\omega_1$\, and $\omega_2$\,, $$\omega_i:=(q_X\mid_{R\Spec A_i})^{-1}(\omega)\in V(\{x_i\}), i=1,2$$
    (observe that $q_X\mid_{R\Spec A_i}, i=1,2$\, is an isomorphism).  By construction, they are both mapped to $\omega$\, under the morphism $q_X$.\\
   But $q_X$ was assumed to be an isomorphism on discrete rank one algebraic valuations, a contradiction.  Hence the assumption was false and  $q$ is injective.\\
    If $q_X$ were not surjective, then the image under $q_X$\, of $R(X)$ would be a proper subset of $R(K(X))$\,. By \ref{lem:L2001}  it follows that $X$ is a separated scheme of finite type. The non-surjectivity of $q_X$ then corresponds to the noncompleteness of $X$. By \cite{Nagata}, there is an integral compactification $\widetilde{X}$\, of $X$.  But there is always a discrete algebraic rank one valuation dominating a closed point $P$\,of $\widetilde{X}\backslash X$\,(take again the  valuation corresponding to any irreducible component of the preimage of the exceptional divisor $E$\, in $\mbox{BL}_P\widetilde{X}$\, in $(\mbox{Bl}_P\widetilde{X})^n$\,). This valuation would not be in $R(X)$\,, contradicting again the assumption that $R^1(X)\cong R^1(K(X))$\,. Thus $q_X$ is an isomorphism.
\end{proof}
\begin{remark}
\mylabel{rem:R105} The  finite type assumption is essential since we need that the exceptional divisor of the blowing up of a point gives a discrete algebraic rank one valuation. This is not true in general. But even if one assumes that $X$ is locally of finite type, the assertion is not true. As in \ref{lem:L1000}, consider with the same notation the sequence of blow-ups $X_{n+1}\longrightarrow X_n$\,. Then $U_n=X_n\backslash \{P_n\}$\, is naturally contained as an open subset in $U_{n+1}=X_{n+1}\backslash \{P_{n+1}\}$\,. If one puts $U:=\bigcup_nU_n$\,, then $U$ is locally of finite type, separated but not complete, as the valuations centered above all the $P_n$ are missing in $R(U)$\,. But as we have shown, $R^1(U)\cong R^1(K(U))$\,.The proof does not work in this case because a separated prescheme locally of finite type may not be compactifyable.\\
This lemma can be seen as a proof of the fact that in order to verify the valuative criterion of properness, it suffices to consider discrete algebraic rank one valuations of the function field if one knows that the prescheme in question is integral and of finite type over a base field. 
\end{remark}
\begin{lemma}\mylabel{lem:L2001} Let $X$ be an integral prescheme  of finite type over the base field $k$. Consider the associated Riemann variety $R(X)$\, with natural morphism $q_X:R(X)\longrightarrow R(K(X))$\,. Then $q_X$ is an isomorphism outside a Zariski closed subset of $R(K(X))$\,.  More generally, if $\mathcal U$\, is a Zariski open subset of $R(K(X))$\,, the same holds for the restriction of $q_X$ to $\mathcal U$\,.
\end{lemma}
\begin{proof} If $X= \bigcup_{i=1}^n\Spec A_i$\, is an affine open covering, consider the open subset $$\mathcal U':=\bigcap_{i=1}^n R\Spec A_i\cap \mathcal U\,\, \mbox{inside}\,\, R(X).$$
 This subset is separated as it is a subset of the Riemann variety of some affine open subset, say $R\Spec A_1$\, and these are always separated. Let  $\mathcal V:=q_X(\mathcal U')$\,. The morphism $q_X$\, restricted to $R\Spec A_1$\, is an isomorphism, so  $\mathcal U'$\, maps isomorphically onto $\mathcal V$\, and this is an open subset of $R(K(X))$\,.\\
  Suppose, that some point $\nu\in R(X)\backslash\mathcal U'$\, maps to a point $\omega\in \mathcal V$\,. The valuation  $\nu$ has to lie in some affine $R(\Spec A_i)$\,. As remarked above there is always a second point $\nu'\in \mathcal U'$\, mapping to $\omega$\,, and this point $\nu'$\, lies a fortiori in $R(\Spec A_i)$\,. As $q_X$ is the identity on $R(\Spec A_i))$\,, we must have $\nu=\nu'$\,, a contradiction. Thus $q_X^{-1}(\mathcal V)=\mathcal U'$\, and $q_X$\,  is an isomorphism over $\mathcal V$\,. Clearly, as $\mathcal V$\, is a finite intersection of open subsets inside $\mathcal U$\,, we have $$\mathcal U\backslash\mathcal U'=\bigcup_{i=1}^n(\mathcal U\backslash (\mathcal U\cap R\Spec A_i))$$ and this is a finite union of Zariski closed subsets inside $\mathcal U$\, and so relatively Zariski closed inside $\mathcal U$\,. 
\end{proof}
\begin{definition}
\mylabel{def:D104} Let $X$ be an integral prescheme and let $$q_X:R(X)\longrightarrow R(K(X))\,\,,p_X:R(X)\longrightarrow X$$ be the natural morphisms. Let $n\geq 1$\, be an integer. We define the following subsets:
$$S_n=S_n(X):=\{\nu\in R(K(X))\mid \sharp \mid q_X^{-1}(\{\nu\})\mid \geq n\},$$
$$T_n=T_n(X):=p_X(q_X^{-1}(S_n(X))).$$
\end{definition}
\begin{remark}
\mylabel{rem:R106} By definition, a point $\nu\in R(K(X))$\, is in $S_n(X)$\, if there exist at least $n$ distinct preimages under the morphism $q_X$\, in $R(X)$\,.\\
A point $x\in X$\, is in $T_n(X)$\, if there exists a valuation $\nu\in R(K(X))$\, and $n$ distinct morphisms $f_i: \Spec A_{\nu}\longrightarrow X, i=1,...,n$\, with (distinct) centers $x_1,...,x_n\in X$\, such that $x_1=x$\,.\\
 If $X$ satisfies the existence condition of the valuative criterion of properness, then  the morphism $q_X$ is surjective and $S_1(X)=R(K(X))$\,. We have $$S_n(X)\supseteq S_{n+1}(X).$$ If $X'$\, is another integral prescheme with the same quotient field $K(X)$\, and $f: X'\longrightarrow X$\, is a proper birational morphism, by \ref{lem:L100}, we have that the natural morphism $R(X')\longrightarrow R(X)$\, is the identity and the morphism $q_{X'}$\, can be identified with the morphism $q_X$\, and it follows that for all $n\in \mathbb N,$\, $S_n(X')=S_n(X)$\,.
\end{remark}
For further use in the sequel, we prove an auxiliary lemma.
 \begin{lemma}
 \mylabel{lem:L103} Let $n\geq 1$\, be an integer and  $\nu\in S_{n}(X)$\, a point with distinct centers $x_1,...,x_n$\, on $X$. Then there is an integral prescheme $X'$ plus a proper birational morphism $f:X'\longrightarrow X$\, such that if $x_{i}'$\, is the unique lift of $x_i$\, for the valuation $\nu$\, via the valuative criterion of properness for the morphism $f$, the residue field extension $\kappa(x_{i}')\subset \kappa (\nu)$\, is algebraic for $i=1,...,n$\,.
\end{lemma}
\begin{proof}
 Let $x_i\in \Spec A_i, i=1,...,n$\,, $\Spec A_i$\, an open affine subset of $X$\, for $i=1,...,n$\,.\\
 We first construct a proper (we may even achieve projective) birational morphism $X_1\longrightarrow X$\, such that if $x_{i1}\in X_1$\, is the unique lift of the center $x_1$\, of $\nu$\, on $X$, then the residue field extension $\kappa(x_{i1})\subset \kappa(\nu)$\, is algebraic.  By \cite{Val}[chapter 2, Proposition 2.3, p.26] this is possible for the separated scheme $\Spec A_1$\,, i.e., we can find an ideal $\mathfrak{a}_1\subset A_1$\, such that on the blow-up $\mbox{Bl}_{\mathfrak{a}_1}(\Spec A_1)$\, the residue field extension is algebraic. We can find on $X$ a sheaf of ideals $\mathcal I_1$\, whose restriction to $\Spec A_1$\, is $\mathfrak{a}_1$\, (we may take $$\mathcal I_1=j_{1*}\mathfrak{a}_1\cap \mathcal O_X, j_1:\Spec A_1\hookrightarrow X$$ being the open immersion). We put $X_1:=\mbox{Bl}_{\mathcal I_1}X$\,. Let 
 $$x_{11}\in\Spec A_{1,1} ,...,x_{n,1}\in \Spec A_{1,n} $$ be the  $n$ distinct centers of $\nu_i$\, on $X_1$\, which are the lifts of $x_1,...,x_n\in X$\,, respectively. Then $\kappa(x_{11})\subset \kappa(\nu_i)$\, is algebraic.\\
Suppose we have constructed the integral prescheme $X_k, (k<n)$\, plus a proper birational morphism $X_k\longrightarrow X$\, such that, if 
$$x_{1k}\in \Spec A_{k,1},...,x_{n,k}\in \Spec A_{k,n}$$ are the distinct centers of $\nu_i$\, on $X_k$ which are the lifts of $x_1,...,x_{n}$\,,respectively, we have that $$\kappa(x_{1k})\subset \kappa(\nu_i),...,\kappa(x_{kk})\subset \kappa(\nu_i)$$
are all algebraic. Let Let $x_{k,k+1}\in \Spec A_{k,k+1}\subset X_k$\,. Then we can find an ideal $\mathfrak{a}_{k,k+1}\subset A_{k,k+1}$\, such that on the blow-up the residue field extension is algebraic. We find an ideal sheaf $\mathcal I_{k}$\, on $X_k$\, that extends $\mathfrak{a}_{k,k+1}$\, on $\Spec A_{k,k+1}$\,. We put $X_{k+1}:=\mbox{Bl}_{\mathcal I_{k}}X_k$\,. Let $x_{1,k+1},...,x_{n,k+1}$\, be the centers of $\nu_i$\, on $X_{k+1}$\, which are the lifts of $x_{1k},...,x_{n,k}$\, on $X_k$\,, respectively. By construction, the residue field extension $\kappa(x_{k+1,k+1})\subset \kappa(\nu_i)$\, is algebraic. As $\kappa(x_{jk})\subset \kappa(x_{j,k+1})$\, and $\kappa(x_{jk})\subset \kappa (\nu_i)$\, is algebraic for $1\leq j\leq k$\,, so is $\kappa(x_{j,k+1})\subset \kappa(\nu_i), 1\leq j\leq k$\,.\\
 We finally put $X':=X_{n}$\,. 
\end{proof} 
\begin{lemma}\mylabel{lem:L2066} Let $X$ be as above an integral prescheme and $q_X:R(X)\longrightarrow R(K(X))$\,  the natural morphism. Suppose that $X$ satisfies the existence condition of the valuative criterion of properness. Then the sets $S_n(X), n\in \mathbb N$\, are closed under specialization.
\end{lemma}
\begin{proof} Recall that a valuation $\omega$\, is a specialization of $\nu$\,iff $\omega\in \overline{\{\nu\}}$\,, where the closure is with respect to the Zariski topology on $R(K(X))$\,, iff $\omega$ is a composition of $\nu$\, with some valuation $\nu_1$\, on the residue field of $\nu$\, (see \cite{Val}[Theorem 2.6]). \\
Suppose to the contrary that there are $n$ distinct morphisms $$f_1,f_2,...,f_n:\Spec A_{\nu}\longrightarrow X$$ lifting the morphism $\Spec K(X)\longrightarrow X$\, but there are at most $n-1$ morphisms $\Spec A_{\omega}\longrightarrow X$\, for a specialization $\omega$\, of $\nu$\,. Let $\omega=\nu\circ \nu_1$\, with $\nu_1$\, a valuation of the residue field of $\nu$\,.\\
  By \ref{lem:L103}, we can always arrange by blowing up $X$ to $p:X'\longrightarrow X$\, that under the (unique) lifts $f_1',f_2',...,f_n'$\, of $f_1,...,f_n$\, with centers $x_1',x_2',...,x_n'$\,, respectively to $X'$\,, the inclusions of residue fields $$\kappa (x_1')\subset \kappa(\nu),\,...,\,\,\kappa(x_n')\subset \kappa(\nu)$$ are algebraic, since $\kappa(\nu)$\, has finite transcendence degree.\\
 Now by the existence condition of the valuative criterion of properness we can lift the morphisms $$f_i'\mid_{\Spec \kappa(\nu)}:\Spec \kappa(\nu)\longrightarrow X', i=1,...,n$$
  to morphisms $$g_i:\Spec A_{\nu_1}\longrightarrow X', i=1,...,n.$$
  Let the images of the closed points in $X'$\, be 
$$y'_1,y_2',...,y'_n\in X'.$$
  Then we can compose the valuations $f_1', f_2',...,f_n'$ with $g_1,g_2,...,g_n$, respectively to get valuations
  $$h_1,h_2,...,h_n:\Spec A_{\omega}\longrightarrow X'$$
  (see \cite{Val}[chapter 1.2, Proposition 1.8, Proposition 1.11, Proposition 1.12]). \\
  By assumption, there is $i,j$ such that $h_i=h_j$\,. The valuation $\nu$\, corresponds to a prime ideal $\mathfrak{p}_{\nu}$\, in $A_{\omega}$\, and one has $A_{\nu}=(A_{\omega})_{\mathfrak{p}_{\nu}}$\,. As $h_i=h_j$\, is both a specialization of $f_i$ and $f_j$, respectively, the morphisms $f_i$ and $f_j$ factor both through the localization morphism $\Spec A_{\nu}\longrightarrow \Spec A_{\omega}$\, followed by $h_i=h_j$\,. Hence $f_i'=f_j'$\,, and composing with $p:X'\longrightarrow X$ we get $f_i=f_j$\,, a contradiction.
\end{proof}
\begin{theorem}\mylabel{thm:A311} Let $X$ be an integral prescheme of finite type and let \linebreak $q_X:R(X)\longrightarrow R(K(X))$\, be the natural morphism. Let $n\geq 2$\, be an integer. We consider the Zariski closure $\overline{S_n(X)}$\, of the subset $S_n(X)\subset R(K(X))$\, above which each point has at least $n$ preimages under $q_X$.\\
Let $Z$ be any complete normal model of $K(X)$\, and $Z_1$ be an irreducible component of the trace $\overline{S_n(X)}_Z$\,. \\
 Let $\nu\in \overline{S_n(X)}$\, be a point of maximal dimension in  $\overline{S_n(X)}\times_ZZ_1$\,, (i.e., $\mbox{trdeg}_k\kappa(\nu)$\, is assumed to be maximal among all points of $\overline{S_n(X)}\times_ZZ_1 $\,).\\
  Then $\nu\in S_n(X)$\,.
\end{theorem}
Before we are  proving  the theorem, in order to motivate it, we give its main application, namely the deduction of the fact that the subsets $$S_n(X)\subset R(K(X)/k)$$ are Zariski closed in case $X$ satisfies the existence condition of the valuative criterion of properness. The proof of \ref{thm:A311} is difficult and contains the major amount of work necessary for the proof of this easily formulated fact. All our efforts to show  directly that the sets $T_n(X)\subset X$\, are Zariski closed, e.g. by using \ref{lem:L2001}, the topological Nagata criterion for openness for the loci $X\backslash T_n(X)$\, (see \cite{Matsumura}[Chapter 8, Paragraph 24, Theorem 24.2, p.186]) and induction on the dimension of $X$, did not succeed. 
Guided by \cite{Ha}[Chapter II, 4, Lemma 4.5] saying that the image of a quasicompact morphism of schemes is closed if it is stable under specialization we are now proving the following
\begin{theorem}\mylabel{thm:T2068} Let $X$ be an integral prescheme  of finite type that satisfies the existence condition of the valuative criterion of properness. Let  $q_X: R(X)\longrightarrow R(K(X))$\, be as above the natural morphism. Then the locus $S_n(X)$\, in \,$R(K(X))$\,  is Zariski closed for all $n\in \mathbb N$\,.
\end{theorem}
\begin{proof} 
We may assume that $n\geq 2$\, since $S_1(X)=R(K(X)/k)$\,( as the morphism $q_X$\, is surjective) and this is trivially Zariski closed.\\
 Let us fix a point\,$\nu \in \overline{S_n(X)}$\,. We can find a finite open affine covering $X=\bigcup_{i=1}^m\Spec A_i$\,, compactify each of the $\Spec A_i$\, to a complete variety $Y_i$\, and  find by \ref{lem:L104}  a complete normal model $Y$\, dominating all the $Y_i$\,, open subsets $U_i\subset Y, $\, and proper surjective morphisms $p_i:U\longrightarrow \Spec A_i,i=1,...,m$\,.\\ We consider the trace $\overline{S_n(X)}_Y$\, of the Zariski closure $\overline{S_n(X)}$\, of $S_n(X)$ in $R(K(X))$\, on $Y$. Let $\eta$\, be the generic point of the irreducible component $Y_1$\, of $\overline{S_n(X)}_Y$\, containing the center $c_Y(\nu)=:y$ of $\nu$\, on $Y$.\\
  We choose $\omega\in \overline{S_n(X)}\times_Y\{Y_1\}$\, of maximal dimension. If $y'=c_Y(\omega)$\,, by \cite{Val}[Proposition 2.3] we may assume by further blowing up $Y$ that the residue field extension $\kappa(y')\subset \kappa(\omega)$\, is algebraic.\\
   By the previous proposition we have $\omega\in S_n(X)$\,. I claim that we have  $c_Y(\omega)=\eta$\,. For, the dimension of $\omega$\, is equal to the transcendence degree of $\kappa(y')/k$\, which is less than the transcendence degree of $\kappa(\eta)/k$\, if $y'\neq \eta$\,. There is always a valuation in $\overline{S_n(X)}$\, with center $\eta$\, on $Y$\, which then has dimension strictly larger than $\omega$\,, a contradiction. By \ref{cor:C88}, we have $\overline{\{\omega\}}_Y=Y_1$\,. Thus there is a valuation $\mu$\, composed with $\omega$\,, $\mu=\omega\circ \xi, \xi\in R(\kappa(\omega)/k)$\, with the same center $y$ on $Y$ as $\nu$\,.\\
    By the specialization property (see \ref{lem:L2066}), we have $\mu\in S_n(X)$\,, i.e., $\mu$\, has at least $n$ centers $x_1,x_2,...,x_n$ on $X$. As each $x_i$\, has to lie in some $$\Spec A_{j_i}, i=1,...,n, j_i\in \{1,2,...,m\}$$ and we have proper surjective morphisms $p_j:U_j\longrightarrow \Spec A_j$\, for all $j$, the unique lifts $y_i$\,of the centers $x_i$\, on $U_{j_i}$\, have by the  completeness of $Y$ all to coincide and to be equal to $y\in Y$\,. Hence $\mathcal O_{Y,y}$\, must  dominate all the local rings $\mathcal O_{X,x_1},...,\mathcal O_{X,x_n}$\,.\\
  As $A_{\nu}$\, dominates $\mathcal O_{Y,y}$\, $\nu$\, has also the centers $x_1,x_2,...,x_n$\, on $X$. Hence $\nu\in S_n(X)$\,. Since this holds for any $\nu\in \overline{S_n(X)}$\,, we must have $S_n(X)=\overline{S_n(X)}$\, and $S_n(X)$ is Zariski closed in $R(K(X))$\,. \\
  \end{proof}
We now give the proof of \ref{thm:A311}. \\
\\
\begin{proof}   
By \ref{prop:P888}, $\overline{S_n(X)}$\, is a proper Zariski closed subset of $R(K(X))$\, for all $n\geq 2$\,. By \ref{prop:P888}, the trace $(\overline{S_n(X)})_Y$\, is a proper Zariski closed subset of $Y$ for all models $Y$ of $K(X)$\, and all $n\geq 2$\,. We fix a normal complete model $Z$ of $K/k$\,, and an irreducible component $Z_1$\, of $\overline{S_n(X)}_Z$\, with generic point $z_1$\, and a point\,$\nu \in \overline{S_n(X)}\times_ZZ_1$\, such that the residue field of the corresponding valuation ring has maximal transcendence degree over $k$ among all points of $\overline{S_n(X)}\times_ZZ_1$. This is possible since the transcendence degree is bounded above by $\mbox{trdeg}(K(X)/k)$\,.\\
  Suppose that $\nu\notin S_n(X)$\,. Let $x_1,...,x_{n-1}$\, be the at most $n-1$\, centers of $\nu$\, on $X$.\, By \ref{lem:L103} there is an integral prescheme $X'$\, plus a proper birational morphism $f: X'\longrightarrow X$\, such that, if $x_1',...,x_{n-1}'$ are the lifts of $x_1,...,x_{n-1}$\,, respectively, that are the  centers of $\nu$\, on $X'$, then the residue field extensions $\kappa(x_j')\subset \kappa(\nu), 1\leq j\leq n-1$\, are all algebraic.\\  
 As $X$ is assumed to be of finite type, so is $X'$\, and we have a finite open covering $X'=\bigcup_{i=1}^{m'} \Spec A_i'$\,. We compactify each $\Spec A_i'$\, to a complete variety $Y_i'$\, and then choose a complete variety $Y$\, dominating all of these finitely many $Y_i'$\,, see \ref{lem:L104}. \\
 Then for each $i=1,...,m'$, we have an open subset $U_i\subset Y$\, and proper surjective morphisms $p_i:U_i\longrightarrow \Spec A_i$\,. In particular, if a  valuation has center $x'$ on some $\Spec A_i'$\,, then it also has a center on $U_i\subset Y$\, lying above $x$. \\
 Note that we do not claim that there is a morphism from $Y$ to $X$.\\
 By \cite{Val}[Proposition 2.3] there is a complete model $Y'\longrightarrow Y$\, such that if $y'$\, is the center of $\nu$\, on $Y'$\,, then the residue field inclusion $\kappa(y')\subset \kappa(\nu)$\, is algebraic. We replace $Y$ with a complete normal model lying over $Y'$ and $Z$. Let $\eta$\, be the center of $\nu$\, on $Y$.\\
  From now on we work with complete models $Y'$ above this fixed $Y$ and consider the traces $\overline{S_n(X)}_{Y'}$\,. The point $\eta$\, must  of course be contained in $\overline{S_n(X)}_{Y}$\, and $\eta$\, must be the generic point of the irreducible component $W$ of $\overline{S_n(X)}_{Y}$\, that maps to $Z_1$\, so that $\eta$\, is mapped to $z_1$.\\
   Otherwise the local ring at the generic point of $W$ would be dominated by a valuation ring $A_{\nu'}$\, with corresponding valuation  $\nu'\in\overline{S_n(X)}\times_ZZ_1$\, and then $$\rm{trdeg}_k\kappa(\nu)=\rm{trdeg}_k\kappa(\eta)<\rm{trdeg}_k\kappa(\eta_W)\leq \rm{trdeg}_k\kappa(\nu'),$$
   a contradiction. For the same reason, if $Y'\longrightarrow Y$\, is any proper birational morphism, the induced morphism $$\overline{S_n(X)}_{Y'}\longrightarrow \overline{S_n(X)}_{Y}$$ must be over $W$ generically finite. Otherwise over $\eta$\, there would lie a point with larger transcendence degree which is again dominated by a point of $\overline{S_n(X)}\times_ZZ_1$\, with larger transcendence degree than $\nu$\,.\\
    Therefore the preimage of $\eta$\, in $\overline{S_n(X)}_{Y'}$\, is a finite set and in the projective limit, the set 
    $$I:=\overline{S_n(X)}\times_{Y}\{\eta\}\subset \overline{S_n(X)}\times_Z\{z_1\}\subset \overline{S_n(X)}\times_ZZ_1$$ must be profinite (containing $\nu$) and consists of valuations $\nu_i, i\in I$\, which have all the same transcendence degree as $\nu$\,. \\
    For, $$\rm{trdeg}_k\kappa(\nu)=\rm{trdeg}_k\kappa(\eta)\leq \rm{trdeg}_k\kappa(\nu_i),$$
     but by assumption $$\rm{trdeg}_k\kappa(\nu_i)\leq \rm{trdeg}_k\kappa(\nu).$$ 
      Of course, if the dimension of the residue field of $\nu$\, is $\rm{dim}X-1$\, and $Y'$\, is normal, then there is exactly one (discrete algebraic) valuation in the fibre. \\
 We want to show that one of the so constructed $\nu_i$ is  in $S_n(X)$. If $\eta$\, is a closed point, then as there is an open subset $\eta\in U\subset Y$\, not containing the other irreducible components of $\overline{S_n(X)}_{Y}$\, we have $$\overline{S_n(X)}\times_{Y}\{\eta\}=\overline{S_n(X)}\times_{Y}U$$ is open in $\overline{S_n(X)}$\, and therefore must contain a point of $S_n(X)$\, which must be one of the $\nu_i$\,. So we may henceforth assume that $\dim(\eta)>0$\,.\\
  Suppose all of these $\nu_i$\, are not in $S_n(X)$. Then they all  have at most $n-1$ centers on $X$. Let $$x_{ij}\in \Spec A_{i_j}\subset X, 1\leq j\leq n-1, i\in I, 1\leq i_j\leq m$$ be the centers of $\nu_i$\, on $X$. Let $$x_{ij}'\in \Spec A_{i_j}', 1\leq j\leq n-1, i\in I, 1\leq i_j\leq m'$$ be their lifts to $X'$\,. By assumption, $Y$ dominates a compactification $Y_{i_j}'$\, of $\Spec A_{i_j}'$\,. The center of $\nu_i$\, on $Y$ is by construction $\eta_W$\,. Also, for each morphism 
  $$f_{ij}:\Spec A_{\nu_i}\longrightarrow \Spec A_{i_j}'\subset Y_{i_j}', 1\leq j\leq n-1, i\in I, 1\leq i_j\leq m'$$
   there is a unique lift to a morphism $$g_{ij}:\Spec A_{\nu_i}\longrightarrow Y,\, 1\leq j\leq n-1$$ because the morphism $p_{i_j}:Y\longrightarrow Y_{i_j}'$\,  is proper. As $Y$ is itself  proper, the centers of the lifts $g_{ij}$\, must be $\eta_W$\,  and $$p_{ij}(\eta_W)=x_{ij}', 1\leq j\leq n-1, i\in I.$$
    But then, as $\eta_W$\, dominates each $x_{ij}'$\,, each $\nu_i$\, has center on each $x_{ij}'$\, and as each $\nu_i$\, was assumed to have at most $n-1$ centers on $X$, hence on $X'$\,, the union of the sets $\{x_{ij}',1\leq j\leq n-1\}$\, must have at most $(n-1)$ elements. Hence all the centers of all the $\nu_i$\, are contained in one set $$\{x_1',...,x_{n-1}'\}\subset X'$$
and each $\nu_i$\, has precisely the centers $x_1',x_2',...,x_{n-1}'$\, on $X'$.  Let $x'_j\in \Spec A'_{i_j}, 1\leq j\leq n-1$.\\
      We consider the Zariski closures $$Z'_j:=\overline{\{x'_j\}}, j=1,...,n-1$$ of $x'_j$\, in $X'$\, which as  closed integral subpreschemes of $X'$\, again satisfy the existence condition  because the  morphism $$ Z'_j\hookrightarrow X'\longrightarrow X$$ is proper.\\
       We consider the Riemann varieties $R(Z'_j), j=1,...,n-1$\, with their natural morphisms $$q_{Z'_j}:R(Z'_j)\longrightarrow R(K(Z'_j)).$$
 By \ref{prop:P888}, the loci $$\overline{S_2(Z'_j)}\subset R(K(Z'_j)), j=1,...,n-1$$
 are proper Zariski closed subsets of $R(K(Z'_j))$\,.\\
We need an auxiliary lemma. Recall that the local ring $\mathcal O_{Y,\eta}$\, dominates the (at most) $n-1$ local rings $\mathcal O_{X',x'_1},...,\mathcal O_{X',x'_{n-1}}$\,, and for each $i\in I$\,, $A_{\nu_i}$\, dominates $\mathcal O_{Y,\eta}$\,. We have therefore the natural inclusions of residue fields  $$\phi_j: K(Z_j')\hookrightarrow \kappa(\eta),1\leq j\leq n-1\,\, \mbox{and}\,\, \psi_i:\kappa(\eta)\hookrightarrow \kappa(\nu_i),\,\, i\in I.$$ 
\begin{lemma}
\mylabel{lem:L105} With notations as above, we have
   $$S_n(X)\cap \overline{\{\nu_i\}}\subset \nu_i\circ\left\langle\bigcup_{j=1}^nR(\psi_i\circ \phi_j)^{-1}(S_2(Z_j'))\right\rangle ,\forall i\in I $$
as subsets of $R(K(X))$\,.
\end{lemma}
\begin{proof}  By the Dirichlet box principle.
The sets $S_n(X)\cap \overline{\{\nu_i\}}$\,, are the sets of valuations composed with $\nu_i$\, that have at least $n$ centers on $X$, (and so on $X'$\,). If $\mu_i=\nu_i\circ \omega_i$\, where $\omega_i$\, is a valuation of $\kappa(\nu_i)$\, and we have $n$ different morphisms $$f_1,f_2,...,f_n:\Spec A_{\mu_i}\longrightarrow X'$$ corresponding to $n$ different centers $z_1',...,z_n'\in X'$\,, by the Dirichlet box principle there are $i,j\in \{1,...,n\}$\, such that
  the induced restrictions of $f_i,f_j$\, to $$\Spec A_{\nu_i}\subset \Spec A_{\mu_i}$$  coincide. This is because $\nu_i$\, was assumed to have at most $(n-1)$ centers on $X'$\,, namely $x'_1,...,x'_{n-1}$\,. We have the following  inclusions of residue fields $$\kappa(x'_j)=K(Z'_j)\stackrel{\phi_{j}}\hookrightarrow \kappa(\eta)\stackrel{\psi_i}\hookrightarrow  \kappa(\nu_i).$$
  We have an induced diagram of morphisms of Riemann varieties 
  $$R(K(X))\stackrel{\nu_i\circ (-)}\hookleftarrow R(\kappa(\nu_i))\stackrel{R(\psi_i)}\longrightarrow R(\kappa(\eta))\stackrel{R(\phi_j)}\longrightarrow R(\kappa(x_j')).$$
   The $n$ centers  $$z'_1,z'_2,...,z'_n\,\, \mbox{of} \,\,f_1,f_2,...,f_n$$
   have to lie in the Zariski closures of $\{x'_1\}, \{x'_2\},...,\{x'_{n-1}\}$\,, respectively, hence in $Z'_1,...,Z'_{n-1}$\,.  It follows that there is $i\in \{1,2,...,n-1\}$\, and $j\neq i$\, such that $z_i',z_j'\in Z_i'$\,. So the  valuation $\omega_i$\, of $\kappa(\nu_i)$\, induces by restriction a valuation $\omega_i'$ on $K(Z'_i)$\, which has at least $2$ centers, namely $z_i',z_j'$ on $Z_i'$\,.  It follows that $\omega_i'\in S_2(Z_i')$\,. \\
Then, $$\omega_i\in R(\psi_i\circ\phi_j)^{-1}(S_2(Z_j'))\subset R(\kappa(\nu_i))\subset R(K(X)).$$   Under the inclusion $$\nu_i\circ(-): R(\kappa(\nu_i))\hookrightarrow R(K(X))$$ we get the result.
\end{proof}
We carry on with the proof of the theorem.\\
By \ref{lem:L100}, the maps
   $$R(\psi_i): R(\kappa(\nu_i))\longrightarrow R(\kappa(\eta))\,\,\mbox{and}\,\,R(\phi_j):R(\kappa(\eta))\longrightarrow R(K(Z_j'))$$ are surjective morphisms of locally ringed spaces (see also the extension theorem for  valuations for algebraic field extensions in \cite{Val}[chapter 1, Proposition 1.14, p.11]). Since $\overline{S_2(Z_j')}$\, is a proper closed subset of $R(K(Z_j')), j=1,...,n-1$\,, the sets  $R(\psi_i\circ \phi_j)^{-1}(\overline{S_2(Z_j'))}$\, are proper Zariski closed subsets of  $R(\kappa(\nu_i))\,\,\forall i\in I, 1\leq j\leq n$\,. Therefore, for each fixed $i\in I$\, the finite union $$\bigcup_{j=1}^nR(\psi_i\circ \phi_j)^{-1}(\overline{S_2(Z_j')})$$ is a proper Zariski closed subset of $R(\kappa(\nu_i))$\,. By \ref{lem:L105},$$\forall i\in I:\,\,S_n(X)\cap \overline{\{\nu_i\}}$$ is a fortiori a proper Zariski closed subset of $$R(\kappa(\nu_i))\cong \overline{\{\nu_i\}}.$$   Since the morphism $R(\psi_i)$\, is closed (see \ref{prop:P501}) and surjective (see \ref{lem:L100}), we have 
  $$R(\psi_i)(R(\psi_i\circ \phi_j)^{-1}(\overline{S_2(Z_j')})=R(\phi_j)^{-1}(\overline{S_2(Z_j')}), \forall i\in I, j=1,...,n-1$$ and 
  $$ \bigcup_{j=1}^{n-1}R(\phi_j)^{-1}(\overline{S_2(Z_j')})$$ is a  Zariski closed subset of $R(\kappa(\eta))$\,.
\\
By \ref{prop:P888} it is even a proper Zariski closed subset of $R(\kappa(\eta)),$\, since
\begin{align*}
\dim(R(\kappa(\nu_i)))>\dim(\bigcup_{j=1}^nR(\psi_i\circ \phi_j)^{-1}(\overline{S_2(Z_j')}))\\ \geq \dim(R(\psi_i)(R(\psi_i\circ \phi_j)^{-1}(\overline{S_2(Z_j')}))=\dim(R(\phi_j)^{-1}(\overline{S_2(Z_j')}))
\end{align*} and $\dim(R(\kappa(\nu_i)))=\dim(R(\kappa(\eta)))$\,.
\\ 
  Observe that $W$ is a birational model of $\kappa(\eta)$\,, so we can calculate the traces with respect to $W$ of closed subsets of $R(\kappa(\eta))$\,. We have 
  \begin{align*}(R(\psi_i)(\bigcup_{i\in I}S_n(X)\cap \overline{\{\nu_i\}}))_W\subset (R( \psi_i)(\bigcup_{i\in I}\bigcup_{j=1}^{n-1}R(\psi_i\circ \phi_j)^{-1}(\overline{S_2(Z_j')})))_W\\
  = (\bigcup_{j=1}^{n-1}R(\phi_j)^{-1}(\overline{S_2(Z_j')}))_W
  \end{align*} 
As $$\bigcup_{j=1}^{n-1}R(\phi_j)^{-1}(\overline{S_2(Z_j')})$$ is a proper closed subset of $R(\kappa(\eta)),$\, by \ref{prop:P888}, its trace $$(\bigcup_{j=1}^{n-1}R(\phi_j)^{-1}(\overline{S_2(Z_j')}))_W$$ must be a proper closed subset of $W$. Otherwise, there should be a valuation $\mu$\, above $\eta$\, in this set and the only valuation in $R(\kappa(\eta))=R(\kappa(\eta_W))$\, is the trivial valuation which is the generic point of $R(\kappa(\eta))$\,. Hence the trace 
$$(\bigcup_{i\in I}R(\psi_i)(S_n(X)\cap \overline{\{\nu_i\}}))_W$$
is a proper Zariski closed subset of $W$.\\
The valuation \,$\nu_i$\, cannot be in the closure of $S_n(X)$ if there are no points of $S_n(X)$ with center in $W\backslash (\bigcup_{j=1}^{n-1}R(\phi_j)^{-1}(\overline{S_2(Z_j')}))_W $\,. We conclude that there is a valuation $\mu\in S_n(X)$\, having center  $$w\in W\backslash (\bigcup_{j=1}^{n-1}R(\phi_j)^{-1}(\overline{S_2(Z_j')}))_W.$$\\
  The valuation $\mu$\, has at least $n$ centers $x_1,x_2,...,x_n$\, on $X$\, which are dominated by centers $y_1,y_2,...,y_n$\, of $\mu$\, on $Y$ which by completeness of $Y$ must be equal to  $w\in W$\,. Now, because $(\overline{\{\nu_i\}})_W=W$\,, there is a valuation $\mu'$\, composed with some $\nu_i$\, that has the same center $w$ on $Y$. But as $\mathcal O_{Y,w}$\,  dominates all $\mathcal O_{X,x_1},\mathcal O_{X,x_2},...,\mathcal O_{X,x_n}$\,  and we have a morphism $\Spec A_{\mu'}\longrightarrow \Spec \mathcal O_{Y,w}$\,, $\mu'$\,  also has at least $n$ centers $x_1,x_2,...,x_n$ on $X$, hence belongs to $S_n(X)\cap \overline{\{\nu_i\}}$\,. But then $\mu'$\, must have center in $$(\bigcup_{j=1}^{n-1}R(\phi_j)^{-1}(\overline{S_2(Z_j')}))_W,$$
      a contradiction.\\
      Hence our assumption was false and there is an $i_0\in I$ such that $\nu_{i_0}$\, belongs to $S_n(X)$\,.\\
 In order to show that also $\nu\in S_n(X)$\,, we know that $\nu_{i_0}$\, has at least $n$ centers $x_1,x_2,...,x_n$ on $X$. By what was said above, $\mathcal O_{Y,\eta}$\,  dominates all the local rings  $\mathcal O_{X,x_1},...,\mathcal O_{X,x_n}$\,. But $A_{\nu}$\, dominates $\mathcal O_{Y,\eta}$\,, hence has also $x_1,...,x_n$\, as  $n$ distinct centers on $X$, meaning $\nu\in S_n(X)$\,.
 \end{proof}

\begin{corollary}\mylabel{cor:C13} With notation as in the above theorem, for each integral\\ prescheme $X$ of finite type satisfying the existence condition of the valuative criterion of properness and for each $n\in \mathbb N$\, the sets $T_n(X),n\in \mathbb N$\, are Zariski-closed subsets of $X$.
\end{corollary}
\begin{proof} By definition, $$T_n(X)=p_X(q_X^{-1}(S_n(X))), n\in \mathbb N.$$ As $q_X$\, is continuous, by the above theorem,  $q_X^{-1}(S_n(X))$\, is closed in $R(X)$\,. As by \ref{thm:A310}, $p_X$\, is a closed map, $T_n(X)=p_X(q_X^{-1}(S_n(X)))$\, is Zariski-closed in $X$. 
\end{proof}

Now we further investigate the structure of the sets $T_n(X),n\in \mathbb N$\, for an integral prescheme $X$\, of finite type.
\begin{lemma}
\mylabel{lem:L401}
Let $X$ be an integral prescheme of finite type satisfying the existence condition of the valuative criterion of properness. Then  there is $N\in \mathbb N$\, such that $T_n(X)=\emptyset, \forall n> N$
\end{lemma}
\begin{proof}
 Let $X=\bigcup_{i=1}^N\Spec A_i$\, be an open affine covering and $n>N$\, be given. If $T_n(X)$\, was nonempty, there would be a valuation $\nu\in R(K(X))$\, with $n$ distinct centers on $X$ and by the Dirichlet box principle, two of them would lie in some $\Spec A_j$\,,  contradicting the separability of $\Spec A_j$\,. Thus $T_n(X)=\emptyset\,\,\forall n>N$\,.
\end{proof}
\begin{definition}
\mylabel{def:D111} Let $X$ be an integral prescheme of finite type satisfying the existence condition of the valuative criterion of properness. We put $$X_n:= T_n(X)\backslash T_{n+1}(X), n\in \mathbb N$$
 and call the collection $(X_n)_{n\in \mathbb N}$\, the canonical stratification of the prescheme $X$.
\end{definition}
\begin{remark}
\mylabel{rem:R401} By \ref{cor:C13}, each $T_n(X), n\in \mathbb N$\, is a Zariski-closed subset. Hence $$X_n=T_n(X)\backslash T_{n+1}(X)=T_n(X)\cap (X\backslash T_{n+1}(X))$$ is a constructible set. By \ref{lem:L401}, $T_n(X)=\emptyset$\, for some $N\in \mathbb N, \forall n>N$\,  and  the canonical stratification is finite.\\
In general, we do not have that $X_{n+1}\subset \overline{X_n}$\,. Consider the following example. Let $X$ be the integral prescheme of finite type obtained by triplicating one closed point of $\mathbb P^1_{\mathbb C}$\,. Then $X_2=\emptyset$\, and $X_3$\, consists of the three triplicated points.
\end{remark}

\bibliography{Riemann}
\bibliographystyle{plain}

\noindent
\textsc{Stefan G\"unther}

\noindent
\emph{E-mail address:} \verb!stef.guenther2@vodafone.de!
\end{document}